%% file: main-torsion-v4-lastupdated.tex
\documentclass{amsart}
\input{command}

\begin{document}

\title{An enhanced Euler characteristic of sutured instanton homology}


\author{Zhenkun Li}
\address{Department of Mathematics, Stanford University}
\curraddr{450 Jane Stanford way, Stanford, 94305, USA}
\email{zhenkun@stanford.edu}
\thanks{}

\author{Fan Ye}
\address{Department of Pure Mathematics and Mathematical Statistics, University of Cambridge}
\curraddr{Wilberforce Road, Cambridge, CB3 0WB, UK}
\email{fanye@math.harvard.edu}
\thanks{}

\date{}
\dedicatory{}
\begin{abstract}
For a balanced sutured manifold $(M,\gamma)$, we construct a decomposition of $SHI(M,\gamma)$ with respect to torsions in $H=H_1(M;\mathbb{Z})$, which generalizes the decomposition of $I^\sharp(Y)$ in previous work of the authors. This decomposition can be regarded as a candidate for the counterpart of the torsion spin$^c$ decompositions in $SFH(M,\gamma)$. Based on this decomposition, we define an enhanced Euler characteristic $\chi_{\rm en}(SHI(M,\gamma))\in\mathbb{Z}[H]/\pm H$ and prove that $\chi_{\rm en}(SHI(M,\gamma))=\chi(SFH(M,\gamma))$. This provides a better lower bound on $\dim_\mathbb{C}SHI(M,\gamma)$ than the graded Euler characteristic $\chi_{\rm gr}(SHI(M,\gamma))$. As applications, we prove instanton knot homology detects the unknot in any instanton L-space and show that the conjecture $KHI(Y,K)\cong \widehat{HFK}(Y,K)$ holds for all $(1,1)$-L-space knots and constrained knots in lens spaces, which include all torus knots and many hyperbolic knots in lens spaces.
\end{abstract}
\maketitle
\tableofcontents

\section{Introduction}

Sutured instanton Floer homology was introduced by Kronheimer and Mrowka in \cite{kronheimer2010knots}. It combines instanton Floer homology with sutured manifold theory and has become a powerful tool since then. In \cite{li2019direct}, Ghosh and the first author constructed a decomposition of the sutured instanton Floer homology $SHI(M,\ga)$ of a balanced sutured manifold with respect to the group $(H_2(M,\partial M;\intg))^*\cong H_1(M;\intg)\slash {\rm Tors}.$ More precisely, a basis of $H_2(M,\partial M;\intg)$ induces a multi-grading on $SHI(M,\ga)$ which can be identified with the group $H_1(M;\intg)\slash {\rm Tors}.$ In \cite{LY2021}, the authors of the current paper studied the Euler characteristics of this decomposition of $SHI(M,\ga)$ and related it to the Euler characteristic of $SFH(M,\ga)$, which is known as the sutured Floer homology introduced by Juh\'asz, and whose Euler characteristic has been understood by work of Friedle, Juh\'asz, and Rasmussen in \cite{Friedl2009}. The study of Euler characteristics was further used by the author to compute the instanton Floer homology of some families of $(1,1)$-knots in a general lens space and was recently further utilized by Xie and Zhang \cite{yixiesu2} to prove that links in $S^3$ all admit irreducible $SU(2)$ representations except for connected sums of Hopf links.

However, only having the decomposition of $SHI(M,\ga)$ along the group $H_1(M;\intg)\slash {\rm Tors}$ is not fully satisfactory for the following two reasons.
\benu
\item Among all known Floer homology theories for sutured manifolds, we have known that sutured monopole Floer homology is isomorphic to sutured Floer homology by work of Lekili \cite{lekili2013heegaard} and Baldwin and Sivek \cite{Baldwin2020}, and it is conjectured that the sutured instanton Floer homology is also isomorphic to sutured Floer homology by Kronheimer and Mrowka \cite{kronheimer2010knots}. However, sutured Floer homology decomposes along spin${}^c$ structures and, in particular, the first Chern classes of torsion spin${}^c$ structures have Poincar\'e dual in the torsion part of the group $H_1(M;\intg)$, so there should be some corresponding decomposition of sutured instanton Floer homology. 

\item The original decomposition along $H_1(M;\intg)\slash {\rm Tors}$ collapses all torsion parts into a single summand of $SHI(M,\ga)$, and some information may lost in this collision; see Example \ref{ex: not sharp}.
\eenu

In this paper, in order to solve this problem coming from collapsing torsion parts, we obtain the following.

\bthm[Main theorem]\label{thm: main}
Suppose $(M,\ga)$ is a balanced sutured manifold and $H=H_1(M;\mathbb{Z})$. Then there is a (possibly noncaonical) decomposition\[SHI(M,\ga)=\bigoplus_{h\in H}SHI(M,\ga,h).\]This decomposition depends on some auxiliary choices. In particular, it is defined up to a global shift of $H$. We define the \textbf{enhanced Euler characteristic} of $SHI$ by \[
    \en(SHI(M,\ga))\deq \sum_{h\in H}\chi(SHI(M,\ga,h))\cdot h\in \mathbb{Z}[H]/\pm H.\]Then we have\begin{equation}\label{eq: enhanced same}\en(SHI(M,\ga))=\chi(SFH(M,\ga))\in \mathbb{Z}[H]/\pm H.
\end{equation}The similar results also hold for $SHM(M,\ga)$. 
\ethm
\brem
If $H_1(M;\mathbb{Z})$ has no torsion, then the decomposition in Theorem \ref{thm: main} is just induced by the multi-grading mentioned above and Equation (\ref{eq: enhanced same}) reduces to \cite[Theorem 1.2]{LY2021}. By results in \cite{Friedl2009}, we have $\chi(SFH(M,\ga))=\tau(M,\ga)$, where $\tau(M,\ga)$ is a (Turaev-type) torsion element that can be calculated by Fox calculus. In particular, if $\partial M$ consists of tori and $\ga$ consists of two parallel copies of a curve $m_i$ with opposite orientations on each boundary component, by \cite[Proof of Lemma 6.1]{Friedl2009} and \cite[Proposition 2.1]{Rasmussen2017}, we have \[\tau(M,\ga)=\tau(M)\cdot \prod_{i}([m_i]-1),\]where $\tau(M)$ is the Turaev torsion of $M$ \cite{Turaev2002}.
\erem
Though the decomposition in Theorem \ref{thm: main} has not been proved to be canonical, we expect it to be well-defined up to a global grading shift of $H$. The following theorem indicates this decomposition is a candidate for the counterpart of the spin$^c$ decomposition. The proof is essentially due to \cite{lekili2013heegaard,Baldwin2020}.

\bthm\label{thm:shm=sfh}
Suppose $(M,\ga)$ is a balanced sutured manifold. Then there is a spin$^c$ structure $\mathfrak{s}_0\in\spin(M,\ga)$ such that for any $\mathfrak{s}\in\spin(M,\ga)$, we have an isomorphism\[SHM(M,\ga,PD(\mathfrak{s}-\mathfrak{s}_0))\cong SFH(M,\ga,\mathfrak{s})\otimes \Lambda,\]
where ${\rm PD}: H^2(M,\partial M;\mathbb{Z})\to H_1(M;\mathbb{Z})$ is the Poincar\'{e} duality map and $\Lambda$ is the Novikov ring.
\ethm

For an element in a group ring $\mathbb{Z}[G]$$$x=\sum_{g\in G} c_g\cdot g,\text{ for }c_g\in\mathbb{Z},$$define\[\norm{x}=\sum_{g\in G}|c_g|.\]This definition is still well-defined for an element in $\mathbb{Z}[G]/\pm G$. By construction of Euler characteristics, we have\begin{equation}\label{eq: lower bound}
    \dim_\mathbb{C}SHI(M,\ga)\ge \norm{\en(SHI(M,\ga))}\ge \norm{\gr(SHI(M,\ga))}.
\end{equation}
To provide an example that the second inequality in (\ref{eq: lower bound}) is strict, and hence $\en$ contains more information than $\gr$, we consider an example from constrained knots studied by the second author \cite{Ye2020}.
\begin{exmp}\label{ex: not sharp}
Consider the 1-cusped hyperbolic manifold $M=m006$ in the \textit{Snappy} program \cite{snappy}. We have $H_1(M;\mathbb{Z})\cong \mathbb{Z}\oplus\mathbb{Z}_5\cong\mathbb{Z}\langle t,r\rangle/(5r)$. By the list of constrained knots in \cite{Ye}, Dehn filling along the slope $(1,0)$ (in the basis from \textit{Snappy}) gives the lens space $L(5,3)$ and the core knot is the constrained knot $C(5,3,4,3,1)$. Suppose $\ga$ consists of two parallel copies of the curve of slope $(1,0)$. Then we have\[\tau(M,\ga)=1 + r + t + rt + r^2t - r^3t - r^4t + rt^2 + r^2t^2,
\]and\[\tau(M,\ga)|_{r=1}=1+1+t+t+t-t-t+t^2+t^2=2+t+2t^2.\] Hence we have\[\norm{\en(SHI(M,\ga))}=\norm{\tau(M,\ga)}=9{\rm ~and~} \norm{\gr(SHI(M,\ga))}=\norm{\tau(M,\ga)|_{r=1}}=5.\]
\end{exmp}
Suppose $K$ is a knot in a closed 3-manifold $Y$. Let $$Y(1)\deq Y\backslash B^3\aand Y(K)\deq Y\backslash {\rm int} N(K).$$Suppose $\delta$ is a simple closed curve on $\partial Y(K)\cong S^2$ and suppose $\ga_K$ is two copies of the meridian of $K$ with opposite orientations.  In Heegaard Floer theory, we have\[SFH(Y(1),\delta)\cong \widehat{HF}(Y)\aand SFH(M,\ga_K)\cong \widehat{HFK}(Y,K),\]where $\widehat{HF}(Y)$ and $\widehat{HFK}$ are the hat versions of the Heegaard Floer homology and the knot Floer homology defined by Ozsv\'ath and Szab\'o \cite{ozsvath2004holomorphic,ozsvath2004holomorphicknot,Rasmussen2003}. In instanton theory, Kronheimer and Mrowka \cite{kronheimer2010knots} defined \[I^\sharp(Y)\deq SHI(Y(1),\delta)\aand KHI(Y,K)\deq SHI(Y(K),\ga_K).\]

An application of Theorem \ref{thm: main} is the following unknot detection result in rational homology spheres.
\bthm\label{thm: unknot detection}
Suppose $K$ is a null-homologous knot in a rational homology sphere $Y$. If \begin{equation}\label{eq: Lspace, KHI}\dim_\mathbb{C}I^{\sharp}(Y)=|H_1(Y;\mathbb{Z})|,\end{equation}then $K$ is the unknot \textit{i.e.}, it bounds a disk in $Y$, if and only if
\begin{equation}\label{eq: floer simple knot, KHI}
    \dim_\mathbb{C} KHI(Y,K)=\dim_\mathbb{C}I^{\sharp}(Y).
\end{equation}
\ethm
\brem
Since instanton theory is closely related to $SU(2)$ representations of fundamental groups, Theorem \ref{thm: unknot detection} may be used to show that for any nontrivial null-homologous knot $K$ in a rational homology sphere $Y$, the fundamental group $\pi_1(Y(K))$ admits an irreducible representation in $SU(2)$ such that the meridian of $K$ is mapped to a traceless element in $SU(2)$. However, the authors do not know how to prove the nondegeneracy results similar to \cite[Section 4.1]{baldwin2018stein} for generators of $KHI(Y,K)$.
\erem
\brem
Rational homology spheres that satisfy (\ref{eq: Lspace, KHI}) are called \textbf{instanton L-spaces}. Theorem \ref{thm: unknot detection} cannot be generalized to knots that are not null-homologous because simple knots in lens spaces also satisfy (\ref{eq: floer simple knot, KHI}) \cite[Proposition 1.9]{LY2020}. It is a natural conjecture that simple knots are the only knots in lens spaces satisfy (\ref{eq: floer simple knot, KHI}) (For Heegaard Floer theory, see \cite[Conjecture 1.5]{Baker2007}).
\erem
Following the similar strategy, we can prove the following theorem for knots whose $KHI$ have small dimensions.

\bthm\label{thm: trefoil detection}
Suppose $K$ is a null-homologous knot in a rational homology sphere $Y$. If
\begin{equation}\label{eq: trefoil-like knot, KHI}\dim_{\mathbb{C}}KHI(Y,K)=\dim_{\mathbb{C}}I^{\sharp}(Y)+2=|H_1(Y;\mathbb{Z})|+2,\end{equation}
then $K$ must be a genus-one-fibred knot.
\ethm
\brem
The only knots in $S^3$ satisfying (\ref{eq: trefoil-like knot, KHI}) are the trefoil and its mirror. Hence Theorem \ref{thm: trefoil detection} is a generalization of the trefoil detection result in $S^3$ \cite[Theorem 1.6]{baldwin2018khovanov}. Both proofs are based on the nonvanishing result on the `next-to-top' grading \cite[Theorem 1.7]{baldwin2018khovanov} for fibred knots. 
\erem
\brem\label{rem: +4 remark}
For a knot $K$ in an instanton L-space $Y$ with $$\dim_{\mathbb{C}}KHI(Y,K)=\dim_{\mathbb{C}}I^{\sharp}(Y)+4=|H_1(Y;\mathbb{Z})|+4,$$the previous arXiv version of this paper said we may still conclude $K$ is fibred using the same strategy. However, later we realize that the strategy does not work in this case by Hongjian Yan's comments. On the other hand, even in the case $Y=S^3$, it is impossible to pin down the genus because there are at least two knots with different genera: the figure-8 knot with genus one and the $T_{(2,5)}$ torus knot with genus two. Recently, there are many new results \cite{baldwin21t25,niyi21t25,niyi21slope} about the Khovanov homology and the knot Floer homology of $T_{(2,5)}$. It is an interesting question if these results can be applied to instanton knot homology.
\erem
\brem
We can prove similar detection results in Heegaard Floer theory; see Section \ref{sec: Knots with small instanton knot homology}.
\erem

Another application of Theorem \ref{thm: main} is to compute $KHI(Y,K)$ for all $(1,1)$-L-space knots and constrained knots in lens spaces. The calculation is based on the following theorem.
\bthm[{\cite[Theorem 1.6]{LY2020}, see also \cite[Theorem 1.1]{BLY2020}}]\label{thm: inequality for 1,1 knots}
Suppose $K\subset Y$ is a $(1,1)$-knot in a lens space (including $S^3$). Then we have
$${\rm dim}_{\mathbb{C}}KHI(Y,K)\leq\dim_\ft\widehat{HFK}(Y,K).$$
\ethm
\bcor\label{cor: 11 HFK=KHI}
Suppose $K\subset Y$ is a $(1,1)$-knot in a lens space (including $S^3$). If $K$ is either a L-space knot, or a constrained knot, then $${\rm dim}_{\mathbb{C}}KHI(Y,K)=\dim_\ft\widehat{HFK}(Y,K).$$
\ecor
\bpf
The theorem follows from comparing the upper bound from Theorem \ref{thm: inequality for 1,1 knots} and the lower bound from $\norm{\en(KHI(Y,K))}$. By \cite[Lemma 3.2]{Rasmussen2017} and \cite[Theorem 2.2]{Greene2018} for L-space knots, and by \cite[Section 4]{Ye2020} for constrained knots, the upper bound matches the lower bound.
\epf
\brem
In the authors' previous work \cite[Corollary 1.11]{LY2021}, we proved Corollary \ref{cor: 11 HFK=KHI} in the case where $H_1(Y(K);\mathbb{Z})\cong\mathbb{Z}$. This is because the lower bound from $\norm{\gr(KHI(Y,K))}$ may not equal to the upper bound in \cite{LY2020} when $H_1(Y\backslash {\rm int}N(K);\mathbb{Z})$ has torsions (\textit{c.f.} Example \ref{ex: not sharp}).
\erem
\brem
The result in Corollary \ref{cor: 11 HFK=KHI} can be generalized to other $(1,1)$-knot $K$ whose $\widehat{HFK}$ is totally determined by $\chi(\widehat{HFK}(Y,K))$, such as $(\pm 2,p,q)$ pretzel knots for odd integers $p$ and $q$ (\textit{c.f.} \cite[Section 5]{Goda2005}).
\erem
Any lens space $Y$ has a standard Heegaard splitting of genus $1$. A knot $K$ in a lens space $Y$ is called a \textbf{torus knot} if $K$ can be isotoped to lie on the Heegaard torus of $Y$. This definition generalizes the usual torus knots in $S^3$. A knot $K$ is called a \textbf{satellite knot} if $Y\backslash{\rm int}N(K)$ has an essential torus. A knot $K$ is called \textbf{hyperbolic} if $Y\backslash K$ admits a hyperbolic metric of finite volume. By Thurston’s Hyperbolization Theorem for Haken 3-manifolds, we have a good classification of knots in lens spaces.
\bprop[{\cite[Proposition 3.1]{Williams2009}}]
Suppose $K$ is a knot in a lens space Y. If $Y(K)$ is irreducible, then $K$ is either a torus knot, a satellite knot, or a hyperbolic knot. 
\eprop
It is straightforward to check that torus knots are $(1,1)$-knots, and their complements are Seifert fibred spaces.
\bprop[{\cite[Theorem 5.1]{Rasmussen2017}}]\label{prop: seifert}
Knots in lens spaces with Seifert fibred complements are L-space knots.
\eprop
Combining Proposition \ref{prop: seifert} and Corollary \ref{cor: 11 HFK=KHI}, we have the following result.
\bcor
For any torus knot $K$ in a lens space $Y$, we have
$${\rm dim}_{\mathbb{C}}KHI(Y,K)=\dim_\ft\widehat{HFK}(Y,K).$$
\ecor
Complements of many constrained knots are orientable hyperbolic 1-cusped manifolds with simple ideal triangulations. In particular, among 286 orientable 1-cusped manifolds that have ideal triangulations with at most five ideal tetrahedra, there are 232 manifolds that are complements of constrained knots. More examples can be found in \cite{Ye}. Indeed, \cite[Conjecture 2]{Ye2020} conjectured that most constrained knots are hyperbolic knots.
\subsection{Organization and sketch of the proofs}\quad

Suppose $(M,\ga)$ is a balanced sutured manifold. To sutured instanon and monopole homology together, we use $\shg(M,\ga)$ to denote both $\shi(M,\ga)$ and $\shm(M,\ga)$ (\textit{c.f.} \cite[Section 2.1]{li2019tau}, see also \cite{baldwin2015naturality}). It is a projectively transitive system, where each space in the system is isomorphic to $SHI(M,\ga)$ (for $\shi(M,\ga)$) and $SHM(M,\ga)$ (for $\shm(M,\ga)$). Note that it is different from the formal sutured homology $\shib(M,\ga)$ and $\shmb(M,\ga)$ for instanton and monopole theory constructed in \cite{LY2021} because formal sutured homology corresponds to the untwisted theory in Baldwin and Sivek's construction \cite{baldwin2015naturality}, where $\shg(M,\ga)$ corresponds to the twisted theory in \cite{baldwin2015naturality}. However, we have the equalities for graded Euler characteristics from \cite[Theorem 7.7 and Theorem 9.21]{baldwin2015naturality}$$\gr(\shi(M,\ga))=\gr(\shib(M,\ga))\aand \gr(\shm(M,\ga))=\gr(\shmb(M,\ga)).$$
Hence the results for graded Euler characteristics of formal sutured homology in \cite{LY2021} can be applied to $\shg(M,\ga)$. In particular, we have 
\begin{equation}\label{eq: twisted euler}
    \gr(\shg(M,\ga))=\gr(SFH(M,\ga)).
\end{equation}
\brem
In the previous version of this paper, we used the formal sutured homology to carry out proofs in the paper, but then noticed that some constructions might involve closures of balanced sutured manifolds of different genera. Moreover, the proof of the functoriality of contact gluing maps in \cite{li2018gluing} involves closures obtained from disconnected auxillary surfaces, which can only be handled by a genus one version of Floer's excision theorem that is available in the twisted theory . Hence in the current version, we use the twisted theory $\shg(M,\ga)$, which relates closures of different genera and closures with possibly disconnected auxillary surfaces.
\erem

In Section \ref{Section: Formal sutured homology}, we review basic properties of $\shg(M,\ga)$, the gradings on $\shg(M,\ga)$ associated to admissible surfaces, and the maps on $\shg(-M,-\ga)$ associated to contact handle attachments (called \textbf{contact gluing maps}). Since we will use contact gluing maps frequently, it is more convenient to consider $(-M,-\ga)$, the sutured manifold with the reverse orientation. Hence we state all results with a minus sign.

In Section \ref{Section: Decomposition associated to tangles}, we generalize the decomposition associated to a rationally null-homologous knot in \cite[Section 4]{LY2020} to a connected rationally null-homologous tangle $\al$, \textit{i.e.} $[\al]=0\in H_1(M,\partial M;\mathbb{Q})$. We write $M_{\al}=M\backslash {\rm int }N(\al)$. Then by Lemma \ref{lem: rk}, we have\[{\rm rk}_\mathbb{Z}H_1(M_\al;\mathbb{Z})={\rm rk}_\mathbb{Z}H_1(M;\mathbb{Z})+1,\]and there is a surjective map \[H_1(M_\al;\mathbb{Z})\to H_1(M_\al;\mathbb{Z})/[m_\al]\cong H_1(M;\mathbb{Z}),\]where $m_\al$ is the meridian of $\al$. Moreover, after picking some suitable $\al$, the pre-images of some torsions in $H_1(M;\mathbb{Z})$ are distinguished in the free part of $H_1(M_\al;\mathbb{Z})$. Since the difference in the free part can be detected by the gradings associated to admissible surfaces, we can decompose $\shg(-M,-\ga)$ by considering direct summands of $\shg(-M_\al,-(\ga\cup m_\al))$ in some gradings whose total dimension is the same as that of $\shg(-M,-\ga)$. The direct sum of these summands are denoted by $\mathcal{SHG}_\al(-M,-\ga)$, which generalizes $\mathcal{I}_+(-\widehat{Y},\widehat{K})$ in \cite[Section 4]{LY2020}.

The above method can be applied iteratively for a tangle $T$ with more than one component and finally we can distinguish all torsions in $H_1(M;\mathbb{Z})$ by choosing $T$ such that $H_1(M_T;\mathbb{Z})$ is torsion-free (\textit{c.f.} Lemma \ref{lem: torsion free}). Similarly, we can identify $\shg(-M,-\ga)$ with a direct summand of $\shg(-M_T,-(\ga\cup m_T))$, where $m_T$ is the union of meridians of tangle components of $T$. Since $H_1(M_T;\mathbb{Z})$ has no torsion, all torsions that are mixed on $H_1(M;\mathbb{Z})$ can be distinguished on $H_1(M_T;\mathbb{Z})$, and this provides the desired decomposition in Theorem \ref{thm: main}.

Suppose $j_*$ is the map on group rings induced by \[j:H_1(M_T;\mathbb{Z})\to H_1(M;\mathbb{Z}).\]Given the construction of $\mathcal{SHG}_T(-M,-\ga)$, Equation (\ref{eq: enhanced same}) reduces to the following equation\begin{equation}\label{eq: A1}
    \en(\shg(-M,-\ga))\deq j_*(\gr(\mathcal{SHG}_T(-M,-\ga)))=\chi(SFH(-M,-\ga)).
\end{equation}We prove this equation by introducing a decomposition $\mathcal{SFH}_T(-M,-\ga)$ of $SFH(-M,-\ga)$ similar to $\mathcal{SHG}_T(-M,-\ga)$. However, the construction of $SFH$ is based on balanced diagrams of balanced sutured manifolds, which is different from the construction of $\shg$ by closures. So we have to show that $SFH$ satisfies the similar setups of $\shg$ to construct $\mathcal{SFH}_T(-M,-\ga)$. This is the main goal of Section \ref{section: Sutured Heegaard Floer homology}, where we collect results for $SFH$ parallel to $\shg$, including gradings associated to admissible surfaces, the surgery exact triangle, the bypass exact triangle, and contact gluing maps.

Since $H_1(M_T;\mathbb{Z})$ is torsion-free, we can apply (\ref{eq: twisted euler}) to $(-M_T,-(\ga\cup m_T))$ to obtain\begin{equation}\label{eq: A2}\gr(\mathcal{SHG}_T(-M,-\ga))=\gr(\mathcal{SFH}_T(-M,-\ga))=\chi(\mathcal{SFH}_T(-M,-\ga)).\end{equation}By discussion on spin$^c$ structures, we show\begin{equation}\label{eq: A3}j_*(\chi(\mathcal{SFH}_T(-M,-\ga)))=\chi(SFH(-M,-\ga)).\end{equation}Equations (\ref{eq: A2}) and (\ref{eq: A3}) imply Equation (\ref{eq: A1}), which finishes the proof of Theorem \ref{thm: main}. The detailed proof can be found in Section \ref{sec: decomposition}.

The proof of Theorem \ref{thm:shm=sfh} is almost straightforward, based on the work of Lekili \cite{lekili2013heegaard}, and Baldwin and Sivek \cite{Baldwin2020}. Since $\mathcal{SHG}_T(-M,-\ga)$ is direct summands of $\shg(-M_T,-\ga\cup m_T)$ in some gradings, it suffices to prove the theorem when $H_1(M;\mathbb{Z})$ is torsion-free. In this case, the decomposition is just induced by admissible surfaces and Theorem \ref{thm:shm=sfh} follows from the isomorphism \[SHM(M,\ga)\cong SFH(M,\ga)\otimes\Lambda.\]The detailed proof can be also found in Section \ref{sec: decomposition}.

In Section \ref{sec: Knots with small instanton knot homology}, we study knots whose $\dim_\mathbb{C}KHI$ are small and prove Theorem \ref{thm: unknot detection}, Theorem \ref{thm: trefoil detection}, and analog theorems in Heegaard Floer theory.


\subsection{Conventions}
If it is not mentioned, all manifolds are smooth, oriented, and connected. Homology groups and cohomology groups are with $\mathbb{Z}$ coefficients, \textit{i.e.} $H_*(M)\deq H_*(M;\mathbb{Z})$ for any manifold $M$. For other coefficients (like $\mathbb{Q})$, we still write $H_*(M;\mathbb{Q})$. We write $\mathbb{Z}_n$ for $\mathbb{Z}/n\mathbb{Z}$. For a simple closed curve on a surface, we do not distinguish between its homology class and itself. The algebraic intersection number of two curves $\al$ and $\be$ on a surface is denoted by $\al\cdot\be$, while the number of intersection points between $\al$ and $\be$ is denoted by $|\al\cap \be|$. A basis $(m,l)$ of $H_1(T^2;\mathbb{Z})$ satisfies $m\cdot l=-1$. The \textbf{surgery} means the Dehn surgery and the slope $q/p$ in the basis $(m,l)$ corresponds to the curve $qm+pl$.




\subsection{Acknowledgement}
 The authors would like to thank Ciprian Manolescu, John A. Baldwin, Steven Sivek, John B. Etnyre, David Shea Vela-Vick, Ian Zemke, Linsheng Wang, Minghao Miao, and Yu Zhao for valuable discussions. The authors are grateful to Sudipta Ghosh and Ian Zemke for pointing out the proof of Lemma \ref{lemma: GZ}. The authors are also grateful to Yi Xie for pointing out the applications about detection results. The authors are grateful to Hongjian Yan for comments on Remark \ref{rem: +4 remark}. The second author would like to thank his supervisor Jacob Rasmussen for patient guidance and helpful comments and thank his parents for support and constant encouragement. The second author is also grateful to Yi Liu for inviting him to BICMR, Peking University.

\section{Twisted sutured homology}\label{Section: Formal sutured homology}

In this section, we collect useful properties of $\shg$. 

\subsection{Notations, gradings, and Euler characteristics}\label{subsec: gradings on SH}\quad

\bdefn[\cite{juhasz2006holomorphic,kronheimer2010knots}]\label{defn_2: balanced sutured manifold}
A \textbf{balanced sutured manifold} $(M,\ga)$ consists of a compact oriented 3-manifold $M$ with non-empty boundary together with a closed 1-submanifold $\ga$ on $\partial{M}$. Let $A(\ga)=[-1,1]\times\ga$ be an annular neighborhood of $\ga\subset \partial{M}$ and let $R(\ga)=\partial{M}\backslash{\rm int}(A(\ga))$. They satisfy the following properties. 
\begin{enumerate}[(1)]
    \item Neither $M$ nor $R(\ga)$ has a closed component.
    \item If $\partial{A(\ga)}=\partial{R(\ga)}$ is oriented in the same way as $\ga$, then we require this orientation of $\partial{R(\ga)}$ induces one on $R(\ga)$. The induced orientation on ${R(\ga)}$ is called the \textbf{canonical orientation}.
    \item Let $R_+(\ga)$ be the part of $R(\ga)$ so that the canonical orientation coincides with the induced orientation on $\partial{M}$, and let $R_-(\ga)=R(\ga)\backslash R_+(\ga)$. We require that $\chi(R_+(\ga))=\chi(R_-(\ga))$. If $\ga$ is clear in the contents, we simply write $R_\pm=R_\pm(\ga)$.
\end{enumerate}

\edefn

The constructions of sutured instanton homology $SHI$ and sutured monopole homology $SHM$ for balanced sutured manifolds were originated by Kronheimer and Mrowka \cite{kronheimer2010knots}. Later, Baldwin and Sivek \cite{baldwin2015naturality} dealt with the naturality problem of these homologies and contructed $\shi$ and $\shm$. After that, several groups of people studies these homologies extensively, see for example \cite{baldwin2016contact,baldwin2018khovanov,li2019direct,li2019decomposition,wang2020cosmetic}. 

In this paper, we will use $\shg$ to denote both $\shi$ and $\shm$ and call it \textbf{twisted sutured homology}. The coefficient field is denoted by $\mathbb{F}$. For closed 3-manifolds and knots with basepoints, we can construct balanced sutured manifolds and then apply twisted sutured homology to them as follows.

\bdefn\label{defn: khg}
Suppose that $Y$ is a closed 3-manifold and $z\in Y$ is a basepoint. Let $Y(1)$ be obtained from $Y$ by removing a 3-ball containing $z$ and let $\delta$ be a simple closed curve on $\partial Y(1)\cong S^2$. Suppose that $K\subset Y$ is a knot and $w$ is a basepoint on $K$. Let $Y(K)$ be the knot complement of $K$ and let $\ga=m\cup (-m)$ consist of two meridians with opposite orientations of $K$ near $w$. Then $(Y(1),\delta)$ and $(Y(K),\ga)$ are balanced sutured manifolds. Define
$$\hhat(Y,z)\deq \shg(Y(1),\delta) ~{\rm and}~\khg(Y,K,w)\deq \shg(Y(K),\ga).$$
\edefn
\begin{conv}
Different choices of the basepoints give isomorphic vector spaces. Since in the current paper we only care about the isomorphism class of the vector spaces, we omit the basepoints and simply write $\hhat(Y)$ and $\khg(Y,K)$ instead.
\end{conv}

If $S$ is a properly embedded surface in $M$ with some admissible conditions. The first author \cite{li2019direct} constructed a grading on $\shg(M,\ga)$ (with some pioneering work done by Kronheimer and Mrowka \cite{kronheimer2010instanton} and Baldwin and Sivek \cite{baldwin2018khovanov}).
\bdefn[{\cite{li2019decomposition}}]\label{defn_2: admissible surfaces}
Suppose $(M,\ga)$ is a balanced sutured manifold and $S\subset M$ is a properly embedded surface. The surface $S$ is called an \textbf{admissible surface} if the followings hold.
\begin{enumerate}[(1)]
    \item Every boundary component of $S$ intersects $\ga$ transversely and nontrivially.
    \item $\frac{1}{2}|S\cap \ga|-\chi(S)$ is an even integer.
\end{enumerate}
\edefn

\bthm[{\cite{li2019direct,li2018gluing}}]\label{thm_2: grading in SHG}
Suppose $(M,\ga)$ is a balanced sutured manifold and $S\subset (M,\ga)$ is an admissible surface. Then there is a $\intg$-grading on $\shg(M,\ga)$ induced by $S$, which we write as
$$\shg(M,\ga)=\bigoplus_{i\in\intg}\shg(M,\ga,S,i).$$
This decomposition satisfies the following properties.

\begin{enumerate}[(1)]
    \item Suppose $n=\frac{1}{2}|\partial S\cap\ga|$. If $|i|>\frac{1}{2}(n-\chi(S))$, then $\shg(M,\ga,S,i)=0.$
    \item If there is a sutured manifold decomposition $(M,\ga)\stackrel{S}{\leadsto}(M',\ga')$ in the sense of Gabai \cite{gabai1983foliations}, then we have
$$\shg(M,\ga,S,\frac{1}{2}(n-\chi(S)))\cong \shg(M',\ga').$$
    \item For any $i\in\intg$, we have
$$\shg(M,\ga,S,i)=\shg(M,\ga,-S,-i).$$
    \item For any $i\in\intg$, we have
$$\shg(M,-\ga,S,i)\cong \shg(M,\ga,S,-i).$$
    \item For any $i\in\intg$, we have
$$\shg(-M,\ga,S,i)\cong {\rm Hom}_\mathbb{F}(\shg(M,\ga,S,-i),\mathbb{F}).$$
\end{enumerate}
\ethm


If $S\subset(M,\ga)$ is not admissible, then we can perform an isotopy on $S$ to make it admissible. 

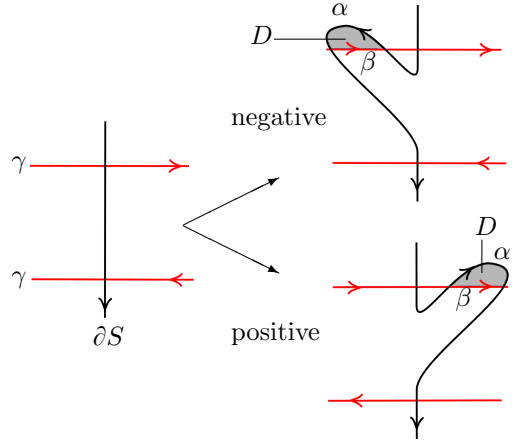
\begin{figure}[ht]
\centering
\begin{overpic}[width=2.5in]{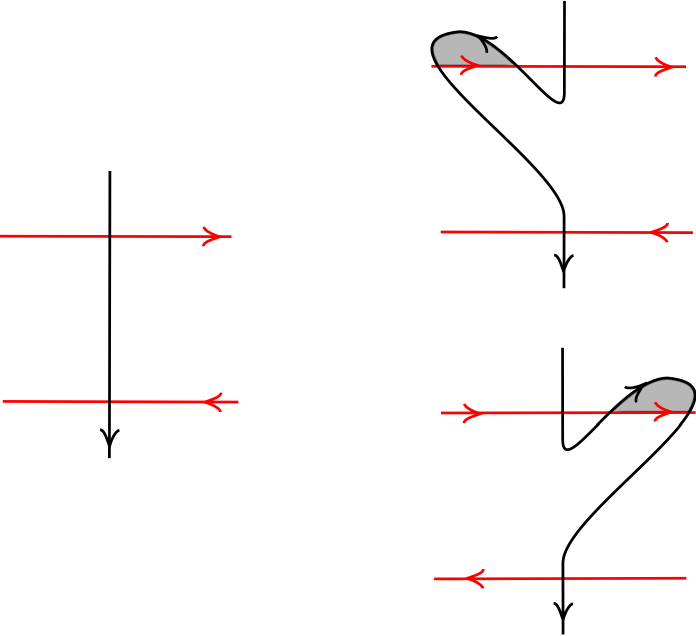}
    \put(13,20){$\partial{S}$}
    \put(-4,57){$\ga$}
    \put(-4,33){$\ga$}
    \put(32,45){\vector(2,1){20}}
    \put(32,45){\vector(2,-1){20}}
    \put(42,21){positive}
    \put(42,66){negative}
    \put(46,83){$D$}
    \put(51,84){\line(1,0){15}}
    \put(63,89){$\al$}
    \put(69,77.5){$\be$}
    \put(93,43){$D$}
    \put(94.5,42){\line(0,-1){7}}
    \put(97,38){$\al$}
    \put(89,28){$\be$}
\end{overpic}
\vspace{0.05in}
\caption{The positive and negative stabilizations of $S$.}\label{fig: pm_stabilization_of_surfaces}
\end{figure}

\bdefn\label{defn_2: stabilization of surfaces}
Suppose $(M,\ga)$ is a balanced sutured manifold, and $S$ is a properly embedded surface. A \textbf{stabilization} of $S$ is a surface $S^\p$ obtained from $S$ by isotopy in the following sense. This isotopy creates a new pair of intersection points:
$$\partial S'\cap\ga=(\partial{S}\cap\ga)\cup \{p_+,p_-\}.$$
We require that there are arcs $\al\subset \partial{S'}$ and $\be\subset \ga$, oriented in the same way as $\partial{S'}$ and $\ga$, respectively, and the followings hold.
\begin{enumerate}[(1)]
    \item $\partial{\al}=\partial{\be}=\{p_+,p_-\}$.
    \item $\al$ and $\be$ cobound a disk $D$ with ${\rm int}(D)\cap (\ga\cup \partial{S}')=\emptyset$.
\end{enumerate}
The stabilization is called \textbf{negative} if $\partial{D}$ is the union of $\al$ and $\be$ as an oriented curve. It is called \textbf{positive} if $\partial{D}=(-\al)\cup\be$. See Figure \ref{fig: pm_stabilization_of_surfaces}. We denote by $S^{\pm k}$ the surface obtained from $S$ by performing $k$ positive or negative stabilizations, repsectively.
\edefn

\brem\label{rem_2: pm switches according to orientations of the suture}
The definition of stabilizations of a surface depends on the orientations of the suture and the surface. If we reverse the orientation of the suture or the surface, then positive and negative stabilizations switch between each other.
\erem

One can also relate the gradings associated to different stabilizations of a fixed surface.

\bthm[{\cite[Proposition 4.3]{li2019direct}} and {\cite[Proposition 4.17]{wang2020cosmetic}}]\label{thm_2: grading shifting property}
Suppose $(M,\ga)$ is a balanced sutured manifold and $S$ is a properly embedded surface in $M$, which intersects the suture $\ga$ transversely. Suppose $S$ has a distinguished boundary component so that all the stabilizations mentioned below are performed on this boundary component. Then, for any $p,k,l\in \intg$ so that the stabilized surfaces $S^{p}$ and $S^{p+2k}$ are both admissible, we have
$$\shg(M,\ga,S^{p},l)=\shg(M,\ga,S^{p+2k},l+k).$$
Note $S^p$ is a stabilization of $S$ as introduced in Definition \ref{defn_2: stabilization of surfaces}, and, in particular, $S^0=S$.
\ethm

If we have multiple admissible surfaces, then they together induce a multi-grading.

\bthm[{\cite[Proposition 1.14]{li2019decomposition}}]\label{thm: zn grading}
Suppose $(M,\ga)$ is a balanced sutured manifold and $S_1,\dots,S_n$ are admissible surfaces in $(M,\ga)$. Then there exists a $\intg^n$-grading on $\shg(M,\ga)$ induced by $S_1,\dots,S_n$, which we write as
$$\shg(M,\ga)=\bigoplus_{(i_1,\dots,i_n)\in\intg^n}\shg(M,\ga,(S_1,\dots,S_n),(i_1,\dots,i_n)).$$
\ethm

\bthm[{\cite[Theorem 1.12]{li2019decomposition}}]\label{thm: shift}
Suppose $(M,\ga)$ is a balanced sutured manifold and $\al\in H_2(M,\partial M)$ is a nontrivial homology class. Suppose $S_1$ and $S_2$ are two admissible surfaces in $(M,\ga)$ such that$$[S_1]=[S_2]=\al\in H_2(M,\partial M).$$Then, there exists a constant $C$ so that $$\shg(M,\ga,S_1,l)=\shg(M,\ga,S_2,l+C).$$
\ethm

Based on the $\mathbb{Z}^n$ grading from Theorem \ref{thm: zn grading}, we can define the graded Euler characteristic.
\bdefn\label{defn: chi(shg)}
Suppose $(M,\ga)$ is a balanced sutured manifold and $S_1,\dots,S_n$ are admissible surfaces in $(M,\ga)$ such that $[S_1],\dots,[S_n]$ generate $H_2(M,\partial M)$. Let $\rho_1,\dots,\rho_n\in H^\p=H_1(M)/{\rm Tors}$ satisfying $\rho_i\cdot S_j=\delta_{i,j}$. The \textbf{graded Euler characteristic} of $\shg(M,\ga)$ is $$\gr(\shg(M,\ga)) \deq\sum_{(i_1,\dots,i_n)\in\intg^n}\chi(\shg(M,\ga,(S_1,\dots,S_n),(i_1,\dots,i_n)))\cdot (\rho_1^{i_1}\cdots\rho_n^{i_n})\in \mathbb{Z}[H^\p]/\pm H^\p.$$
\edefn
\brem\label{rem: fix closure}
By Theorem \ref{thm: shift}, the definition of graded Euler characteristic is independent of the choices of $S_1,\dots,S_n$ if we regard it as an element in $\mathbb{Z}[H^\p]/\pm H^\p$. If the admissible surfaces $S_1,\dots,S_n$ and a particular closure of $(M,\ga)$ is fixed, then the ambiguity of $\pm H^\p$ can be removed.
\erem


\subsection{Contact handles and bypasses}\label{subsection: Contact handles and bypasses}\quad

Suppose $(M,\ga)\subset(M^\p,\ga^\p)$ is a proper inclusion of balanced sutured manifolds and suppose $\xi$ is a contact structure on $M^\p\backslash{\rm int} M$ with dividing sets $\ga^\p\cup(-\ga)$. For monopole theory and instanton theory, Baldwin and Sivek \cite{baldwin2016contact,baldwin2016instanton} (see also \cite{li2018gluing}) constructed a \textbf{contact gluing map}\[\Phi_\xi:\shg(-M,-\ga)\ra \shg(-M^\p,-\ga^\p)\]based on contact handle decompositions and the first author \cite{li2018gluing} showed that the map is functorial, \textit{i.e.} it is indepedent of the contact handle decompositions and gluing two contact structures induces composite maps. In this subsection, we will describe the maps associated to contact 1- and 2-handle attachments, and bypass attachments (\textit{c.f.} \cite{honda2000classification}).

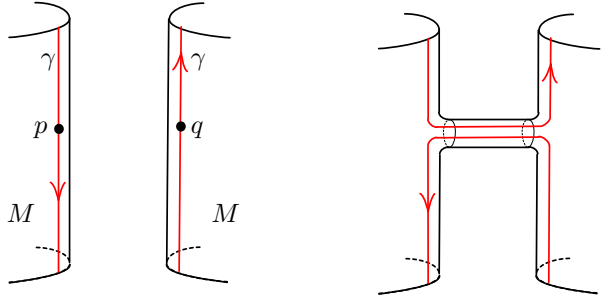
\begin{figure}[ht]
\centering
\begin{overpic}[width=0.6\textwidth]{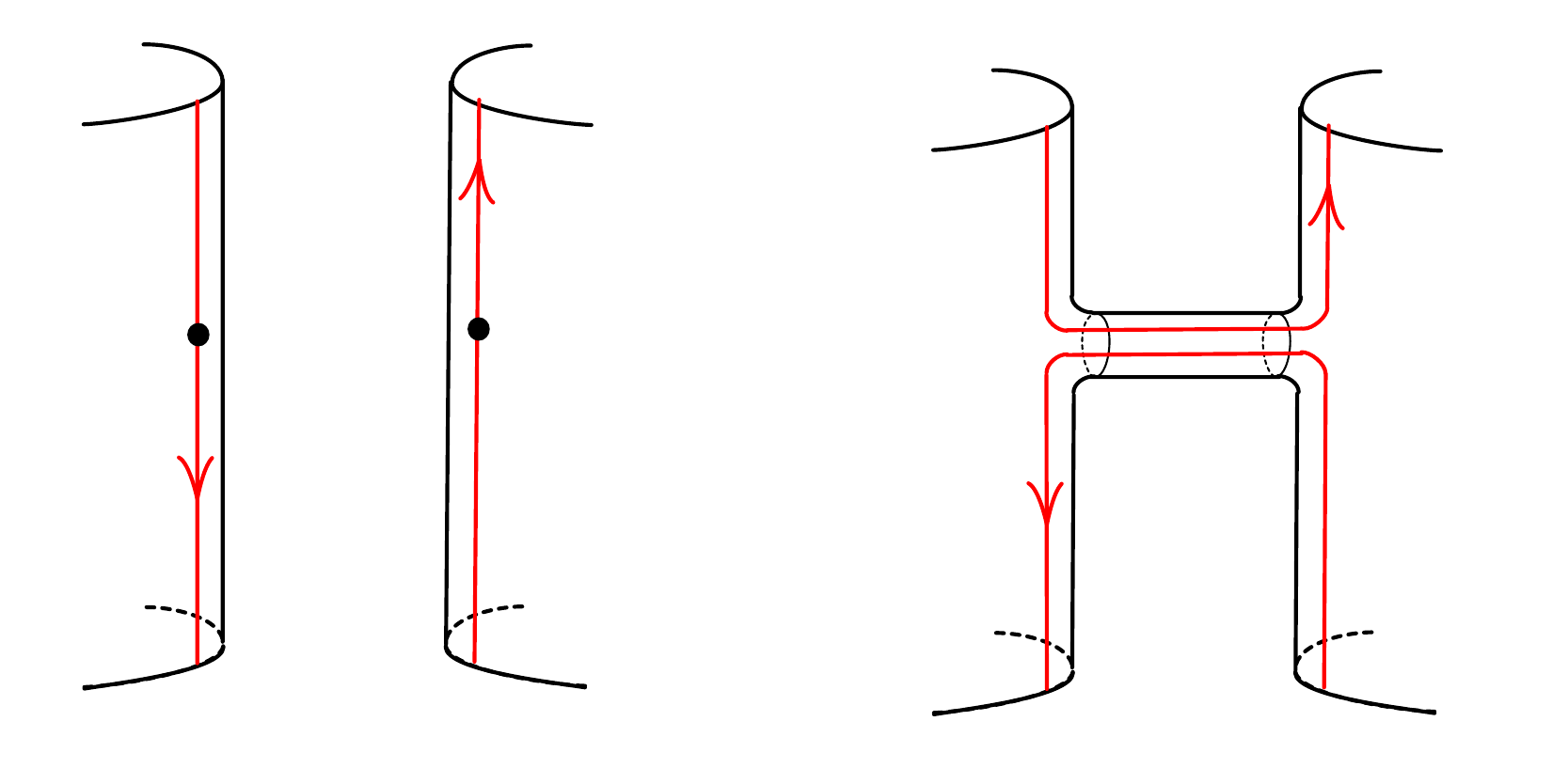}
	\put(9,28){$p$}
	\put(32,28){$q$}
	\put(5,15){$M$}
	\put(35,15){$M$}
	\put(10,38){$\ga$}
	\put(32,38){$\ga$}
\end{overpic}
\caption{Left, the sutured manifold $(M,\ga)$ with two points $p$ and $q$ on the suture. Right, the 1-handle attachment along $p$ and $q$.}
\label{fig:onehandle}
\end{figure}

{\bf Contact 1-handle.} Suppose $D_-$ and $D_+$ are disjoint embedded disks in $\partial M$ which each intersect $\ga$ in a single properly embedded arc. Consider the standard contact structure $\xi_{\rm std}$ on the 3-ball $B^3$. We glue $(D^2\times [-1,1],\xi_{D^2})\cong (B^3,\xi_{\rm std})$ to $(M,\ga)$ by diffeomorphisms \[D^2\times \{-1\} \to D_-~{\rm and}~ D^2\times\{+1\} \to D_+,\]which preserve and reverse orientations, respectively, and identify the dividing sets with the sutures. Then we round corners as shown in Figure \ref{fig:onehandle} (\textit{c.f.} \cite[Figure 2]{baldwin2016instanton}). Let $(M_1,\ga_1)$ be the resulting sutured manifold.

Suppose $(Y,R)$ is a closure of $(M_1,\ga_1)$. By \cite[Section 3.2]{baldwin2016instanton}, it is also a closure of $(M,\ga)$. Define the map associated to the contact 1-handle attachment by the identity map\[C_{h^1}=C_{h^1,D_-,D_+}\deq {\rm id}:\shg(-M,-\ga)\xra{=}\shg(-M_1,-\ga_1).\]

\begin{figure}[ht]
\centering
\begin{overpic}[width=0.6\textwidth]{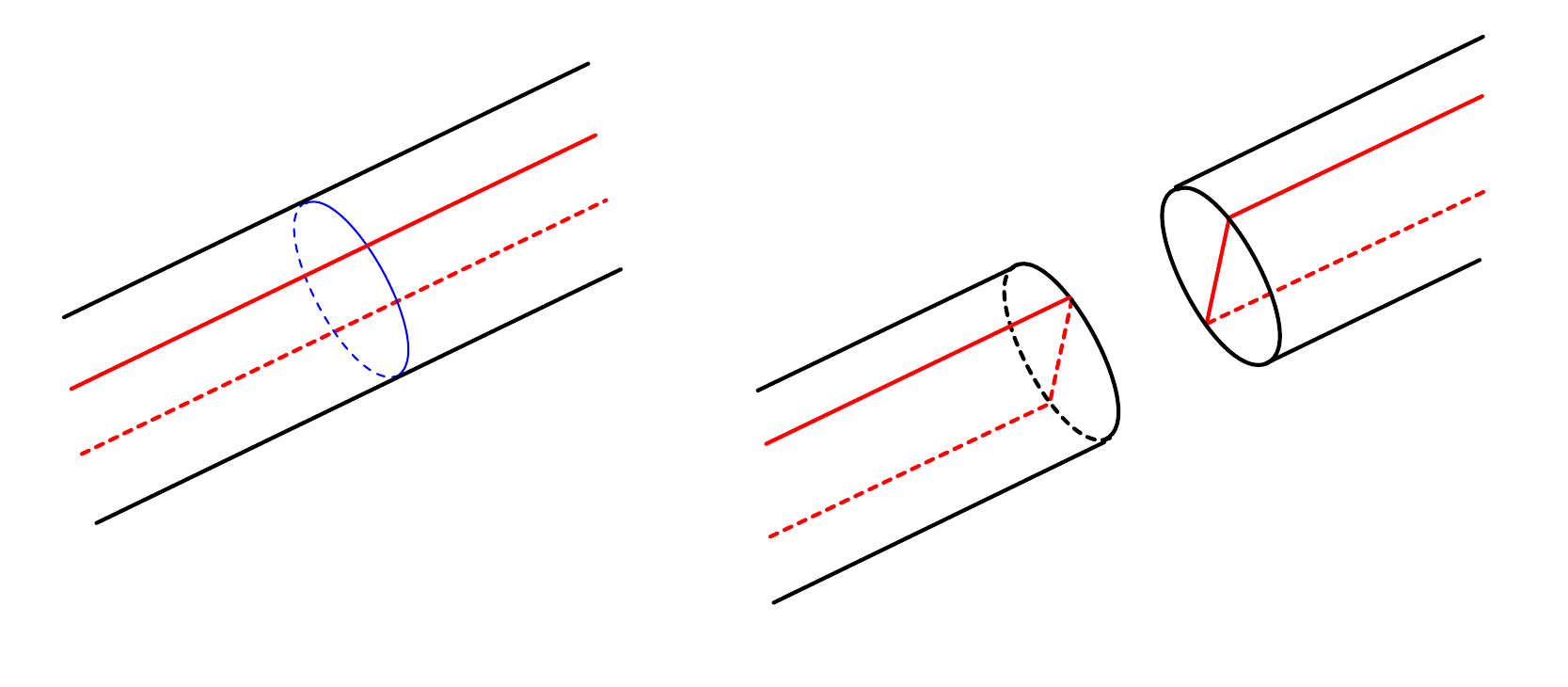}
	\put(2,16){$\ga$}
	\put(17,32){$\mu$}
	\put(5,32){$M$}
	\put(30,14){$M$}
\end{overpic}
\caption{Left, the sutured manifold $(M,\ga)$ and the curve $\be\subset\partial M$ that intersects $\ga$ at two points. Right, the 2-handle attachment along the curve $\mu$.}
\label{fig:twohandle2}
\end{figure}
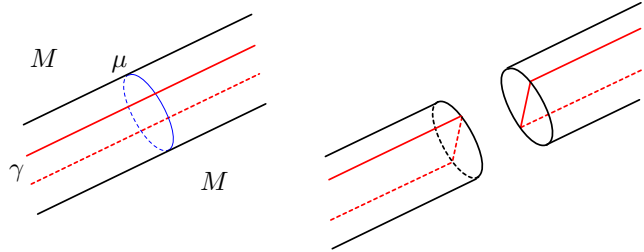

{\bf Contact 2-handle.} Suppose $\mu$ is an embedded curve in $\partial M$ which intersects $\ga$ in two points. Let $A(\mu)$ be an annular neighborhood of $\lambda$ intersecting $\ga$ in two cocores. We glue $(D^2\times [-1,1],\xi_{D^2})\cong (B^3,\xi_{\rm std})$ to $(M,\ga)$ by an orientation-reversing diffeomorphism \[\partial D^2\times [-1,1]\to A(\mu),\]which identifies positive regions with negative regions. Then we round corners as shown in Figure \ref{fig:twohandle2} (\textit{c.f.} \cite[Figure 3]{baldwin2016instanton}). Let $(M_2,\ga_2)$ be the resulting sutured manifold.

We construct the map associated to the contact 2-handle attachment as follows. Let $\mu^\p$ be the knot obtained by pushing $\mu$ into $M$ slightly. Suppose $(N,\ga_{N})$ is the manifold obtained from $(M,\ga)$ by a $0$-surgery along $\mu^\p$ with respect to the framing from $\partial N$. By \cite[Section 3.3]{baldwin2016instanton}, the sutured manifold $(N,\ga_{N})$ can be obtained from $(M_2,\ga_2)$ by attaching a contact 1-handle. Since $\mu'\subset {\rm int}(M)$, the construction of the closure of $(M,\ga)$ does not affect $\mu'$. Thus, we can construct a cobordism between closures of $(M,\ga)$ and $(N,\ga_{N})$ by attaching a $4$-dimensional 2-handle associated to the surgery on $\mu'$. This cobordism induces a cobordism map
$$C_{\mu'}:\shg(-M,-\ga)\to \shg(-N,-\ga_{N}).$$Consider the identity map$$\iota:\shg(-M_2,-\ga_2)\xra{=}\shg(-N,-\ga_{N}).$$Define the the map associated to the contact 2-handle attachment as$$C_{h^2}=C_{h^2,\mu}\deq\iota^{-1}\circ C_{\mu'}:\shg(-M,-\ga)\to\shg(-M_2,-\ga_2).$$

\begin{figure}[ht]
\centering
\begin{overpic}[width=0.7\textwidth]{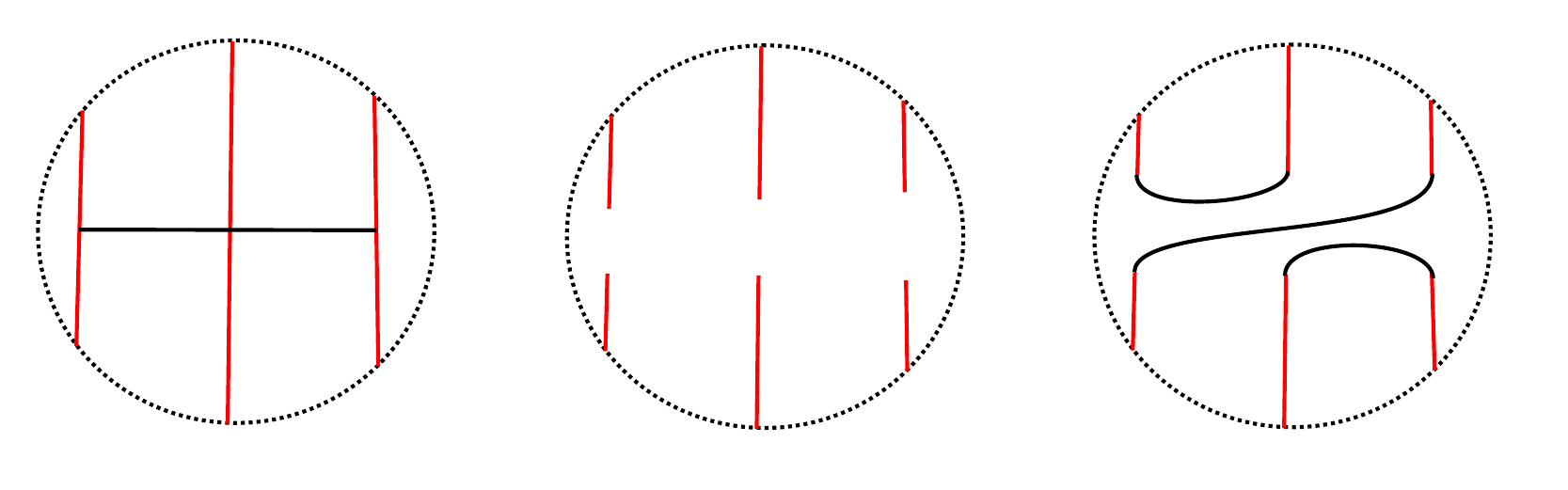}
    \put(18,17.5){$\alpha$}
\end{overpic}
\caption{The bypass arc and the bypass attachment, where the orientation of $\partial M$ is pointing out.}\label{fig: the bypass arc}
\end{figure}

{\bf Bypass attachment.} Suppose $\al$ is an embedded arc in $\partial M$ which intersects $\ga$ in three points. Let $D$ be a disk neighborhood of $\al$ intersecting $\ga$ in three arcs. There are six endpoints after cutting $\ga$ along $\al$. We replace three arcs in $D$ with another three arcs as shown in Figure \ref{fig: the bypass arc}. Let $(M,\ga^\p)$ be the resulting sutured manifold. The arc $\al$ is called a \textbf{bypass arc} and this procedure is called \textbf{bypass attachment} along $\al$.

By Ozbagci \cite[Section 3]{ozbagci2011contact}, the bypass attachment can be recovered by contact handle attachments as follows. First, one can attach a contact 1-handle along two endpoints of $\al$. Then one can attach a contact 2-handle along a circle that is the union of $\al$ and an arc on the attached 1-handle. Topologically, the 1-handle and the 2-handle form a canceling pair, so the diffeomorphism type of the 3-manifold does not change. However, the contact structure is changed, and the suture $\ga$ is replaced by $\ga^\p$. We define the \textbf{bypass map} associated to the bypass attachment as$$\psi_\al\deq C_{h^2}\circ C_{h^1}:\shg(-M,-\ga)\ra\shg(-M,-\ga^\p).$$

In \cite{baldwin2016contact,baldwin2018khovanov}, Baldwin and Sivek proved the bypass exact triangle for sutured monopole Floer homology and sutured instanton Floer homology.

\begin{figure}[htbp]
\centering
\begin{minipage}[t]{0.35\textwidth}
\centering
\begin{overpic}[width=0.9\textwidth]{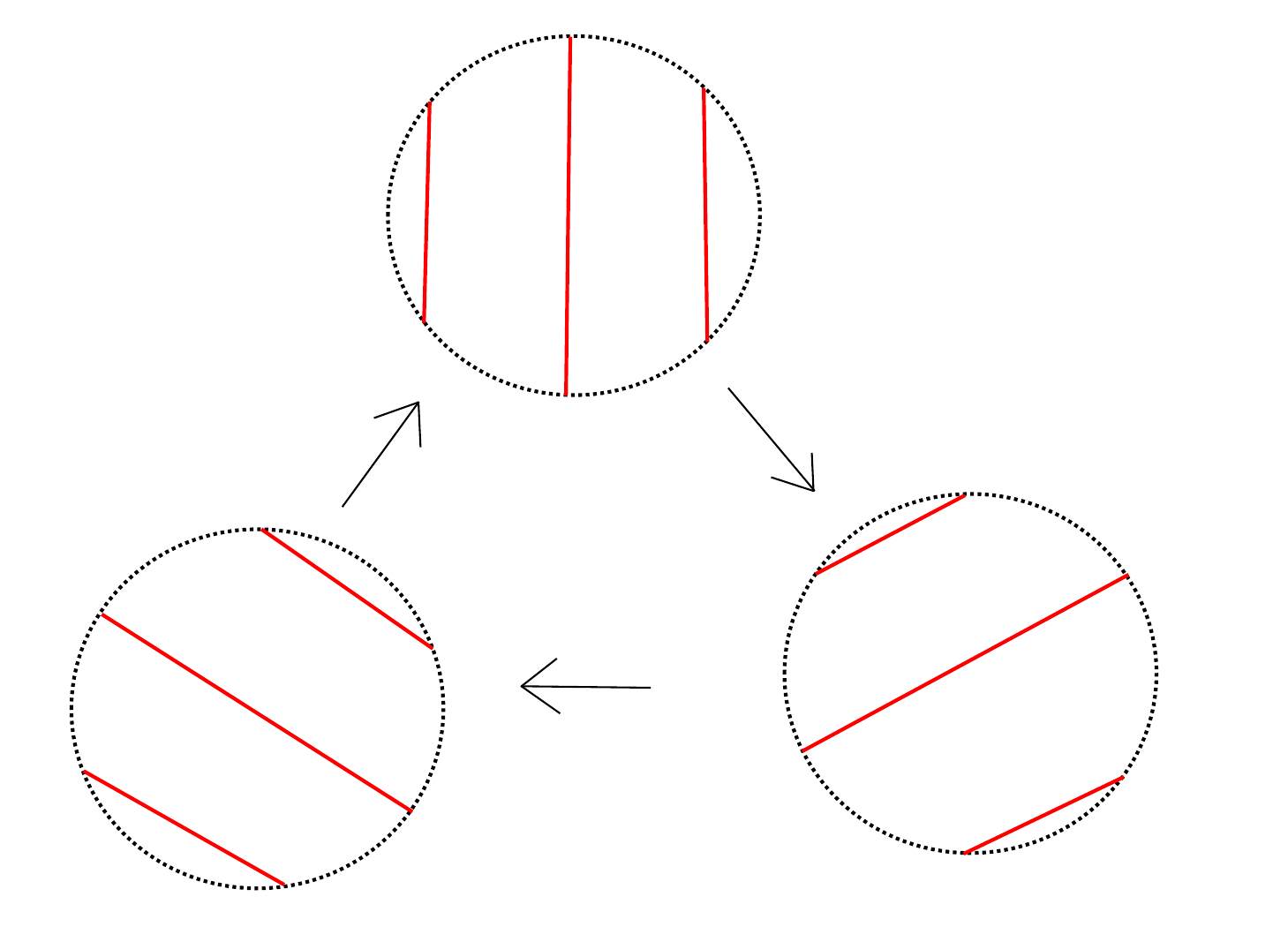}
    \put(16,22){$\ga_3$}
    \put(38,58){$\ga_1$}
    \put(73,28){$\ga_2$}
\end{overpic}
\vspace{-0.05in}
\caption{The bypass triangle.}\label{fig: the bypass triangle}
\end{minipage}
\begin{minipage}[t]{0.4\textwidth}
\centering
\begin{overpic}[width=0.9\textwidth]{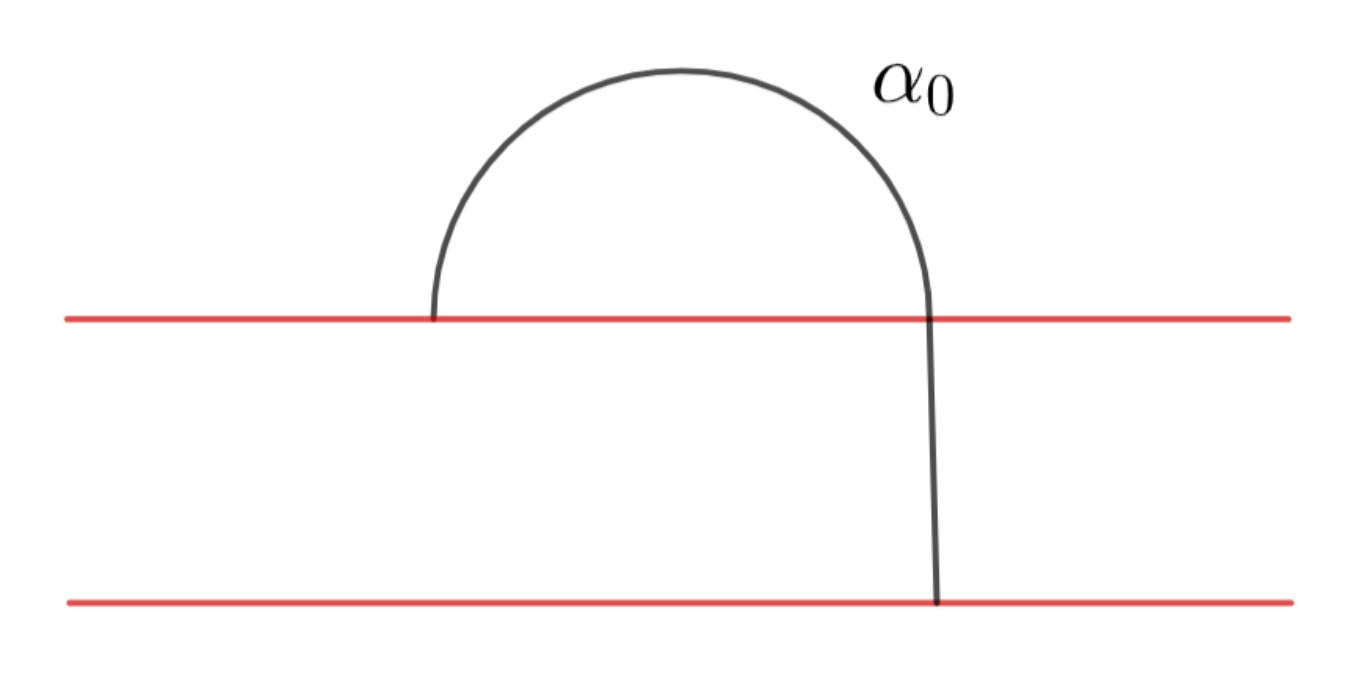}
\end{overpic}
\caption{A trivial bypass.}\label{trivial bypass}
\end{minipage}
\end{figure}

\bthm[{\cite[Theorem 5.2]{baldwin2016contact} and \cite[Theorem 1.21]{baldwin2018khovanov}}]\label{thm_2: bypass exact triangle on general sutured manifold}
Suppose $(M,\ga_1)$, $(M,\ga_2)$, $(M,\ga_3)$ are balanced sutured manifolds such that the underlying 3-manifolds are the same, and the sutures $\ga_1$, $\ga_2$, and $\ga_3$ only differ in a disk shown in Figure \ref{fig: the bypass triangle}.  Then there exists an exact triangle
\begin{equation}\label{eq: bypass exact triangle}
\xymatrix@R=6ex{
\shg(-M,-\ga_1)\ar[rr]^{\psi_1}&&\shg(-M,-\ga_2)\ar[dl]^{\psi_2}\\
&\shg(-M,-\ga_3)\ar[lu]^{\psi_3}&
}    
\end{equation}

Moreover, the maps $\psi_i$ are induced by cobordisms, hence is homogeneous with respect to the relative $\mathbb{Z}_2$ grading on $\shg(M,\ga_i)$.
\ethm

The following proposition is straightforward from the description of the bypass map.
\bprop\label{prop: bypass maps preserves gradings}
Suppose $(M,\ga)$ is a balanced sutured manifold and $S\subset (M,\ga)$ is an admissible surface. Suppose the disk as in Figure \ref{fig: the bypass triangle}, where we perform the bypass change, is disjoint from $\partial S$. Let $\ga_2$ and $\ga_3$ be the resulting two sutures. Then all the maps in the bypass exact triangle (\ref{eq: bypass exact triangle}) are grading preserving, \textit{i.e.}, for any $i\in\intg$, we have an exact triangle
\begin{equation*}
\xymatrix@R=6ex{
\shg(-M,-\ga_1,S,i)\ar[rr]^{\psi_{1,i}}&&\shg(-M,-\ga_2,S,i)\ar[dl]^{\psi_{2,i}}\\
&\shg(-M,-\ga_3,S,i)\ar[lu]^{\psi_{3,i}}&
}    
\end{equation*}
where $\psi_{k,i}$ are the restriction of $\psi_k$ in (\ref{eq: bypass exact triangle}).
\eprop
A special bypass arc $\al_0$ is depicted in Figure \ref{trivial bypass}, where the bypass attachment along $\al$ is called a \textbf{trivial bypass} (\textit{c.f.} \cite[Section 2.3]{honda2002gluing}). Attaching a trivial bypass does not change the suture on $\partial M$ and induces a product contact structure on $\partial M\times I$. The functoriality of the contact gluing maps indicates the following proposition.
\bprop[]\label{prop: trivial bypass}
A trivial bypass on $(M,\ga)$ induces an identity map on $\shg(M,\ga)$.
\eprop

\section{Decomposition associated to tangles}\label{Section: Decomposition associated to tangles}

Suppose $K$ is a rationally null-homologous knot in a closed 3-manifold $Y$, \textit{i.e.} $[K]=0\in H_1(Y;\mathbb{Q})$. Suppose $q$ is the order of $[K]$, \textit{i.e.} $q$ is the smallest number satisfying $q[K]=0\in H_1(Y;\intg)$. In \cite[Section 4]{LY2020}, we construct a decomposition$$I^\sharp(Y)\cong \bigoplus_{i=0}^q I^\sharp(Y,i).$$This decomposition provides a candidate for the counterpart of the torsion spin$^c$ decompositions in monopole theory and Heegaard Floer theory. 

In this section, we generalize this decomposition to rationally null-homologous tangles in balanced sutured manifolds. There is no essential difference between the proofs for knots and tangles. All arguments apply to both sutured instanton and monopole homology, so we can safely use $\shg$. 

\subsection{Basic setups}\quad
\label{Subsection: Basic setups}

In this subsection, we review the construction for tangles and collect important lemmas in \cite[Section 3.2]{LY2020}, with mild modifications.

Suppose $(M,\ga)$ is a balanced sutured manifold. Suppose $T=T_1\cup \cdots\cup T_m$ is a \textbf{vertical tangle} in $(M,\ga)$ (\textit{c.f.} \cite[Definition 1.1]{xie2019tangle}), \textit{i.e.} a properly embedded 1-submanifold with \[|T_i\cap R_+(\ga)|=|T_i\cap R_-(\ga)|=1.\]  Let $T_i$ be oriented from $R_+(\ga)$ to $R_-(\ga)$. Throughout this subsection, we consider one component $\al$ of $T$ and assume it is rationally null-homologous, \textit{i.e.} $[\al]=0\in H_1(M,\partial M;\mathbb{Q})$. Without loss of generality, suppose $\al=T_1$. 

We can construct a new balanced sutured manifold $(M_T,\ga_T)$ as follows. Let $M_T$ be obtained from $M$ by removing a neighborhood $N(T)=\bigcup_{i=1}^m N(T_i)$ of $T$. Suppose $\ga_i$ is a positively oriented meridian of $T_i$ on $\partial N(T_i)$. Define \[\ga_T=\ga\cup \ga_1\cup\cdots\cup \ga_m.\] 


Since $\al$ is rationally null-homologous, there exists a surface $S$ in $M$ with $\partial S$ consisting of arcs $\be_1$,\dots, $\be_k$ and $q$ copies of $\al$ for some integers $k$ and $q$. Here $q$ is the order of $\al$.

The surface $S$ can be modified into a properly embedded surface $S_T$ in $M_T$ as follows. First, for $q$ arcs in $\partial S$ parallel to $\al$, we isotop them to be on $\partial N(\al)$. Then $\be_1,\dots,\be_k$ can be regarded as arcs on $\partial M_T$. Second, We can isotop $S$ to make it intersect $T_2,\dots,T_m$ transversely. Then removing disks in $N(T_i)\cap S$ for all $i=2,\dots,m$ induces a properly embedded surface $S_T$ in $M_T$. Note that $\partial S_T$ intersects $\ga_1$ at $q$ points, one for each arc parallel to $\al$, and the part of $\partial S_T$ on $\partial N(T_i)$ consists of circles parallel to $\ga_i$ for $i=2,\dots,m$.

Suppose $p_+$ and $p_-$ are the endpoints of $\al$ on $R_+(\ga)$ and $R_-(\ga)$, respectively. Choose an arc $\zeta_+\subset R_+(\ga)$ connecting $p_+$ and $\ga$. The arc $\zeta_+$ induces an arc on $R_+(\ga_T)$ connecting $\ga_1$ to $\ga$ such that the part on $\partial N(\al)$ is parallel to $\al$. We still denote this arc by $\zeta_+$ for simplicity. Similarly we can choose an arc $\zeta_-\subset R_-(\ga_T)$ connecting $\ga_1$ to $\ga$.

\begin{figure}[ht]
\centering
\begin{overpic}[width=0.99\textwidth]{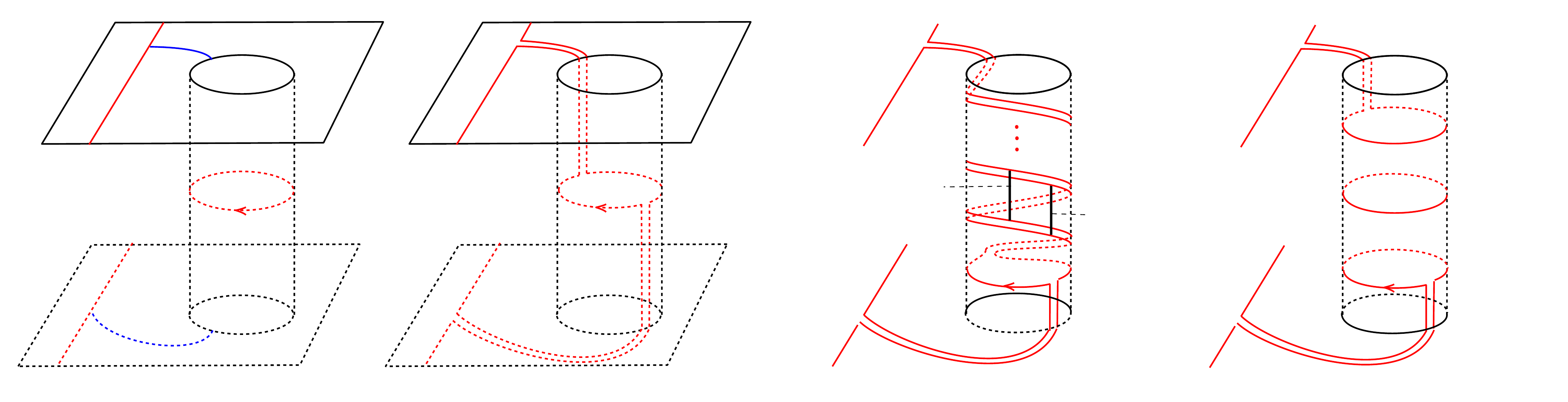}
    \put(9,12){$\ga_1$}
    
    \put(17,22){$R_{+}(\ga_T)$}
    \put(17,3){$R_{-}(\ga_T)$}

    \put(9,20){$\zeta_+$}
    \put(8,5){$\zeta_-$}
    
    \put(12,0){$\Ga_-=\ga_T$}
    \put(36,0){$\Ga_0$}
    \put(62,0){$\Ga_n$}
    \put(86,0){$\Ga_+$}
    
    \put(58,14){$\eta_+$}
    \put(70,12){$\eta_-$}
\end{overpic}
\caption{The arcs $\zeta_+,\zeta_-$, the sutures $\Ga_-,\Ga_0,\Ga_n,\Ga_+$, and the bypass arcs $\eta_+,\eta_-$.}\label{fig: the new Ga_0}
\end{figure}

Let $\Ga_0$ be obtained from $\ga_T$ by band sum operations along $\zeta_+$ and $\zeta_-$. Then let $\Ga_n$ be obtained from $\Ga_0$ by twisting along $(-\ga_1)$ for $n$ times. Moreover, let $\Ga_+$ be the suture as depicted in Figure \ref{fig: the new Ga_0} and let $\Ga_-=\ga_T$.

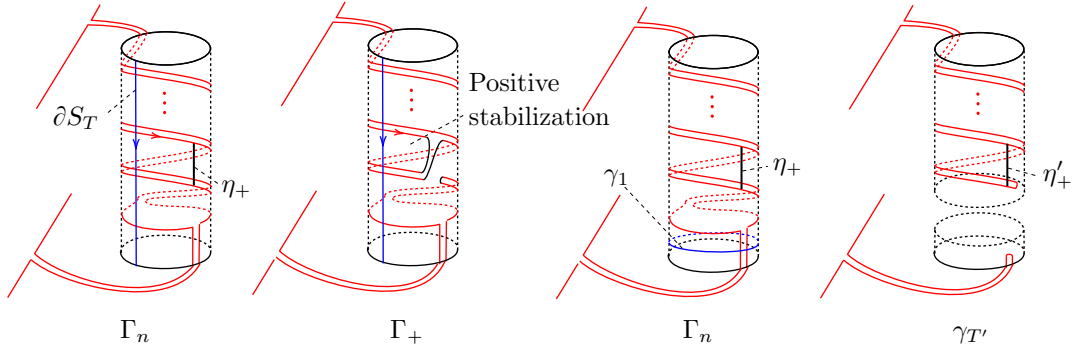
\begin{figure}[ht]
\centering
\begin{overpic}[width=0.99\textwidth]{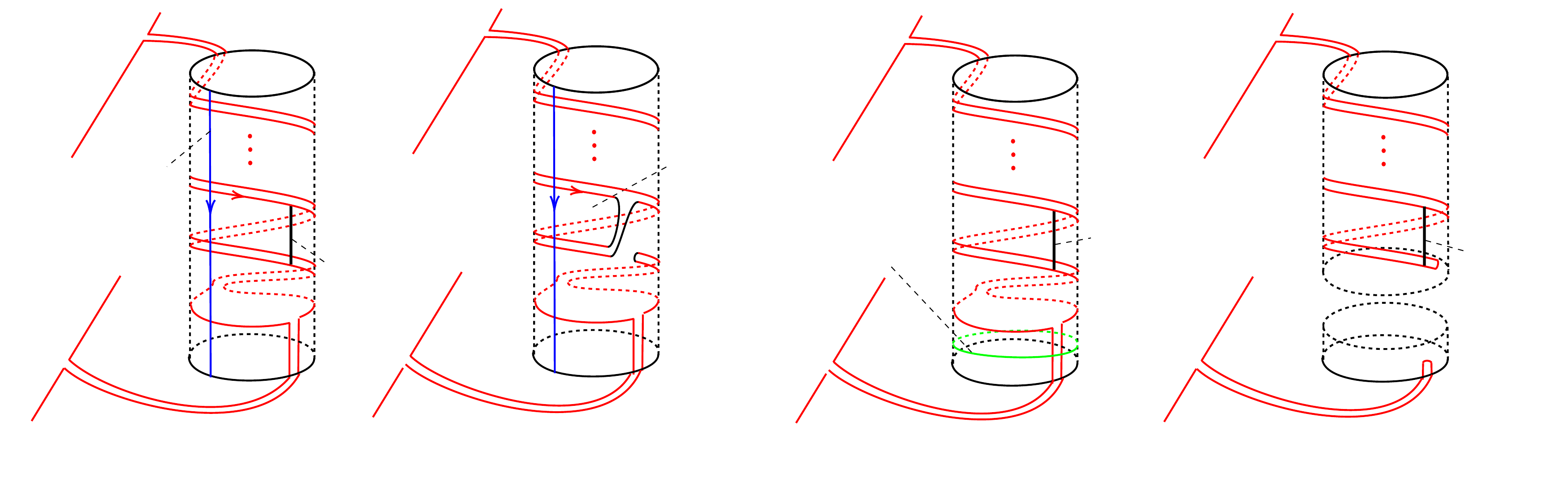}
    \put(21,13){$\eta_+$}
    \put(6,19){$\partial S_T$}

    \put(43,22){Positive}
    \put(43,19){stabilization}
    \put(12,0){$\Ga_n$}
    \put(36,0){$\Ga_+$}
    \put(62,0){$\Ga_n$}
    \put(86,0){$\ga_{T^\p}$}
    
    \put(55,14){$\ga_1$}
    \put(70,15){$\eta_+$}
    \put(94,14){$\eta_+^\p$}
\end{overpic}
\caption{Left two subfigures, the bypass attachment along $\eta_+$. Right two subfigures, the bypass arcs before and after the contact 2-handle attachment.}\label{fig: the new Ga_n with 2-handle attached}
\end{figure}

\brem
The construction of $\zeta_+$ and $\zeta_-$ here is a little different from the one in \cite[Section 3.2]{LY2020}, where we used $\be_1$ to construct $\zeta_\pm$ and removed a trivial tangle from $M_T$ to obtain a manifold $M_{T_0}$. Hence the construction of $\Ga_n,\Ga_\pm$ is also different. In particular, they were on $M_{T_0}$ in the construction of \cite[Section 3.2]{LY2020}. However, it turns out that removing the trivial tangle is not necessary and we can decompose $M_{T_0}$ along a product disk to recover $M_T$ in \cite[Section 3.2, Step 3]{LY2020}. Thus, we can consider sutures on $M_T$ and all results in \cite[Section 3.2]{LY2020} apply without essential change. Also, the conditions that $\zeta_\pm$ are disjoint from $\be_1,\dots,\be_k$ are not essential.
\erem

There are two straightforward choices of bypass arcs on $\Ga_n$ in Figure \ref{fig: the new Ga_0}, denoted by $\eta_+$ and $\eta_-$, respectively. It is straightforward to check that these two bypass arcs induce the following bypass exact triangles from Theorem \ref{thm_2: bypass exact triangle on general sutured manifold} (\textit{c.f.} the left two subfigures of Figure \ref{fig: the new Ga_n with 2-handle attached}).

\begin{equation}\label{eq_3: bypass triangle, pm}
    \xymatrix@R=6ex{
    \shg(-M_{T},-\Ga_{n-1})\ar[rr]^{\psi_{\pm,n}^{n-1}}&&\shg(-M_{T},-\Ga_{n})\ar[dl]^{\psi_{\pm,\pm}^{n}}\\
    &\shg(-M_{T},-\Ga_{\pm})\ar[ul]^{\psi_{\pm,n-1}^{\pm}}&
    }
\end{equation}
The bypasses are attached along $\eta_+$ and $\eta_-$ from the exterior of the 3-manifold $M_{T_0}$, though the point of view in Figure \ref{fig: the new Ga_0} is from the interior of the manifold. So readers have to take extra care when performing these bypass attachments.

Since the bypass arcs $\eta_+$ and $\eta_-$ are disjoint from $\partial S_T$, the bypass maps in the exact triangles (\ref{eq_3: bypass triangle, pm}) preserve gradings associated to $S_T$ by Proposition \ref{prop: bypass maps preserves gradings}. We describe it precisely as follows.
\bdefn\label{defn_3: shifting the grading}
Suppose $(M,\ga)$ is a balanced sutured manifold and $S$ is an admissible surface in $(M,\ga)$. For any $i,j\in\intg$, define
$$\shg(M,\ga,S,i)[j]=\shg(M,\ga,S,i-j).$$
\edefn 

\blem[{\cite[Lemma 3.15 and Lemma 3.18]{LY2020}}]\label{lem_3: graded bypass triangle}
For any $j\in\mathbb{N}\cup\{+,-\}$, there exists admissible surfaces $S_j$ with respect to $(M_T,\Ga_j)$ obtained from $S_T$ by stabilizations and integers $i_{max}^j$ and $i_{min}^j$ such that
\[\lim_{n\to+\infty}i_{max}^n=+\infty,\lim_{n\to+\infty}i_{min}^n=-\infty,\]and
\[\shg(-M_T,-\Ga_j,S_j,i)=0\text{ for }i\not\in [i_{min}^j,i_{max}^j].\]Moreover, for any $n\in\mathbb{N}$, there are two exact triangles
\begin{equation*}\label{eq_6: graded bypass, +}
\xymatrix{
\shg(-M_{T},-{\Ga}_{n},S_n)[{i}_{min}^{n+1}-{i}_{min}^{n}]\ar[rr]^{\quad\quad\psi^{n}_{+,n+1}}&&\shg(-M_{T},-{\Ga}_{n+1},S_{n+1})\ar[dll]^{\psi^{n+1}_{+,+}}\\
\shg(-M_{T},-{\Ga}_{+},S_{+})[{i}_{max}^{n+1}-{i}_{max}^{+}]\ar[u]^{\psi^{+}_{+,n}}&&
}    
\end{equation*}
and
\begin{equation*}\label{eq_3: graded bypass, -}
\xymatrix{
\shg(-M_{T},-{\Ga}_{n},S_{n})[{i}_{max}^{n+1}-{i}_{max}^n]\ar[rr]^{\quad\quad\psi^{n}_{-,n+1}}&&\shg(-M_{T},-{\Ga}_{n+1},S_{n+1})\ar[dll]^{\psi^{n+1}_{-,-}}\\
\shg(-M_{T},-{\Ga}_{-},S_{-})[{i}_{min}^{n+1}-{i}_{min}^{-}]\ar[u]^{\psi^{-}_{-,n}}&&
}.   
\end{equation*}
Furthermore, all maps in the above two exact triangles are grading preserving.
\elem
\brem\label{rem: imin imax}
For $j\in\mathbb{N}\cup\{+,-\}$, the surfaces $S_j$ and integers $i_{max}^j,i_{min}^j$ were defined explicitly in \cite[Step 2 in Section 3.2]{LY2020} by some stabilizations. However, three conditions about $S_T$ at the start of \cite[Step 2 in Section 3.2]{LY2020} are not necessary. We can choose $S_j$ to be either $S_T$ or $S_T^{-1}$ (the negative stabilization of $S_T$ with respect to $\Ga_j$, \textit{c.f.} Definition \ref{defn_2: stabilization of surfaces}), which is admissible with respect to $\Ga_j$. The choice is denoted by $S_j=S_T^{\tau(j)}$ for $\tau(j)\in\{0,-1\}$. Explicitly, $\tau(j)=0$ if $S_j$ is admissible and $\tau(j)=-1$ if $S_j$ is not. For the definitions of $i_{max}^j,i_{min}^j$, consider the closure $(Y_j,R_j)$ of $(M_T,\Ga_j)$ such that $S_j$ extends to a closed surface $\bar{S}_j\subset Y_j$. Define
\begin{equation*}
i^j_{max}=-\frac{1}{2}\chi(\bar{S}_j), ~{\rm and}~i^j_{min}= \frac{1}{2}\chi(\bar{S}_j)-\tau(j).
\end{equation*}
Moreover, we have \[\chi(\bar{S}_j)=\chi(S_j)-\frac{1}{2}|S_j\cap \Ga_j|\]Hence	\[\chi(\bar{S}_-)=\chi(\bar{S}_+)-q+\tau(-)~{\rm and}~\chi(\bar{S}_n)=\chi(\bar{S}_+)-nq+\tau(n)\text{ for }n\in \mathbb{N},\]where $q$ is the order of the tangle $\al$. Since the surfaces after stabilizations are disjoint from the bypass arcs, the bypass maps are grading preserving by Proposition \ref{prop: bypass maps preserves gradings}. The vanishing results follow from term (2) of Theorem \ref{thm_2: grading in SHG}, and \textit{a priori} we do not know if $\shg(-M_T,-\Ga_j,S_j,i)$ is non-vanishing for $i\in\{i_{min}^j,i^{j}_{min}\}$.
\erem

From the vanishing results and the exact triangles in Lemma \ref{lem_3: graded bypass triangle}, the following lemma is straightforward. For any $i\in\intg,n\in\mathbb{N}$, let $\psi_{\pm,n+1}^{n,i}$ be the restriction of $\psi_{\pm,n+1}^{n}$ on the $i$-th grading associated to $S_n$.

\begin{lem}[{\cite[Lemma 3.20]{LY2020}}]\label{lem_3: iso at particular gradings}
	The map
$$\psi_{+,n+1}^{n,i}:\shg(-M_{T},-\Ga_n,S_n,i)\ra\shg(-M_{T},-\Ga_{n+1},S_{n+1},i-(i^{n}_{min}-i^{n+1}_{min}))$$
is an isomorphism if
$ i<P_n\deq~i_{min}^n+(n+1)q-\tau(+).$
Similarly, the map
$$\psi_{-,n+1}^{n,i}:\shg(-M_{T},-\Ga_n,S_n,i)\ra\shg(-M_{T},-\Ga_{n+1},S_{n+1},i+(i^{n+1}_{max}-i^{n}_{max}))$$
is an isomorphism if 
$ i>\rho_n\deq ~ i^{n}_{max}-nq.$
\end{lem}
\bpf
Note that \beq  i^{n+1}_{max}+(i^{n}_{min}-i^{n+1}_{min})-(i^{+}_{max}-i^{+}_{min})
=&~i_{min}^n+(i_{max}^{n+1}-i_{min}^{n+1})-(i^{+}_{max}-i^{+}_{min})\\
=&~i_{min}^n+(-\chi(\bar{S}_+)+(n+1)q)-(-\chi(\bar{S}_+)+\tau(+))\\
=&~i_{min}^n+(n+1)q-\tau(+).
\eeq
\beq  i^{n+1}_{min}-(i^{n+1}_{max}-i^{n}_{max})+(i^{-}_{max}-i^{-}_{min})
=&~ i^{n}_{max}-(i^{n+1}_{max}-i^{n+1}_{min})+(i^{-}_{max}-i^{-}_{min})\\
=&~ i^{n}_{max}-(-\chi(\bar{S}_+)+(n+1)q)+(-\chi(\bar{S}_+)-q)\\
=&~ i^{n}_{max}-nq.\eeq
\epf
There is another important exact triangle induced by the surgery exact triangle.
\blem[{\cite[Lemma 3.21]{LY2020}}]\label{lem_3: surgery triangle}
Suppose $T'=T\backslash\al=T_2\cup\cdots\cup T_m$. Then for any $n\in\mathbb{N}$, there is an exact triangle
\begin{equation}\label{eq_6: surgery exact triangle}
\xymatrix@R=6ex{
\shg(-M_{T},-{\Ga}_{n})\ar[rr]&&\shg(-M_{T},-{\Ga}_{n+1})\ar[dl]^{F_{n+1}}\\
&\shg(-M_{T'},-\ga_{T'})\ar[ul]^{G_n}&
}    
\end{equation}
Furthermore, we have two commutative diagrams related to $\psi^{n}_{+,n+1}$ and $\psi^{n}_{-,n+1}$, respectively
\begin{equation*}
\xymatrix@R=6ex{
\shg(-M_{T},-{\Ga}_{n})\ar[rr]^{\psi_{\pm,n+1}^n}&&\shg(-M_{T},-{\Ga}_{n+1})\\
&\shg(-M_{T'},-\ga_{T'})\ar[ul]^{G_n}\ar[ur]_{G_{n+1}}&
}    
\end{equation*}
\elem
The proof of this lemma is used in Section \ref{section: Sutured Heegaard Floer homology}. So we sketch the proof for the reader's convenience.

\bpf[Sketch of the proof of Lemma \ref{lem_3: surgery triangle}]
Recall that $\ga_1$ is the meridian of $\al$ on $\partial M_T$. Let $\ga_1^\p$ be the curve obtained by pushing $\ga_1$ into the interior of $M_{T}$, with the framing induced by $\partial M_T$. The exact triangle in (\ref{eq_6: surgery exact triangle}) comes from the surgery exact triangle along $\ga_1$: $1$-surgery twists the suture and leads to $\Gamma_{n+1}$; $0$-surgery corresponds to a $2$-handle attachment along $\ga_1$ and hence fills the tangle which leads to the sutured manifold $(M_{T^\p},\ga_{T^\p})$. The commutativity coming from the fact that the surgery curve $\ga_1^\p$ and the bypass arc $\eta_+^\p$ are disjoint from each other, and the fact that the corresponding bypass after the $0$-surgery along $\ga_1^\p$ is trivial.
\epf

\brem
Indeed, we have two more commutative diagrams about $F_n$:
\begin{equation*}
\xymatrix@R=6ex{
\shg(-M_{T},-{\Ga}_{n})\ar[rr]^{\psi_{\pm,n+1}^n}\ar[dr]^{F_n}&&\shg(-M_{T},-{\Ga}_{n+1})\ar[dl]_{F_{n+1}}\\
&\shg(-M_{T'},-\ga_{T'})&
}    
\end{equation*}
The proofs are also similar.
\erem

Combining Lemma \ref{lem_3: graded bypass triangle}, Lemma \ref{lem_3: iso at particular gradings}, and Lemma \ref{lem_3: surgery triangle}, we get the following result.

\begin{lem}[{\cite[Lemma 3.22]{LY2020}}]\label{lem_3: G_n is zero}
	For a large enough integer $n$, the map $G_n$ in Lemma \ref{lem_3: surgery triangle} is zero. Hence $F_{n+1}$ is surjective by the exact triangle (\ref{eq_6: surgery exact triangle}).
\end{lem}

The map $F_{n+1}$ in Lemma \ref{lem_3: surgery triangle} is a map associated to a contact 2-handle attachment. We have the following grading preserving result.

\blem[{\cite[Section 4.2]{LY2021}}]\label{lem: contact 2-handle preserve gradings}
Suppose $(M,\ga)$ is a balanced sutured manifold and $S\subset (M,\ga)$ is an admissible surface. Suppose $\al\subset M$ is a properly embedded arc that intersects $S$ transversely and $\partial \al\cap \partial S=\emptyset$. Let $N=M\backslash{\rm int}N(\al),S_N=S\cap N$, and let $\mu\subset \partial N$ be the meridian of $\al$ which is disjoint from $\partial S_N$. Suppose $\ga_N$ is a suture on $\partial N$ such that $(N,\ga_N)$ is a balanced sutured manifold and attaching a contact 2-handle along $\mu$ gives $(M,\ga)$. Let $C_{h^2,\mu}$ be the map associated to the contact 2-handle attachment. Then for any $i\in\mathbb{Z}$, we have\[C_{h^2,\mu}(\shg(-N,-\ga_N, S_N,i))\subset \shg(-M,-\ga,S,i).\] 
\elem

\subsection{One tangle component}\label{subsection: One tangle component}\quad

In this subsection, we apply lemmas in Section \ref{Subsection: Basic setups} to obtain a decomposition of twisted sutured homology associated to one tangle component. The results in this subsection are a generalization of \cite[Section 4.3]{LY2020}, where we dealt with rationally null-homologous knots. The proofs are almost identical, so we omit details and only point out the difference.

We adapt the notations in Subsection \ref{Subsection: Basic setups}. Suppose $(M,\ga)$ is a balanced sutured manifold and suppose $T\subset (M,\ga)$ is a vertical tangle with only one component $\al=T_1$, which is rationally null-homologous of order $q$. Let $M_T$ be the manifold obtained from $M$ by removing a neighborhood of $T$ and let $\ga_T=\ga\cup m_\al$, where $m_\al$ is a positively oriented meridian of $\al$. 

We start with the following lemma, which roughly says the summands in the `middle' gradings of $\shg(-M_T,-\Ga_n)$ associated to $S_n$ are cyclic of order $q$. 
\blem\label{lem_4: q-cyclic}
Suppose $n\in\mathbb{N}$ and $i_1,i_2\in \intg$ satisfying $i_1,i_2\in(\rho_n,P_n)$ and $i_1-i_2=q$, where $\rho_n$ and $P_n$ are constants in Lemma \ref{lem_3: iso at particular gradings}:

$$\rho_n=i^{n}_{max}-nq{\rm~and~} P_n=i_{min}^n+(n+1)q-\tau(+).$$
Then we have
$$\shg(-M_T,-\Ga_n,S_n,i_1)\cong \shg(-M_T,-\Ga_n,S_n,i_2).$$
\elem
\bpf
Based on Lemma \ref{lem_3: iso at particular gradings}, the proof is similar to that of \cite[Lemma 4.20]{LY2020}. Here we only include some key steps as follows.

Since $i_1< P_n$, by Lemma \ref{lem_3: iso at particular gradings}, we know
$$\shg(-M_T,-\Ga_n,S_n,i_1)\cong\shg(-M_T,-\Ga_{n+1},S_n,i_1-(i_{min}^n-i_{min}^{n+1})).$$
Similarly, since $i_2> \rho_n$, we know that
$$\shg(-M_T,-\Ga_n,S_n,i_2)\cong\shg(-M_T,-\Ga_{n+1},S_n,i_2+(i_{max}^{n+1}-i_{max}^n)).$$
By definitions of $i_{min}^n$ and $i_{max}^n$ in Remark \ref{rem: imin imax}, we have 
\beq
i_1-i_{min}^n+i_{min}^{n+1}=~&i_2-i_{max}^n+i_{max}^{n+1}+q+(i_{max}^n-i_{min}^n)-(i_{max}^{n+1}-i_{min}^{n+1})\\
=~&i_2-i_{max}^n+i_{max}^{n+1}+q+(-\chi(\bar{S}_{n})-\tau(n))-(-\chi(\bar{S}_{n+1})-\tau(n+1))\\
=~&i_2-i_{max}^n+i_{max}^{n+1}+q+(-\chi(\bar{S}_{+})+nq)-(-\chi(\bar{S}_{+})+(n+1)q)\\
=~&i_2-i_{max}^n+i_{max}^{n+1}.
\eeq
Hence we obtain the desired result.
\epf
Note that 
\beq
P_n-\rho_n=~&(i_{min}^n+(n+1)q-\tau(+))-(i^{n}_{max}-nq)\\
=~&-(i^{n}_{max}-i^{n}_{min})-\tau(+)+(2n+1)q\\
=~&-(-\chi(\bar{S}_+)+nq)-\tau(+)+(2n+1)q\\
=~&\chi(\bar{S}_+)-\tau(+)+(n+1)q.\\
\eeq
Thus, the difference of $P_n$ and $\rho_n$ can be infinitely large.
\bdefn\label{defn_4: essential component}
Define $Q_n=P_n-q+\tau(+)$. Suppose $n\in\mathbb{N}$ satisfies $Q_n-\rho_n>q$, define
$$\mathcal{SHG}_\al(-M,-\ga,i)\deq\shg(-M_T,-\Ga_n,S_n,Q_n-i),$$
and
$$\mathcal{SHG}_\al(-M,-\ga)\deq\bigoplus_{i=1}^{q}\mathcal{SHG}_\al(-M,-\ga,i).$$
\edefn
\brem\label{rem: sizes}
From the definitions of $Q_n,P_n,\rho_n$, and the fact $$i_{max}^n-i_{min}^n=-\chi(\bar{S}_n)+\tau(n)=-\chi(\bar{S}_+)+nq$$in Remark \ref{rem: imin imax}, we have
\begin{equation}\label{eq: difference}
    i_{max}^n-Q_n=i_{max}^n-(P_n-q+\tau(+))=i_{max}^n-i_{min}^n-nq=-\chi(\bar{S}_+)=\rho_n-i_{min}^n
\end{equation} 
Those equations are used in the proof of Lemma \ref{lem_4: the dimension equals}. This motivates the definition of $Q_n$, which is only for the convenience of the computation. The following remark implies we can freely choose $Q_n$ so that $Q_n-i\in(\rho_n,P_n)$ for $i=1,\dots,q$ to carry out the construction.
\erem
\brem\label{rem: grading shift}
From Lemma \ref{lem_3: iso at particular gradings} and the fact \[P_{n+1}-P_n=i_{min}^{n+1}-i_{min}^n+q=i_{max}^{n+1}-i_{max}^n,\]the isomorphism class of $\mathcal{SHG}_\al(-M,-\ga,i)$ is independent of the choice of the large integer $n$. Also, by Lemma \ref{lem_4: q-cyclic}, the isomorphism class of $\mathcal{SHG}_\al(-M,-\ga)$ would be the same (up to a $\mathbb{Z}_q$ grading shift) if we consider arbitrary $q$ many consecutive gradings within the range $(\rho_n,P_n)$. 
\erem
\brem\label{rem: case knot}
For a rationally null-homologous knot $\widehat{K}\subset \widehat{Y}$ with a basepoint $p$, we can remove a neighborhood of $p$ add a suture $\delta$ on $\partial N(p)$ such that two intersection points of $\widehat{K}$ and $\partial N(p)$ lie on $R_+(\ga)$ and $R_-
(\ga)$, respectively. Then $\widehat{K}$ becomes a vertical tangle $\al$ in $(\widehat{Y}-{\rm int}N(p),\delta)$ which is rationally null-homologous. In this case, $\mathcal{SHG}_\al(\widehat{Y}-{\rm int}N(p),\delta,i)$ reduces to $\mathcal{I}_+(-\widehat{Y},\widehat{K},i)$ in \cite[Definition 4.21]{LY2020}, up to a $\mathbb{Z}_q$ grading shift.
\erem

\blem\label{lem_4: the dimension equals}
Let $\mathcal{SHG}_\al(-M,-\ga)$ be defined as in Definition \ref{defn_4: essential component}. We have
$$\dim_{\mathbb{F}}\mathcal{SHG}_\al(-M,-\ga)=\dim_{\mathbb{F}}\shg(-M,-\ga).$$
\elem
\bpf
Based on Lemma \ref{lem_3: graded bypass triangle}, the proof is similar to that of \cite[Lemma 4.25]{LY2020}. Now we split the bypass exact triangle of $(-\Ga_+,-\Ga_n,-\Ga_{n+1})$ into five blocks of sizes \[q, -\chi(\bar{S}_+)+1,\chi(\bar{S}_+)+(n-1)q-1,q,-\chi(\bar{S}_+)+1,\]respectively, and split the bypass exact triangle of $(-\Ga_-,-\Ga_n,-\Ga_{n+1})$ into five blocks of sizes\[ -\chi(\bar{S}_+)+1,q,\chi(\bar{S}_+)+(n-1)q-1,-\chi(\bar{S}_+)+1,q,\]respectively. Remark \ref{rem: sizes} ensures that the proof of \cite[Lemma 4.25]{LY2020} applies verbatim. Here we only include some key steps as follows.

Suppose $n\in\mathbb{N}$ satisfies $Q_n-\rho_n>q$. We can apply Proposition \ref{lem_3: graded bypass triangle}. Using blocks, we have the following. (There is no enough room for writing down the whole notation for formal sutured homology, so we will only write down the sutures to denote them.)
\begin{equation*}
\xymatrix@R=0.2ex{
{\rm size}&-\Ga_{+}\ar[rr]^{\psi_{+,n}^{+}}&&-\Ga_n\ar[rr]^{\psi_{+,n+1}^{n}}&&-\Ga_{n+1}\ar[rr]^{\psi_{+,+}^{n+1}}&&-\Ga_{+}\\
q&G&&&&X_1&&G\\
-\chi(\bar{S}_+)+1&H&&A&&X_2&&H\\
\chi(\bar{S}_+)+(n-1)q-1&&&E&&X_3&&\\
q&&&F&&X_4&&\\
-\chi(\bar{S}_+)+1&&&D&&X_5&&\\
}    
\end{equation*}
The empty block implies the summands in the block are zeros. Note that \[i_{max}^+-i_{min}^++1=-\chi(\bar{S}_+)+\tau(+)+1\le q+(-\chi(\bar{S}_+)+1).\]From the exactness, we know that
$$X_1=G,~X_3=E,~X_4=F,{\rm~and~}X_5=D.$$

There is another bypass exact triangle, and similarly we have
\begin{equation*}
\xymatrix@R=0.2ex{
{\rm size}&-\Ga_{-}\ar[rr]^{\psi_{-,n}^{-}}&&-\Ga_n\ar[rr]^{\psi_{-,n+1}^{n}}&&-\Ga_{n+1}\ar[rr]^{\psi_{-,-}^{n+1}}&&-\Ga_{-}\\
-\chi(\bar{S}_+)+1&&&A&&A&&\\
q&&&B&&B&&\\
\chi(\bar{S}_+)+(n-1)q-1&&&C&&C&&\\
-\chi(\bar{S}_+)+1&I&&D&&X_6&&I\\
q&J&&&&J&&J\\
}
\end{equation*}
Note that \[i_{max}^--i_{min}^-+1=-\chi(\bar{S}_+)-q+1\le q+(-\chi(\bar{S}_+)+1).\]
Comparing the two expressions of $\shg(-M_T,-\Ga_{n+1},S_n)$, we have
\begin{equation*}
\left(
\begin{array}{c}
    G\\
    X_2\\
    E\\
    F\\
    D
\end{array}
\right)=\shg(-M_T,-\Ga_{n+1},S_n)=\left(
\begin{array}{c}
    A\\
    B\\
    C\\
    X_6\\
    J
\end{array}
\right).
\end{equation*}
Taking sizes into consideration, we know that
\begin{equation*}
\left(
\begin{array}{c}
    G\\
    X_2
\end{array}
\right)=\left(
\begin{array}{c}
    A\\
    B
\end{array}
\right),~E=C,{\rm~and~}\left(
\begin{array}{c}
    F\\
    D
\end{array}
\right)=\left(
\begin{array}{c}
    X_6\\
    J
\end{array}
\right).
\end{equation*}
Thus, we know that
\begin{equation*}
\shg(-M_T,-\Ga_{n+1},S_n)=\left(
\begin{array}{c}
    A\\
    B\\
    E\\
    F\\
    D
\end{array}
\right).
\end{equation*}
By construction, we have
\beq
\dim_{\mathbb{F}}\mathcal{SHG}_\al(-M,-\ga)=&\dim_{\mathbb{F}}B\\
=&\dim_{\mathbb{F}}\shg(-M_T,-\Ga_{n+1})-\dim_{\mathbb{F}}\shg(-M_T,-\Ga_n)\\
=&\dim_{\mathbb{F}}\shg(-M,-\ga).
\eeq
Note that the last equality follows from Lemma \ref{lem_3: G_n is zero}. Hence we obtain the desired result.
\epf

\brem
The essential difference for the case of tangles is that $\Ga_+$ is not equal to $\Ga_-$, though it is true in the case of knots in Remark \ref{rem: case knot}. 
\erem
\bprop\label{prop_4: F_n is an isomorphism}
Suppose $n\in\mathbb{N}$ is large enough. Then the map $F_n$ in Lemma \ref{lem_3: surgery triangle} restricted to $\mathcal{SHG}_\al(-M,-\ga)$ is an isomorphism, \textit{i.e.}
$$F_n|_{\mathcal{SHG}_\al(-M,-\ga)}:\mathcal{SHG}_\al(-M,-\ga)\xra{\cong} \shg(-M,-\ga).$$
\eprop

\bpf

Based on Lemma \ref{lem_3: G_n is zero}, the proof is similar to that of \cite[Proposition 4.26]{LY2020}. It suffices to show that the restriction of $F_n$ is surjective. Here we only include some key steps as follows.

By Lemma \ref{lem_3: G_n is zero}, we know that $F_n$ is surjective. Then it suffices to show that $F_n$ remains surjective when restricted to $\mathcal{SHG}_\al(-M,-\ga)$. For any $x\in \shg(-M,-\ga)$, let $y\in\shg(-M_T,-\Ga_n)$ be an element so that $F_n(y)=x.$ Suppose
$$y=\sum_{j\in\intg}y_j,\text{ where }y_j\in\shg(-M_T,-\Ga_n,S_n,j).$$
For any $y_j$, we want to find $y'_j\in\mathcal{SHG}_\al(-M,-\ga)$ so that $F_n(y_j)=F_n(y_j').$

To do this, we first assume that $j\ge Q_n$. Then there exists an integer $m$ so that
$$Q_n-q\leq j-mq\leq Q_n-1.$$
We can take
\begin{equation}\label{eq: many bypasses}
    y_j'=(\psi_{-,n+1}^{n,j-mq})^{-1}\circ\cdots\circ(\psi_{-,n+m}^{n,i_{max}^{n+m}-i_{max}^n+j-mq})^{-1}\circ\psi_{+,n+m}^{n+m-1}\circ\dots\circ\psi_{+,n+1}^n(y_j).
\end{equation}
From Lemma \ref{lem_3: iso at particular gradings}, all the negative bypass maps involved in (\ref{eq: many bypasses}) are isomorphisms so the inverses exist. Also, we have
$$y_j'\in\shg(-M_T,-\Ga_n,S_n,j-mq)\subset \mathcal{SHG}_\al(-M,-\ga).$$
Finally, from commutative diagrams in Lemma \ref{lem_3: surgery triangle}, we know that $F_n(y'_j)=F_n(y_j).$ 

For $$j\in [ Q_n-q, Q_n-1],$$we can simply take $y_j'=y_j$. 

For $j<Q_n-q$, we can pick $y_j'$ similarly, while switching the roles of $\psi_{+,*}^*$ and $\psi_{-,*}^*$ in (\ref{eq: many bypasses}). 

In summary, we can take
$$y'=\sum_{j\in\intg}y_j'\in\mathcal{SHG}_\al(-M,-\ga)\text{ with }F_n(y')=F_n(y)=x.$$
Hence $F_n$ is surjective, and we obtain the desired result.
\epf

\brem
In Definition \ref{defn_4: essential component}, we use a large enough integer $n$ to define $\mathcal{SHG}_\al(-M,-\ga)$. We can also define $\Ga_{-n}$ from $\Ga_0$ by twisting along $\ga_1$ for $n$ times. For a large enough integer $n$, we can define a vector space $\mathcal{SHG}_\al^\p(-M,-\ga)$ generalizing $\mathcal{I}_-(-\widehat{Y},\widehat{K})$ in \cite[Definition 4.27]{LY2020}. However, from the discussion in \cite[Section 4.4 in ArXiv version 2]{LY2020} between $\mathcal{I}_+(-\widehat{Y},\widehat{K})$ and $\mathcal{I}_-(-\widehat{Y},\widehat{K})$, we expect that $\mathcal{SHG}_\al^\p(-M,-\ga)$ is isomorphic to $\mathcal{SHG}_\al(-M,-\ga)$ up to a $\mathbb{Z}_q$ grading shift. Hence there is no new information and we skip the discussion here. 
\erem

\subsection{More tangle components}\quad

In this subsection, we obtain a decomposition of twisted sutured homology associated to more tangle components. Suppose $(M,\ga)$ is a balanced sutured manifold. For a vertical tangle $T$ in $M$, let $M_T=M\backslash{\rm int}N(T)$ and let $\ga_T$ be the union of $\ga$ and positively oriented meridians of components of $T$. 

First, we prove some lemmas about homology groups.
\blem\label{lem: rk}
For any connected tangle $\al$ in $M$, we have\[{\rm rk}_{\mathbb{Z}}H_1(M_\al)=\begin{cases}
{\rm rk}_{\mathbb{Z}}H_1(M) & \text{if }[\al]\neq 0\in H_1(M,\partial M;\mathbb{Q}),\\{\rm rk}_{\mathbb{Z}}H_1(M)+1&\text{if }[\al]= 0\in H_1(M,\partial M;\mathbb{Q}).
\end{cases}\]
\elem
\bpf
Consider the long exact sequence assoicated to the pair $(M,M_\al)$:
\begin{equation}
\label{eq: long exact sequence}
H^1(M,M_\al)\xra{p^*_1} H^1(M)\xra{i^*_1} H^1(M_\al)\xra{\delta^*_1}H^2(M,M_\al)\xra{p^*_2} H^2(M)\xra{i^*_2} H^2(M_\al)\xra{\delta^*_2}H^3(M,M_\al).
\end{equation}
By the excision theorem, we have \[H^*(M,M_\al)\cong H^j(N(\al),\partial N(\al)\cap M_\al)\cong H^j(D^2,\partial D^2)\cong \begin{cases}\mathbb{Z}&j=2,\\0&j=1,3.\end{cases}\]Since $H^2(N(\al),\partial N(\al)\cap M_\al)$ is generated by the disk that is  the Poincar\'{e} dual of $[\al\cap N(\al)]$ and $p^*_2$ is induced by the projection, the image of $p^*_2$ is generated by the Poincar\'{e} dual of $[\al]$. Since $H^1(M)$ and $H_1(M)$ always have the same rank, we obtain the rank equation from (\ref{eq: long exact sequence}).
\epf

\blem\label{lem: torsion free}
Suppose $(M,\ga)$ is a balanced sutured manifold. There exists a (possibly empty) tangle $T=T_1\cup\cdots \cup T_m$ in $M$, such that ${\rm Tor}(H_1(M_T))=0$ and for any $T^\p\subset T$ and $T_i\subset T\backslash T^\p$, we have\begin{equation}
    \label{eq: all rational tangle}
[T_i]=0\in H_1(M_{T^\p},\partial M_{T^\p} ;\mathbb{Q}).
\end{equation}
\elem
\bpf
Suppose $\al$ is a connected tangle in $M$. From (\ref{eq: long exact sequence}) and the proof of Lemma \ref{lem: rk}, we have\[\mathbb{Z}\langle \phi_\al\rangle\xra{p_2^*}H^2(M)\xra{i_2^*}H^2(M_\al)\ra 0,\]
where $\phi_\al$ is the Poincar\'{e} dual of $[\al]$. By the universal coefficient theorem, the torsion subgroups of $H^2(M)$ and $H_1(M)$ are isomorphic. In particular, ${\rm Tor}(H^2(M))=0$ if and only if ${\rm Tor}(H_1(M))=0$. Let $\al$ be a rationally null-homologous tangle, then \[{\rm Tor}(H^2(M_\al))\cong {\rm Tor}(H^2(M))/{\rm PD}(\al).\]Thus, we can always choose connected tangles \[T_1\subset M,T_2\subset M_{T_1},T_3\subset M_{T_1\cup T_2},\dots,T_m\subset M_{T_1\cup\cdots\cup T_{m-1}}\]that are rationally null-homologous to kill the whole torsion subgroup. In other word, for $T=T_1\cup\cdots\cup T_m$, we have ${\rm Tor}(H_1(M_T))=0$.

By Lemma \ref{lem: rk}, we have\begin{equation}
    \label{eq: rk equal 2}{\rm rk}_\mathbb{Z}H_1(M_T)={\rm rk}_\mathbb{Z}H_1(M)+m.
\end{equation}Hence for any $T^\p$ and any $T_i$ satisfy the assumption, (\ref{eq: all rational tangle}) holds, otherwise it contradicts with the rank equality (\ref{eq: rk equal 2}).
\epf

\brem\label{rem: vertical tangle}
Since moving the endpoints of a tangle on the boundary of the ambient 3-manifold does not change the homology class of the tangle, we can suppose the tangle $T$ in Lemma \ref{lem: torsion free} is a vertical tangle. Moreover, when $M$ has connected boundary, we can suppose endpoints of $T$ all lie in a neighborhood of a point on the suture $\ga$.
\erem

\blem\label{lem: generate}
Suppose $(M,\ga)$ is a balanced sutured manifold and suppose $\al$ is a connected rationally null-homologous tangle of order $q$. Let $S_{\al}$ be a Seifert surface of $T_i$, \textit{i.e.}, $\partial S_{i}$ consists of $q$ parallel copies of $\al$ and arcs on $\partial M$. Suppose $S_{1},\dots,S_{n}$ are admissible surfaces in $(M,\ga)$ generating $H_2(M,\partial M)$. Then the restrictions of $S_1,\dots,S_n$ and $S_\al$ on $M_T$ generate $H_2(M_T,\partial M_T)$.
\elem
\bpf
From (\ref{eq: long exact sequence}) and the proof of Lemma \ref{lem: rk}, we have
\[0\ra H^1(M)\xra{i_1^*} H^1(M_\al)\xra{\delta^*_1} \mathbb{Z}\langle \phi_\al \rangle\xra{p_2^*}H^2(M),\]
where $\phi_\al$ is the Poincar\'{e} dual of $[\al]$. It is straighforward to calculate
\begin{equation}
    \label{eq: a}
\delta_1^*({\rm PD}([S_\al]))=q\phi_\al.
\end{equation}Since $H^1(M)\cong H_2(M,\partial M)$, we have
\begin{equation}
\label{eq: b}\begin{aligned}
H_2(M_\al,\partial M_\al)/H_2(M,\partial M)\cong& H^1(M_\al)/H^1(M)\cong H^1(M_\al)/{\rm im}(i_1^*)\\\cong H^1(M_\al)/{\rm ker}(\delta_1^*)\cong& {\rm im}(\delta_1^*)\cong {\rm ker}(p_2^*).
\end{aligned}
\end{equation}Since the image of $p_2^*$ is the Poincar\'{e} dual of $[\al]$, we have 
\begin{equation}
    \label{eq: c}{\rm ker}(p_2^*)\cong \langle q\phi_\al\rangle.
\end{equation}Combining (\ref{eq: a}), (\ref{eq: b}), and (\ref{eq: c}), we know that $[S_\al]$ generates $H_2(M_\al,\partial M_\al)/H_2(M,\partial M)$. Thus, we conclude the desired property.

\epf



In the rest of this subsection, we suppose $(M,\ga)$ is a balanced sutured manifold and $T=T_1\cup\cdots \cup T_m$ is a vertical tangle satisfying Lemma \ref{lem: torsion free}. Suppose the order of the first component $T_1$ in $H_1(M)$ is $q_1$ and suppose $S_1$ is a Seifert surface of $T_1$. 

\begin{conv}
We will still use $S_1$ to denote its restriction on $M_{T_1}$. This also applies to other Seifert surfaces mentioned below.
\end{conv}
We adapt the construction in Subsection \ref{Subsection: Basic setups}. Applying results in Subsection \ref{subsection: One tangle component}, we have\[\mathcal{SHG}_{T_1}(-M,-\ga)\deq \bigoplus_{i=1}^{q_1}\shg(-M_{T_1},-\Ga_{n},(S_1)_{n},Q_n-i)\cong \shg(-M,-\ga),\]where $n$ is a large integer, $(S_1)_n$ is a (possibly empty) stabilization of $S_1$, and $Q_n$ is a fixed integer. For simplicity, we choose a large integer $n_1$ such that $(S_1)_{n_1}=S_1$ and write \[\Ga_{n_1}^1=\Ga_{n}|_{n=n_1}\aand Q_{n_1}^1=Q_{n}|_{n=n_1}.\]

For the second component $T_2$, suppose $S_2$ is its Seifert surface in $M_{T_1}$ with $\partial S^2$ containing $q_2$ copies of $T_2$. Now we can apply the construction in Subsection \ref{Subsection: Basic setups} and the results in Subsection \ref{subsection: One tangle component} to $(M,\Ga_{n_1}^1)$. For a large integer $n_2$ such that $(S_2)_{n_1}=S_2$, we define
\[\mathcal{SHG}_{T_1\cup T_2}(-M,-\ga)\deq \bigoplus_{i_1=1}^{q_1}\bigoplus_{i_2=1}^{q_2}\shg(-M_{T_1\cup T_2},-\Ga_{n_2}^2,(S_1,S_2),(Q_{n_1}^1-i_1,Q_{n_2}^2-i_2))\cong \shg(-M,-\ga).\]
Iterating this procedure, we have the following definition.

\bdefn \label{defn: decomposition}
For $i=1,\dots,m$, suppose the component $T_k$ is rationally null-homologous of order $q_k$ in $M_{T_1\cup\cdots\cup T_{k-1}}$. Inductively, for $k=1,\dots,m$, we choose a large integer $n_k$, a suture $\Ga_{n_k}^k\subset \partial M_{T_1\cup\cdots\cup T_{k}}$, a Seifert surface $S_k=(S_k)_{n_k}\subset M_{T_1\cup\cdots\cup T_{k}}$, and an integers $Q_{n_k}^k$, such that $n_{k},\Ga_{n_k}^k,S_k,Q_{n_k}^k$ depend on the choices for the first $(k-1)$ tangles. Define\[\mathcal{SHG}_T(-M,-\ga)\deq \bigoplus_{i_1\in[1,q_1],\dots,i_m\in[1,q_m]}\shg(-M_T,-\Ga_{n_m}^m,(S_1,\dots,S_m),(Q_{n_1}^1-i_1,\cdots,Q_{n_m}^m-i_m)).\]
\edefn
\brem\label{rem: not canonical}
Though we only use the subscript $T$ in the notation $\mathcal{SHG}_T(-M,-\ga)$, it is not known if $\mathcal{SHG}_T(-M,-\ga)$ is independent of the choices of all constructions. In particular, we have to choose an order of the components to define $\mathcal{SHG}_T(-M,-\ga)$.
\erem
Applying results in Subsection \ref{subsection: One tangle component} for $m$ times, the following proposition is straightforward.
\bprop\label{prop: isomorphism}
$\mathcal{SHG}_T(-M,-\ga)\cong \shg(-M,-\ga).$
\eprop

The map $H_1(M_{T_1})\to H_1(M)$ is surjective. The $q_1$ direct summands of $\shg_{T_1}(-M,-\ga)$ correspond to the order $q_1$ torsion subgroup generated by\[[T_1]\in {\rm Tor}(H_1(M,\partial M))\cong {\rm Tor}(H^2(M))\cong {\rm Tor}(H_2(M))\]Hence the summands of $\shg_{T_1}(-M,-\ga)$ provide a decomposition of $\shg(-M,-\ga)$ with respect to the torsion subgroup generated by $[T_1]$. By induction and the fact that ${\rm Tor}(H_1(M_T))=0$, we can regard summands in $\mathcal{SHG}_T(-M,-\ga)$ as a decomposition of $\shg(-M,-\ga)$ with respect to ${\rm Tor}(H_1(M))$.

To provide a decomposition of $\shg(-M,-\ga)$ with respect to the whole $H_1(M)$ as in Theorem \ref{thm: main}, we choose admissible surfaces $S_{m+1},\dots,S_{m+n}$ generating $H_2(M,\partial M)$. By Lemma \ref{lem: generate}, the restrictions of $S_1,\dots,S_{m+n}$ generate $H_2(M_T,\partial M_T)$. By Lemma \ref{lem: contact 2-handle preserve gradings}, the gradings associated to these surfaces behave well under restriction.

\bdefn\label{defn: enhanced euler}
Consider the construction as above. For $i=1,\dots,m+n$, let $\rho_1,\dots,\rho_{m+n}\in H_1(M_T)=H_1(M_T)/\tor$ be the class satisfying $\rho_i\cdot S_j=\delta_{i,j}$. Consider$$j_*:\mathbb{Z}[H_1(M_T)]\to \mathbb{Z}[H_1(M)].$$We write \[H=H_1(M),\bs{S}=(S_1,\dots,S_{m+n}),-i^\p_{k}=Q^k_{n_k}-i_{n+k}\text{ for }k=1,\dots,m,\]and\[-\bs{i}^\p=(-i_1^\p,\dots,-i_m^\p,-i_{m+1},\dots,-i_{m+n}),\bs{\rho}^{-\bs{i}^\p}=\rho_1^{-i_1^\p}\cdots \rho_n^{-i_m^\p}\cdot \rho_{m+1}^{-i_{n+1}}\cdots \rho_{m+n}^{-i_{m+n}}.\]The \textbf{enhanced Euler characteristic} of $\shg(-M,-\ga)$ is \beq\en(\shg&(-M,-\ga))=j_*(\chi(\mathcal{SHG}_T(-M,-\ga))) \\&\deq j_*(\sum_{\substack{i_{1}\in[1,q_1],\dots,i_{m}\in[1,q_m]\\(i_{m+1},\dots,i_{m+n})\in\intg^n}}\chi(\shg(-M_T,-\ga_T,\bs{S},-\bs{i}^\p))\cdot \bs{\rho}^{-\bs{i}^\p})\in \mathbb{Z}[H]/\pm H.\eeq
For $h\in H_1(M)$, let $\shg(-M,-\ga,h)$ be image of the summand of $\mathcal{SHG}_T(-M,-\ga)$ under the isomorphism in Propsition \ref{prop: isomorphism} whose corresponding element in $\en(\shg(-M,-\ga))$ is $h$.
\edefn
\brem
As mentioned in Remark \ref{rem: not canonical}, the definition of $\shg(-M,-\ga,h)$ is not canonical, \textit{i.e.} it may depend on many auxiliary choices. After fixing these choices, it is still only well-defined up to a global grading shift by multiplication by an element in $h_0\in H_1(M)$. However, by Theorem \ref{thm: chi equal, main}, the enhanced Euler characteristic $\en(\shg(-M,-\ga))$ only depends on $(M,\ga)$.
\erem

\section{Sutured Heegaard Floer homology}\label{section: Sutured Heegaard Floer homology}

In this section, we discuss properties of sutured (Heegaard) Floer homology $SFH$ that are similar to those for twisted sutured homology, so we can apply results in Section \ref{Section: Decomposition associated to tangles} to $SFH$. Since $SFH$ is not defined by closures of sutured manifolds, the maps associated to surgeries and contact handle attachments are different from those for $\shg$.

\subsection{Construction and gradings}\label{subsection: sfh gradings}\quad

In this subsection, we describe the definition of $SFH$ and discuss the gradings on $SFH$ associated to admissible surfaces.

\bdefn[{\cite[Section 2]{juhasz2006holomorphic}}]
A \textbf{balanced diagram} $\mch=(\Sigma,\al,\be)$ is a tuple  satisfying the following.
\begin{enumerate}[(1)]
    \item $\Sigma$ is a compact, oriented surface with boundary.
    \item $\al = \{ \al_1, \dots, \al_n \}$ and $\be = \{ \be_1,\dots, \be_n \}$ are two sets of pairwise disjoint simple closed curves in the interior of $\Sigma$.
    \item The maps $\pi_0(\partial \Sigma)\to \pi_0(\Sigma\backslash \al)$ and $\pi_0(\partial \Sigma)\to \pi_0(\Sigma\backslash\be)$ are surjective.
\end{enumerate}

For such triple, let $N$ be the 3-manifold obtained from $\Sigma\times [-1,1]$ by attaching 3–dimensional 2–handles along $\al_i \times \{-1\}$ and $\be_i \times \{1\}$ for $i=1,\dots,n$ and let $\nu=\partial \Sigma\times \{0\}$. A balanced diagram $(\Sigma,\al,\be)$ is called \textbf{compatible} with a balanced sutured manifold $(M,\ga)$ if the balanced sutured manifold $(N,\nu)$ is diffeomorphic to $(M,\ga)$. 

\edefn
Suppose $\mch=(\Sigma,\al,\be)$ is a balanced diagram with $g=g(\Sigma)$ and $n=|\al|=|\be|$. 

\begin{conv}
In this paper, we always suppose balanced diagrams satisfy the admissible condition in \cite[Section 3]{juhasz2006holomorphic}.
\end{conv}

Consider two tori$$\mathbb{T}_\al\deq \al_1\times\cdots\times\al_n~{\rm and}~\mathbb{T}_\be\deq \be_1\times\cdots\times\be_n$$in the symmetric product $${\rm Sym}^n\Sigma\deq (\prod_{i=1}^n\Sigma)/S_n.$$

The chain complex $SFC(\mch)$ is a free $\ft$-module generated by intersection points $\bs{x}\in\mathbb{T}_\al\cap \mathbb{T}_\be$. Let $\pi_2(\bs{x},\bs{y})$ be the set of homology classes of Whitney disks connecting intersection points $\bs{x}$ and $\bs{y}$. Choose a generic path of almost complex structures $J_s$ on $\sym^n\Sigma$. For $\phi\in \pi_2(\bs{x},\bs{y})$, let $\mathcal{M}_{J_s}(\phi)$ be the moduli space of $J_s$-holomorphic maps \[u:[0,1]\times \mathbb{R}\to \sym^n\Sigma\]which represent $\phi$ and let $\mu(\phi)$ be the expected dimension of $\mathcal{M}_{J_s}(\phi)$. The moduli space $\mathcal{M}_{J_s}(\phi)$ has a natural action of $\mathbb{R}$, corresponding to reparametrization of the source. We write
\[\widehat{\mathcal{M}}_{J_s}(\phi)\deq \mathcal{M}_{J_s}(\psi)/\mathbb{R}.\]

Based on the above construction, Juh\'{a}sz \cite{juhasz2006holomorphic} defined a differential on $SFC(\mch)$ by$$\partial_{J_s}(\bs{x})=\sum_{\bs{y}\in \mathbb{T}_\al\cap \mathbb{T}_\be}\sum_{\substack{\phi\in\pi_2(\bs{x},\bs{y})\\\mu(\phi)=1}}\#\widehat{\mathcal{M}}_{J_s}(\phi)\cdot \bs{y}.$$


\bthm[{\cite{juhasz2006holomorphic,Juhasz2012}}]\label{thm: SFH}
Suppose $(M,\ga)$ is a balanced sutured manifold. Then there is an admissible balanced diagram $\mch$ compatible with $(M,\ga)$. The vector spaces $H(SFC(\mch),\partial_{J_s})$ for different choices of $\mch$ and $J_s$, together with some canonical maps, form a transitive system $SFH(M,\ga)$ over $\mathbb{F}_2$.
\ethm
For a balanced sutured manifold $(M,\ga)$, we can decompose $SFH(M,\ga)$ along spin$^c$ structures.

Fix a Riemannian metric $g$ on $M$. Let $v_0$ be a nowhere vanishing vector field along $\partial M$ that points into $M$ along $R_-(\ga)$, points out of $M$ along $R_+(\ga)$, and on $\ga$ it is the gradient of the height function $A(\ga)\times I\to I$. The space of such vector fields is contractible, so the choice of $v_0$ is not important.

Suppose $v$ and $w$ are nowhere vanishing vector fields on $M$ that agree with $v_0$ on $\partial M$. They are called \textbf{homologous} if there is an open ball $B\subset{\rm int}M$ such that $v$ and $w$ are homotopic on $M\backslash B$ through nowhere vanishing vector
fields ${\rm rel}~\partial M$. Let $\spin(M,\ga)$ be the set of homology classes of nowhere vanishing vector fields $v$ on $M$ with $v|_{\partial M}=v_0$. Note that ${\rm Spin}^c(M,\ga)$ is an affine space over $H^2(M,\partial M)$.

Suppose $\mch=(\Sigma,\al,\be)$ is a balanced diagram compatible with $(M,\ga)$. For each intersection point $\bs{x}\in\mathbb{T}_\al\cap\mathbb{T}_\be$, we can assign a spin$^c$ structure $\mathfrak{s}(\bs{x})\in{\rm Spin}^c(M,\ga)$ as follows (\textit{c.f.} \cite[Section 4]{juhasz2006holomorphic}).

we choose a self-indexing Morse function $f:M\to [-1,4]$ such that \[f^{-1}(\frac{3}{2})=\Sigma\times\{0\}.\] Moreover, curves $\al,\be$ are intersections of $\Sigma\times\{0\}$ with the ascending and descending manifolds of the index 1 and 2 critical points of $f$, respectively. Then any intersection point of $\al_i\subset \al$ and $\be_j\subset\be$ corresponds to a trajectory of ${\rm grad}f$ connecting a index 1 critical point to a index 2 critical point. For $\bs{x}\in\mathbb{T}_\al\cap\mathbb{T}_\be$, let $\ga_{\bs{x}}$ be the multi-trajectory corresponding to intersection points in $\bs{x}$.

In a neighborhood $N(\ga_{\bs{x}})$, we can modify ${\rm grad}f$ to obtain a nowhere vanishing vector field $v$ on $M$ such that
$v|_{\partial M}=v_0$. Let $\mathfrak{s}(\bs{x})\in\spin(M,\ga)$ be the homology class of this vector field $v$.

From the assignment of the spin$^c$ structure, we have the following proposition.
\bprop\label{prop: one cycle}
For any $\bs{x},\bs{y}\in \mathbb{T}_\al\cap\mathbb{T}_\be$, we have\[\mathfrak{s}(\bs{x})-\mathfrak{s}(\bs{y})={\rm PD}([\ga_{\bs{x}}-\ga_{\bs{y}}]),\]where ${\rm PD}:H_1(M)\to H^2(M,\partial M)$ is the Poincar\'{e} duality map.
\eprop
It can be shown that there is no differential between generators corresponding to different spin$^c$ structures. Hence we have the following decomposition.
\bprop[\cite{juhasz2006holomorphic}]
 For any balanced sutured manifold $(M,\ga)$, there is a decomposition$$SFH(M,\ga)=\bigoplus_{\mathfrak{s}\in \spin(M,\partial M)}SFH(M,\ga,\mathfrak{s}).$$
\eprop
Suppose $S\subset(M,\ga)$ is an admissible surface $S$. To associate a $\mathbb{Z}$-grading on $SFH(M,\ga)$ similar to Subsection \ref{subsec: gradings on SH}, we need to suppose $(M,\ga)$ is \textbf{strongly balanced}, \textit{i.e.} for every component $F$ of $\partial M$, we have\[\chi(F\cap R_+(\ga))=\chi(F\cap R_-(\ga)).\]

\brem\label{rem: strongly balanced}
If $\partial M$ is connected, then it is automatically strongly balanced. For any balanced sutured manifold $(M,\ga)$, we can obtain a strongly balanced manifold $(M^\p,\ga^\p)$ by attaching contact 1-handles \cite[Remark 3.6]{juhasz2008floer}. In Subsection \ref{subsection: Contact handles and bypasses SFH}, we will show\[SFH(M^\p,\ga^\p)\cong SFH(M,\ga)\]and this isomorphism respects spin$^c$ structures. Hence we can always deal with a strongly balanced manifold without losing any information.
\erem
\begin{conv}
When discussing the $\mathbb{Z}$-grading on $SFH(M,\ga)$ associated to an admissible surface $S\subset (M,\ga)$, we always suppose $(M,\ga)$ is strongly balanced.
\end{conv}

The following construction is based on \cite[Section 3]{juhasz2008floer}.

Let $v_0^\perp$ be the plane bundle perpendicular to $v_0$ under the fixing Riemannian metric $g$. Suppose $v_0^\perp$ is oriented so that $v_0$ is the positive normal vector field. The strongly balanced condition on $(M,\ga)$ ensures that $v_0^\perp$ is trivial (\textit{c.f.} \cite[Proposition 3.4]{juhasz2008floer}). Let $t$ be a trivialization of $v_0^\perp$. Since any spin$^c$ structure $\mathfrak{s}\in \spin(M,\ga)$ can be represented by a nonvanishing vector field $v$ on $M$ with $v|_{\partial M}=v_0$, we can define \[c_1(\mathfrak{s},t)\deq c_1(v^\perp,t)\in H^2(M,\partial M)\]to be the relative Euler class of the plane bundle $v^\perp$ with respect to the trivialization $t$, where $v^\perp$ is perpendicular to $v$. In other words, the class $c_1(\mathfrak{s},t)$ is the obstruction to extending $t$ from $\partial M$ to a trivialization of $v^\perp$ over $M$.

Let $v_S$ be the positive unit normal field of $S$. For a generic $S$, we can suppose $v_S$ is nowhere parallel to $v_0$ along $\partial S$. Let $p(v_S)$ be the projection of $v_S$ into $v^\perp_0$. Note that $p(v_S)|_{\partial S}$ is nowhere zero. Suppose the components of $\partial S$ are $T_1,\dots,T_k$, oriented by the boundary orientation. 

For $i=1,\dots,k$, let $r(T_i,t)$ be the rotation number $p(v_S)|_{T_i}$ with respect to the
trivialization $t$ as we go around $T_i$. Moreover, define\[r(S,t)\deq \sum_{i=1}^kr(T_i,t).\]Suppose $T_1,\dots,T_k$ intersect $\ga$ transversely. Define\begin{equation}\label{eq: defn c}
    c(S,t)=\chi(S)-\frac{1}{2}|\partial S\cap \ga|-r(S,t).
\end{equation}
\brem
The original definition of $c(S,t)$ in \cite[Section 3]{juhasz2008floer} involves the index $I(S)$, which is equal to $\frac{1}{2}|\partial S\cap \ga|$ when $T_1,\dots,T_k$ intersect $\ga$ transversely (\textit{c.f.} \cite[Lemma 3.9]{juhasz2008floer}).
\erem
Suppose $t_S$ is the trivialization of $v_0^\perp$ induced by $p(v_S)|_{\partial S}$. Then for any $v^\perp$ with $v^\perp|_{\partial M}=v^\perp_0$ and any trivialization $t$ of $v_0^\perp$, we have \begin{equation}\label{eq: rst}
    \langle c_1(v^\perp,t_S)-c_1(v^\perp,t),[S]\rangle=r(S,t)
\end{equation}
(\textit{c.f.} Proof of \cite[Lemma 3.10]{juhasz2008floer}; see also \cite[Lemma 3.11]{juhasz2010polytope}). In particular, we have $r(S,t_S)=0$.
\bdefn\label{defn: SFH grading}
 Consider the construction as above. Define 
 \begin{equation}
    SFH(M,\ga,S,i)\deq\bigoplus_{\substack{\mathfrak{s}\in\spin(M,\ga)\\\langle c_1(\mathfrak{s},t_S),[S]\rangle=-2i}}SFH(M,\ga,\mathfrak{s}).
\end{equation}
\edefn
\brem
The minus sign of $(2i)$ is to make this definition parallel to the $\mathbb{Z}$-grading on $\shg(M,\ga)$ associated to $S$. See the proofs of the following propositions.
\erem
\bprop
The decomposition in Definition \ref{defn: SFH grading} satisfies properties in Theorem \ref{thm_2: grading in SHG}, replacing $\shg$ by $SFH$.
\eprop
\bpf
Term (1) follows from the adjunction inequality in \cite[Theorem 2]{juhasz2010polytope}. Note that if $2i=|\partial S\cap \ga|-\chi(S)$, then for $\mathfrak{s}$ corresponds to $SFH(M,\ga,S,i)$, we have
\begin{equation}\label{eq: term}
    \langle c_1(\mathfrak{s},t_S),[S]\rangle=\chi(S)-|\partial S\cap \ga|=c(S,t_S),
\end{equation}where the last equality follows from (\ref{eq: defn c}) and (\ref{eq: rst}).

Term (2) follows from \cite[Lemma 3.10]{juhasz2008floer} and (\ref{eq: term}).

\begin{figure}[ht]
\centering
\begin{overpic}[width=1\textwidth]{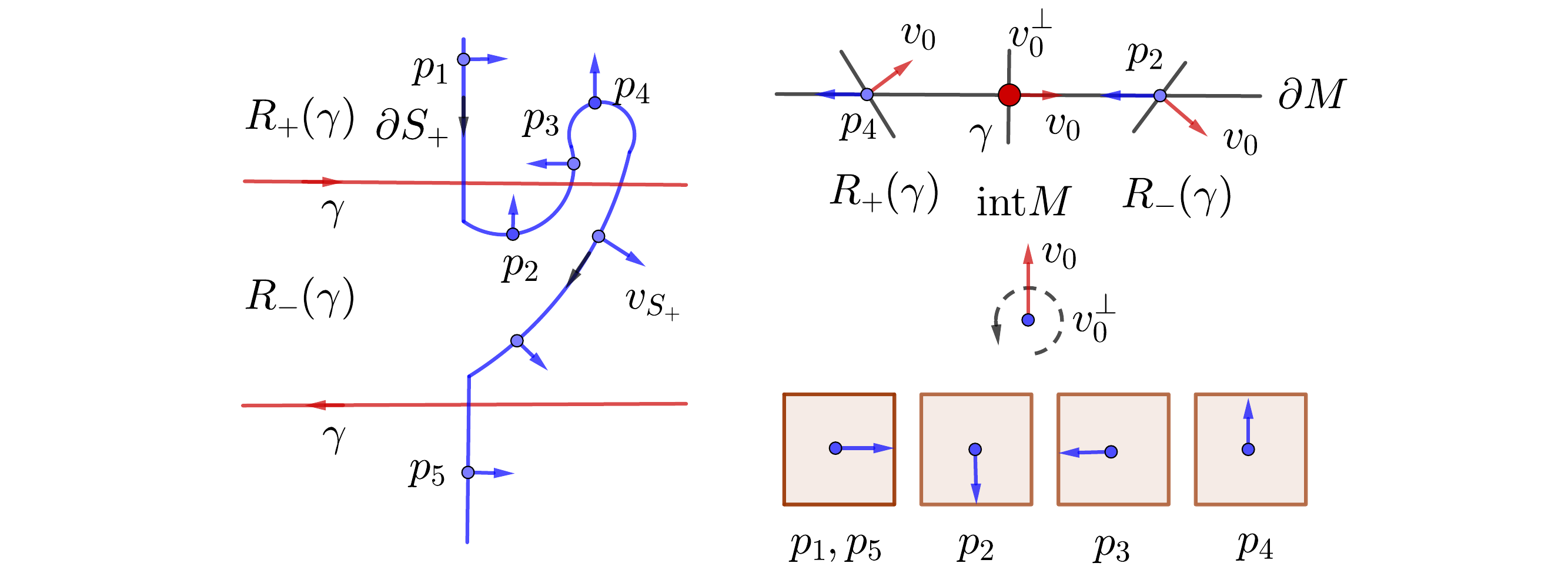}
\end{overpic}
\caption{Left subfigure, the positive stabilization $S_+$ and the vector $v_{S_+}$, where ${\rm int} M$ is inside the page. Right top subfigure, the vector field $v_0$ and the plane field $v_0^\perp$ from another viewpoint. Right middle subfigure, the orientation of $v_0^\perp$, where $v_0$ points out the page. Right bottom subfigures, the projections of $v_{S_+}$ on $v_0^\perp$ on given points.\label{fig: rotation}}
\end{figure}

Terms (3)-(5) follow from definitions and symmetry on balanced diagrams.
\epf
\bprop
Consider the stabilized surfaces $S^p$ and $S^{p+2k}$ in Theorem \ref{thm: zn grading}. Then for any $l\in\mathbb{Z}$, we have
\[SFH(M,\ga,S^p,l)=SFH(M,\ga,S^{p+2k},l+k).\]
\eprop
\bpf
Suppose $S^+$ and $S^-$ are positive and negative stabilizations of $S$. Since the stabilization operation is local, we have the following equation by direct calculation\[r(S^+,t)=r(S,t)-1\]for any trivialization $t$ of $v_0^\perp$; see Figure \ref{fig: rotation}.
Note that $[S^+]=[S].$ Hence for $\mathfrak{s}\in\spin(M,\ga)$ corresponds to $SFH(M,\ga,S,i)$, we have\beq\langle c_1(\mathfrak{s},t_{S^+}),[S^+]\rangle=&\langle c_1(\mathfrak{s},t_{S}),[S^+]\rangle+r(S^+,t_{S})\\
=&\langle c_1(\mathfrak{s},t_{S}),[S^+]\rangle+r(S,t_{S})-1\\
=&\langle c_1(\mathfrak{s},t_{S}),[S^+]\rangle-1\\
=&\langle c_1(\mathfrak{s},t_{S}),[S]\rangle-1\\
=&-2i-1.\eeq
Applying this calculation for $(2k)$ times gives the desired result.
\epf
\bprop
Suppose $S_1$ and $S_2$ are two admissible surfaces in $(M,\ga)$ such that$$[S_1]=[S_2]=\al\in H_2(M,\partial M).$$Then there exists a constant $C$ so that $$\shg(M,\ga,S_1,l)=\shg(M,\ga,S_2,l+C).$$
\eprop
\bpf
This follows directly from the definition.
\epf
\subsection{Euler characteristics}\quad

Then we can consider the Euler characteristic of $SFH$ with respect to spin$^c$ structures.
\bdefn[]\label{defn: chi(sfh)}
For a balanced sutured manifold $(M,\ga)$, let the $\mathbb{Z}_2$ grading of $SFH(M,\ga)$ be induced by the sign of intersection points of $\mathbb{T}_\al$ and $\mathbb{T}_\be$ for some compatible diagram $\mch=(\Sigma,\al,\be)$ (\textit{c.f.} \cite[Section 3.4]{Friedl2009}). Suppose $H=H_1(M)$ and choose any $\mathfrak{s}_0\in Spin^c(M,\ga)$. The \textbf{Euler characteristic} of $SFH(M,\ga)$ is $$\chi(SFH(M,\ga)) \colonequals \sum_{\substack{\mathfrak{s}\in \spin(M,\ga)\\\mathfrak{s}-\mathfrak{s}_0=h\in H^2(M,\partial M)}} \chi(SFH(M,\ga,\mathfrak{s})) \cdot {\rm PD}(h)\in \mathbb{Z}[H]/\pm H,$$where ${\rm PD}: H^2(M,\partial M)\to H_1(M)$ is the Poincar\'{e} duality map.
\edefn
\bthm[\cite{Friedl2009}]\label{thm: chi SFH}Suppose $(M,\ga)$ is a balanced sutured manifold. Then$$\chi(SFH(M,\ga))=\tau (M,\ga),$$where $\tau(M,\ga)$ is a (Turaev-type) torsion element computed from the map $$\pi_1(R_-(\ga))\to \pi_1(M)$$ by Fox calculus. In particular, if $(M,\ga)=(Y(K),\ga_K)$ for a knot $K$ in $Y$, then $$\tau(M,\ga)=(1-[m])\tau(Y(K)),$$where $m$ is the meridian of $K$ and $\tau(Y(K))$ is the Turaev torsion defined in \cite{Turaev2002}.
\ethm

\subsection{Surgery exact triangle}\quad

Suppose $(M,\ga)$ is a balanced sutured manifold and $K$ is a knot in $M$. Consider three balanced sutured manifolds $(M_i,\ga_i)$ for $i=1,2,3$ obtained from $(M,\ga)$ by Dehn surgeries along $K$. If the Dehn filling curves $\eta_1,\eta_2,\eta_3 \subset \partial (M\backslash{\rm int}\partial N(K))$ satisfy$$\eta_1\cdot \eta_2=\eta_2\cdot\eta_3=\eta_3\cdot \eta_1=-1,$$then we have the following exact triangle for twisted sutured homology from the surgery exact triangle in the closure of $(M_i,\ga_i)$
\begin{equation}\label{eq_2: Floer's triangle sutured}
\xymatrix{
\shg(M_1,\ga_1)\ar[rr]&&\shg(M_2,\ga_2)\ar[dl]\\
&\shg(M_3,\ga_3)\ar[ul]&
}	
\end{equation}

In this subsection, we show the exact triangle (\ref{eq_2: Floer's triangle sutured}) is also true when replacing $\shg$ by $SFH$.

First, we quickly review Juh\'asz's construction of the cobordism map associated to a Dehn surgery (\textit{c.f.} \cite[Section 6]{juhasz2016cobordisms}, see also \cite{Ozsvath2006} for Dehn surgeries on closed 3-manifolds).

For simplicity, suppose $\eta_1$ is the meridian of $K$. Choose an arc $a$ connecting $K$ to $R_+(\ga)$. We can construct a sutured triple diagram $(\Sigma,\al,\be,\delta)$ satisfying the following properties.

\begin{enumerate}
    \item $|\al|=|\be|=|\ga|=d$.
    \item $(\Sigma,\al,\{\be_{2},\dots,\be_d\})$ is a diagram of $(M^\p,\ga^\p)=(M\backslash N(K\cup a),\ga)$.
    \item $\delta_{2},\dots,\delta_d$ are obtained from $\be_2,\dots,\be_d$ by small isotopy, respectively.
    \item After compressing $\Sigma$ along $\be_{2},\dots,\be_d$, the induced curves $\be_1$ and $\delta_1$ lie in the punctured torus $\partial N(K)\backslash N(a)$.
    \item $\be_1$ represents  the meridian $\eta_1$ of $K$ and $\delta_1$ represents the curve $\eta_2$. In particular, $\be_1$ intersects $\delta_1$ transversely at one point.
\end{enumerate}

Then we can construct a 4-manifold $\mathcal{W}_{\al,\be,\delta}$ associated to $(\Sigma,\al,\be,\delta)$ such that it is a cobordism from $(M,\ga)=(M_1,\ga_1)$ to \[(M_2,\ga_2)\sqcup (R_+\times I\times \partial R_+\times I)\#^{d-n}(S^2\times S^1),\]where $R_+=R_+(\ga)$ and different copies of $S^2 \times S^1$ might be summed along different
components of $R_+ \times I$.

Choose a top dimensional generator $\Theta_{\be,\delta}$ of \[SFH (R_+\times I\times \partial R_+\times I)\#^{d-n}(S^2\times S^1)\cong \Lambda^* H^1 (\#^{d-n}(S^2\times S^1)). \]
Note that $(\Sigma,\al,\be)$ is a balanced diagram of $(M_1,\ga_1)$ and $(\Sigma,\al,\delta)$ is a balanced diagram of $(M_2,\ga_2)$. There is a map\[F_{\al,\be,\ga}:SFH(\Sigma,\al,\be)\otimes SFH(\Sigma,\be,\delta)\ra SFH(\Sigma,\al,\de)\]obtained by counting holomorphic triangles in $(\Sigma,\al,\be,\de)$. Then define the cobordism map as\[\begin{aligned}
F_1:SFH(M_1,\ga_1)&\ra SFH(M_2,\ga_2)\\
F_1(x)&=F_{\al,\be,\de}(x,\Theta_{\be,\de})
\end{aligned}\]Similarly, we can define the cobordism maps $F_2$ and $F_3$.

\bthm[Surgery exact triangle]\label{thm: surgery exact triangle SFH}
Consider $(M_i,\ga_i)$ and cobordism maps $F_i$ for $i=1,2,3$ as above. Then we have an exact triangle

\begin{equation}\label{eq_2: Floer's triangle SFH}
\xymatrix{
SFH(M_1,\ga_1)\ar[rr]^{F_1}&&SFH(M_2,\ga_2)\ar[dl]^{F_2}\\
&SFH(M_3,\ga_3)\ar[ul]^{F_3}&
}	
\end{equation}

\ethm
\bpf
The proof follows the proof of \cite[Theorem 9.12]{Ozsvath2004c} without essential changes (see also \cite{ozsvathbranch,lectureos06}). Since the cobordism maps $F_i$ are well-defined on $SFH$, we can verify the exact triangle for any diagram. We can construct a diagram $(\Sigma,\al,\be,\de,\zeta)$ such that $(\Sigma,\al,\be,\de)$ defines $F_1$, $(\Sigma,\al,\de,\zeta)$ defines $F_2$, and $(\Sigma,\al,\zeta,\be)$ defines $F_3$. Then we can verify the assumptions of the triangle detection lemma \cite[Lemma 4.2]{ozsvathbranch} by counting holomorphic squares and pentagons and then this lemma induces the desired exact triangle.
\epf

\subsection{Contact handles and bypasses}\label{subsection: Contact handles and bypasses SFH}\quad

Suppose $(M,\ga)\subset(M^\p,\ga^\p)$ is a proper inclusion of balanced sutured manifolds and suppose $\xi$ is a contact structure on $M^\p\backslash{\rm int} M$ with dividing sets $\ga^\p\cup(-\ga)$. Honda, Kazez, and Mati{\'{c}} \cite{honda2008contact} defined a map \[\Phi_{\xi}:SFH(M,\ga)\ra SFH(M^\p,\ga^\p),\]which is indeed the motivation of Baldwin and Sivek's construction in Subsection \ref{subsection: Contact handles and bypasses}. 

Originally, this map is defined by partial open book decompositions, and there are some technical conditions. Juh\'asz and Zemke \cite{juhasz1803contact} provided an alternative description of this map by contact handle decompositions. Their description is explicit on balanced diagrams of sutured manifolds. We will follow this alternative definition and describe the maps for contact 1- and 2-handle attachments. 

It is also worth mentioning that Zarev \cite{zarev2010joining} defined a gluing operation for sutured manifolds and conjectured the map associated to contact structures above can be recovered by the gluing operation. This was proved by Leigon and Salmoiraghi \cite{LS2020contact}.

\begin{figure}[htbp]
\centering
\begin{minipage}[t]{0.48\textwidth}
\centering
\begin{overpic}[width=5cm]{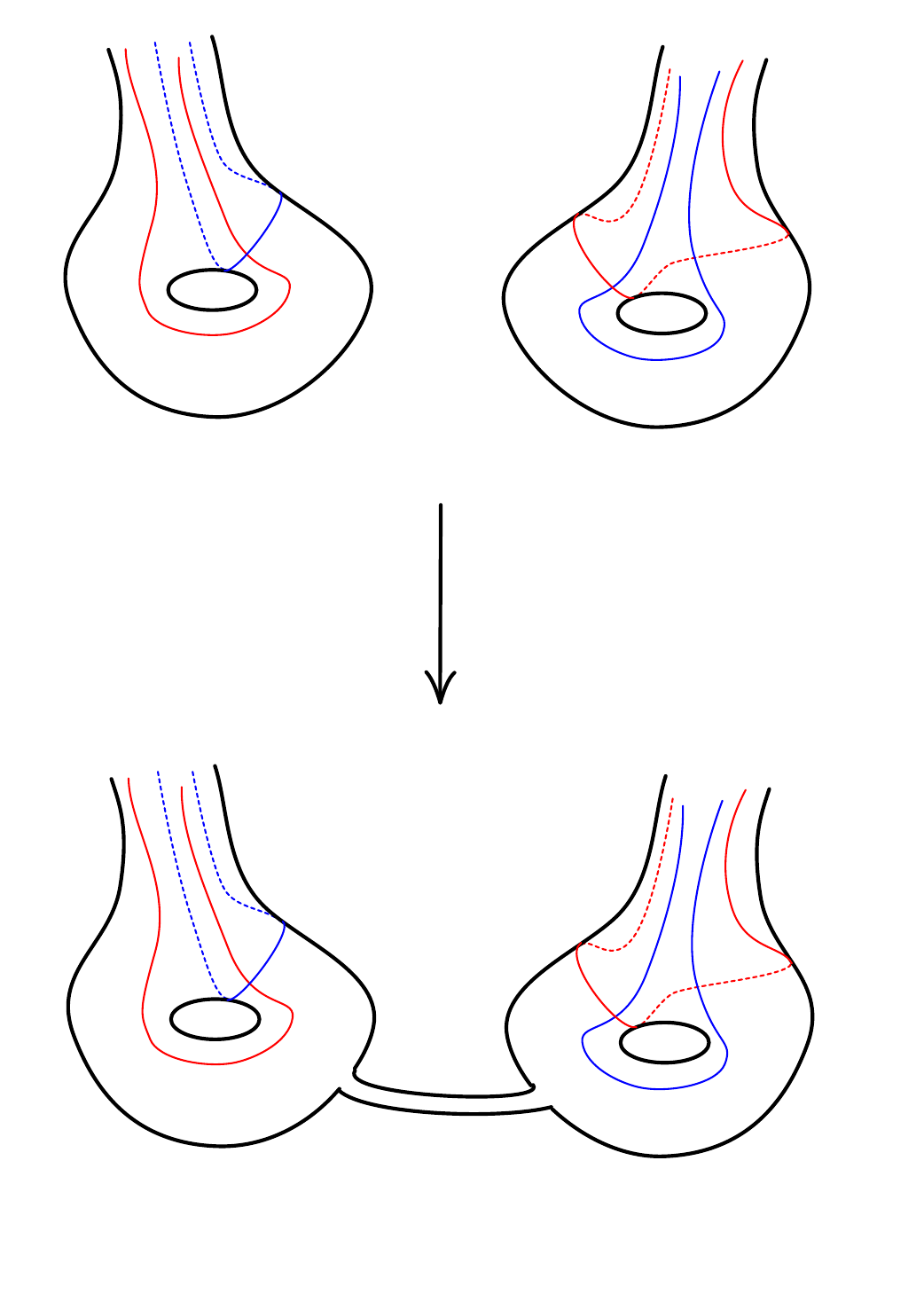}
\put(37,52){$1$-handle}
\end{overpic}
\vspace{-0.0in}
\caption{Contact 1-handle.\label{SFH 1-handle}}
\end{minipage}
\begin{minipage}[t]{0.48\textwidth}
\centering
\begin{overpic}[width=5cm]{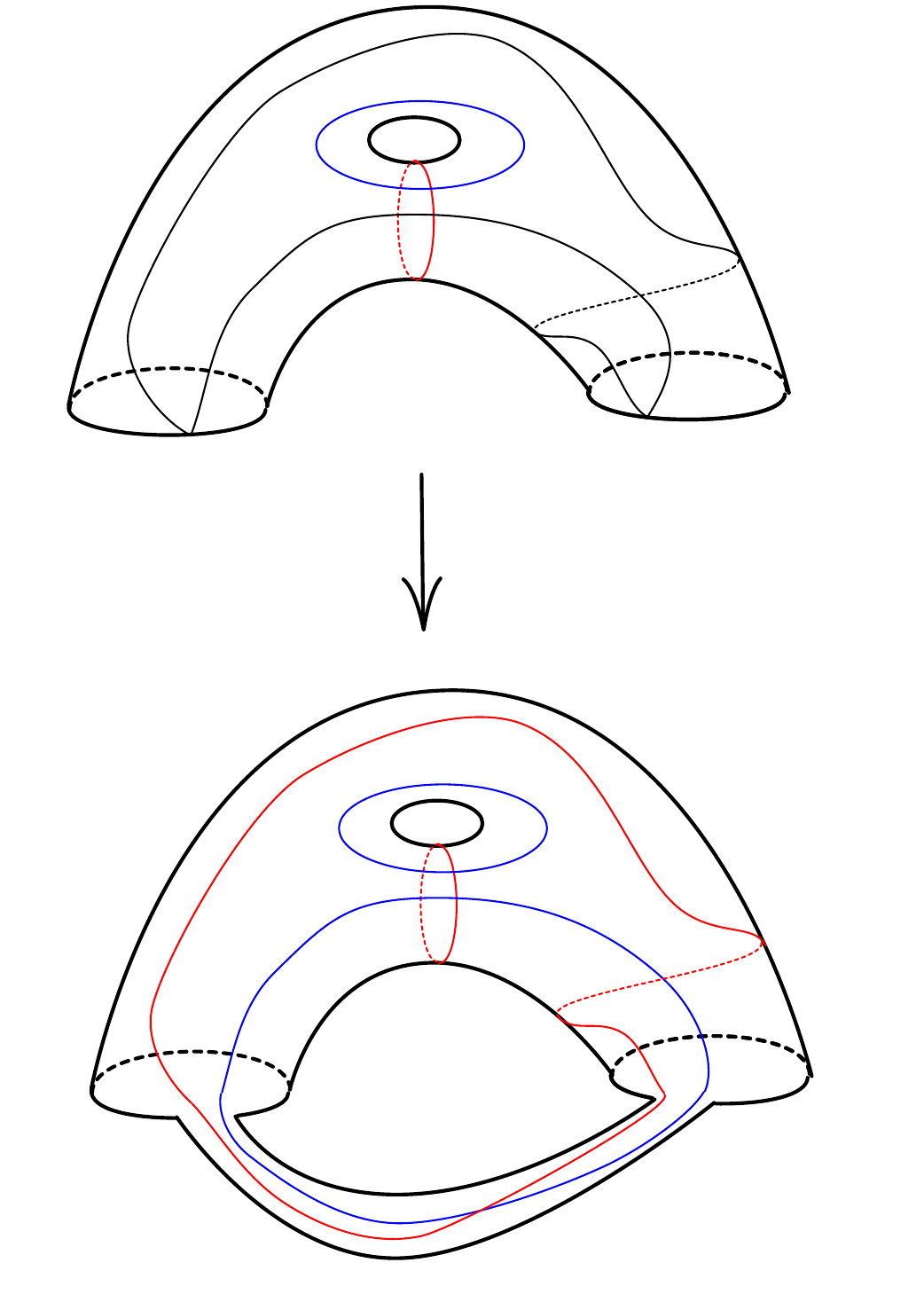}
\put(37,54){$2$-handle}
\end{overpic}
\vspace{-0.0in}
\caption{Contact 2-handle.\label{SFH 2-handle}}
\end{minipage}
\end{figure}

Juh\'asz and Zemke's construction can be shown in Figure \ref{SFH 1-handle} and Figure \ref{SFH 2-handle} (\cite[Figure 1.1]{juhasz1803contact}). Note that for all maps associated to contact structures, we should reverse the orientations of the manifold and the suture.

Let $(\Sigma,\al,\be)$ be a balanced diagram compatible with $(M,\ga)$. Then $(-\Sigma,\al,\be)$ is a balanced diagram compatible with $(-M,-\ga)$. Attaching a (3-dimensional) contact 1-handle along $D_+$ and $D_-$ corresponds to attaching a 2-dimensional 1-handle along $D_+\cap \ga$ and $D_-\cap\ga$ in $\partial \Sigma$. This operation does not change the sutured Floer chain complex and we define $C_{h^1}=C_{h^1,D_+,D_-}$ as the tautological map on intersection points.

For a contact 2-handle attachment along $\mu\subset \partial M$, note that $|\mu\cap \ga|=2$. Suppose $\lambda_+$ and $\lambda_-$ are arcs corresponding to $\mu\cap R_+(\ga)$ and $\mu\cap R_-(\ga)$, respectively. After isotopy, we can suppose $\lambda_+$ and $\lambda_-$ are propertly embedded arcs on $\Sigma$. We glue a 2-dimensional 1-handle $h$ along $\partial \Sigma$ to obtain $\Sigma^\p$, and construct two curves $\al_0$ and $\be_0$ that intersect at one point $c$ in H, and such that \[\al_0\cap\Sigma=\lambda_+,\be_0\cap\Sigma=\lambda_-.\]Consider the balanced diagram $(\Sigma^\p,\al\cup\{\al_0\},\be\cup\{\be_0\})$ and define the map associated to the contact 2-handle attachment as\[C_{h^2}(\bs{x})=C_{h^2,\mu}(\bs{x})\deq \bs{x}\times \{c\}\]for any $\bs{x}\in\mathbb{T}_\al\cap\mathbb{T}_\be$.

Since a bypass attachment can be regarded as a composition of a contact 1-handle and 2-handle attachment (\textit{c.f.} Subsection \ref{subsection: Contact handles and bypasses}), we can define the bypass map by $C_{h^2}\circ C_{h^1}$.

Honda \cite{honda2000bypass} proposed an exact triangle associated to bypass maps for $SFH$, which is indeed the motivation of the bypass exact triangle in Theorem \ref{thm_2: bypass exact triangle on general sutured manifold}. A proof of the exact triangle based on bordered sutured Floer homology was provided by Etnyre, Vela-Vick, and Zarev \cite{etnyre2017sutured}.
\bthm[{Bypass exact triangle, \cite[Section 6]{etnyre2017sutured}}]\label{thm: bypass SFH}
Suppose $(M,\ga_1)$, $(M,\ga_2)$, $(M,\ga_3)$ are balanced sutured manifolds such that the underlying 3-manifolds are the same, and the sutures $\ga_1$, $\ga_2$, and $\ga_3$ only differ in a disk shown in Figure \ref{fig: the bypass triangle}. Then there exists an exact triangle
\begin{equation*}\label{eq: bypass exact triangle, SFH}
\xymatrix@R=6ex{
SFH(-M,-\ga_1)\ar[rr]^{\psi_1}&&SFH(-M,-\ga_2)\ar[dl]^{\psi_2}\\
&SFH(-M,-\ga_3)\ar[lu]^{\psi_3}&
}    
\end{equation*}
where $\psi_1,\psi_2,\psi_3$ are bypass maps associated to the corresponding bypass arcs.
\ethm

From Juh\'asz and Zemke's description of contact gluing maps, it is obvious that the maps respect the decomposition of $SFH$ by spin$^c$ structures. We describe this fact explicitly as follows. 
\blem\label{lem: preserve spinc}
Suppose $(M,\ga)$ is a balanced sutured manifold and suppose $(M^\p,\ga^\p)$ is the resulting sutured manifold after either a contact 1-handle or 2-handle attachment. For any spin$^c$ structure $\mathfrak{s}\in\spin(-M,-\ga)$, suppose $\mathfrak{s}^\p\in\spin(-M^\p,-\ga^\p)$ is its extension corresponding to handle attachments. Then we have\[C_{h^{i}}(SFH(-M,-\ga,\mathfrak{s}))\subset SFH(-M^\p,-\ga^\p,\mathfrak{s}^\p),\]where $i\in\{1,2\}$.
\elem
\bpf
We prove the claim on the chain level. After fixing a spin$^c$ structure $\mathfrak{s}_0$ on $(M,\ga)$, we can identify $\spin(M,\ga)$ with $H^2(M,\partial M)\cong H_1(M)$. Moreover, we can represent the difference of two spin$^c$ structures by a one-cycle in Proposition \ref{prop: one cycle}. 

We can extend $\mathfrak{s}_0$ to a spin$^c$ structure $\mathfrak{s}_0^\p$ on $(M,\ga)$ and identify $\spin(M^\p,\ga^\p)$ with $H_1(M^\p)$. The inclusion $i:M\to M^\p$ induces a map \[i_*:H_1(M)\to H_1(M^\p).\]

For any $\bs{x},\bs{y}\in\mathbb{T}_\al\cap\mathbb{T}_\be$, the one cycle $\ga_{\bs{x}}-\ga_{\bs{y}}$ defined in Proposition \ref{prop: one cycle} lies in the interior of $M$.

For a contact 1-handle, since the associated map $C_{h^1}$ is tautological on intersection points, the homology class $i_*([\ga_{\bs{x}}-\ga_{\bs{y}}])$ characterizes the difference of spin$^c$ structures on $(M^\p,\ga^\p)$ for $\bs{x}$ and $\bs{y}$.

For a contact 2-handle, since $\ga_{\bs{x}\times\{c\}}$ is the union of multi-trajectory $\ga_{\bs{x}}$ and the trajectory associated to $c$, we have\[[\ga_{\bs{x}\times\{c\}}-\ga_{\bs{y}\times\{c\}}]=i_*([\ga_{\bs{x}}-\ga_{\bs{y}}]).\]This implies the desired proposition.
\epf
\brem
The reader can compare Lemma \ref{lem: preserve spinc} with Lemma \ref{lem: contact 2-handle preserve gradings}. Note that when $H_1(M)$ has torsions, preserving the spin$^c$ structures is stronger than preserving the gradings associated to an admissible surface.
\erem
\bcor\label{cor: bypass preserve grading}
Suppose $\al$ is a bypass arc on a balanced sutured manifold $(M,\ga)$. Suppose $(M,\ga^\p)$ is the resulting manifold after the bypass attachment along $\al$. Then the bypass map $\psi_\al$ for $SFH$ respects spin$^c$ structures, \textit{i.e.}, for any $\mathfrak{s}\in\spin(M,\ga)$ and its extension $\mathfrak{s}^\p\in\spin(M,\ga^\p)$, we have\[\psi_\al(SFH(-M,-\ga,\mathfrak{s}))\subset SFH(-M,-\ga^\p,\mathfrak{s}^\p).\] 
\ecor
\bpf
This follows directly from Lemma \ref{lem: preserve spinc} by the fact that a bypass attachment is a composition of a contact 1-handle and 2-handle attachment.
\epf
\brem
By Corollary \ref{cor: bypass preserve grading}, if we consider the $\mathbb{Z}$-grading associated to an admissible surface $S$ in Subsection \ref{defn: SFH grading}, then the bypass exact triangle in Theorem \ref{thm: bypass SFH} satisfies the similar grading shifting property to that in Lemma \ref{lem_3: graded bypass triangle}.
\erem

For twisted sutured homology, the map associated to a contact 2-handle is defined by the composition of the inverse of a contact 1-handle map and the cobordism map of a 0-surgery. The following proposition shows that we can define the map $C_{h^2}$ for $SFH$ in the same way.

\begin{figure}[ht]
\centering
\begin{overpic}[width=0.8\textwidth]{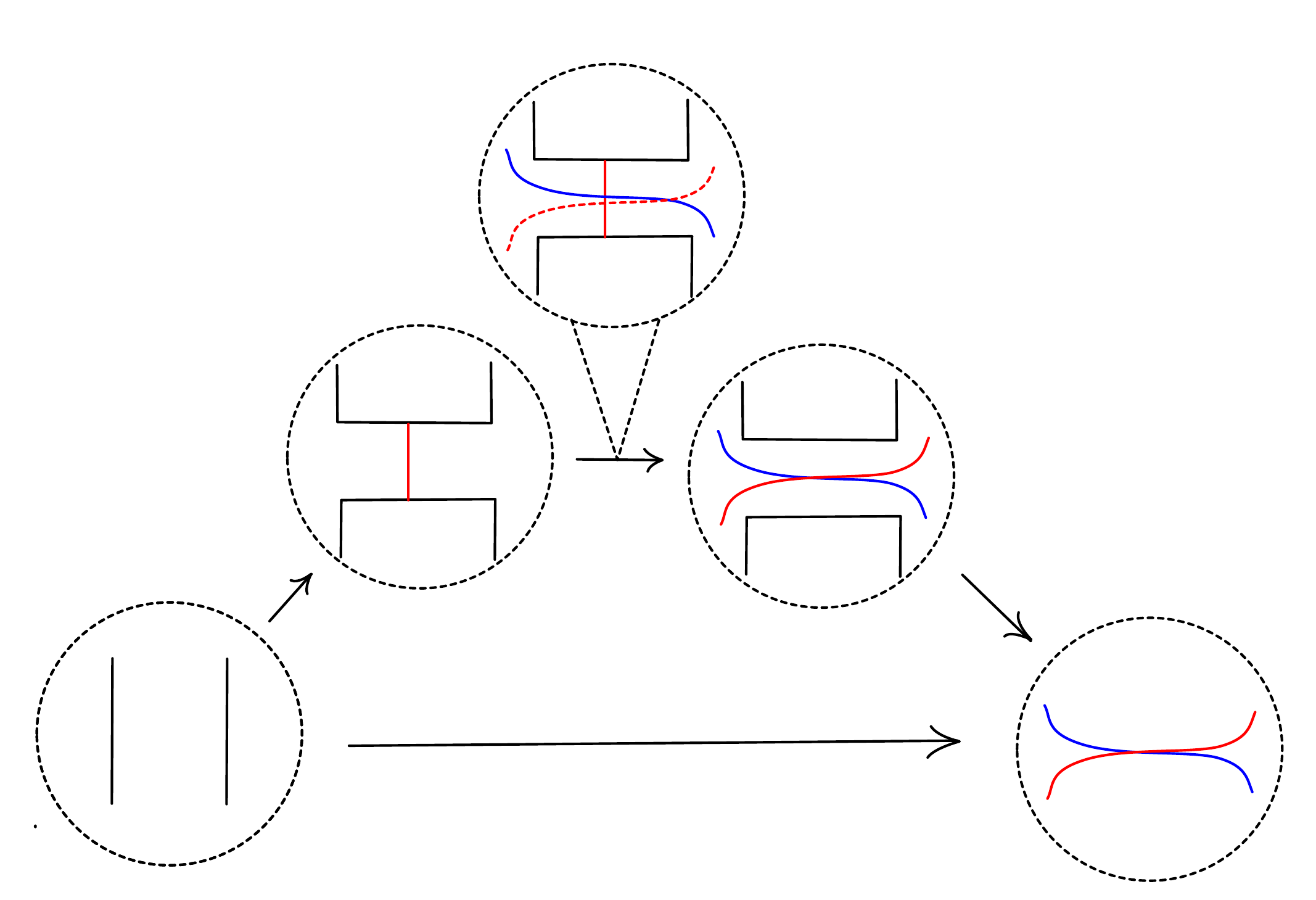}
	\put(5,12){$\Sigma$}
	\put(19,12){$\Sigma$}
\end{overpic}
\caption{Realizing the contact 2-handle map (bottom-most long arrow) as a composition of a compound stabilization (top), followed by a 4-dimensional 2-handle map (middle left), followed by a product disk map (middle right). A holomorphic triangle of the 2-handle map is indicated in the top subfigure.}
\label{GZ}
\end{figure}

\blem[\cite{GZ2021}]\label{lemma: GZ}
Suppose $(M,\ga)$ is a balanced sutured manifold and $(M^\p,\ga^\p)$ is the resulting sutured manifold after  a contact 2-handle attachment along $\mu\subset \partial M$. Let $\mu^\p$ be the framed knot obtained by pushing $\mu$ into the interior of $M$ slightly, with the framing induced from $\partial M$. Suppose $(N,\ga_N)$ is the sutured manifold obtained from $(M,\ga)$ by a 0-surgery along $\mu^\p$. Let \[F_{\mu^\p}:SFH(-M,-\ga)\ra SFH(-N,-\ga_N)\]be the associated map. Let $D\subset N$ be the product disk which is the union of the annulus bounded by $\mu\cup\mu^\p$ and the meridian disk of the filling solid torus. Let \[C_D: SFH(-N,-\ga_N)\ra SFH(-M^\p,-\ga^\p)\]be the map associated to the decomposition along $D$ (\textit{i.e.} the inverse of a contact 1-handle map). Then we have
\[C_{h^2,\mu} = C_D \circ F_{\mu^\p}:SFH(-M,-\ga)\ra SFH(-M^\p,-\ga^\p).\]
\elem
\bpf
Since all maps are well-defined on $SFH$, we can verify the claim by any diagram. Suppose $(\Sigma,\al,\be)$ is a balanced diagram compatible with $(M,\ga)$. We note that the map associated to the 0-surgery along $\mu^\p$ may be achieved by first performing a compound stabilization and then computing a triangle map. The resulting diagram leaves an extra band which is deleted by $C_D$. By \cite[Theorem 9.4]{ozsvath2004holomorphic}, the claim then follows from a model computation in the stabilization region, as shown in Figure \ref{GZ}. 
\epf

Combining the surgery exact triangle in Theorem \ref{thm: surgery exact triangle SFH} with Lemma \ref{lemma: GZ}, we obtain similar results in Lemma \ref{lem_3: surgery triangle} for $SFH$.

\bprop\label{lem_3: surgery triangle, SFH}
Consider the setups in Subsection \ref{Subsection: Basic setups}. Suppose $T'=T\backslash\al=T_2\cup\cdots\cup T_m$. Then for any $n\in\mathbb{N}$, there is an exact triangle
\begin{equation}\label{eq_6: surgery exact triangle SFH}
\xymatrix@R=6ex{
SFH(-M_{T},-{\Ga}_{n})\ar[rr]&&SFH(-M_{T},-{\Ga}_{n+1})\ar[dl]^{F_{n+1}}\\
&SFH(-M_{T'},-\ga_{T'})\ar[ul]^{G_n}&
}    
\end{equation}
The map $F_{n+1}$ is induced by the contact 2-handle attachment along the meridian of $\al$. Furthermore, we have commutative diagrams related to $\psi^{n}_{+,n+1}$ and $\psi^{n}_{-,n+1}$, respectively
\begin{equation*}
\xymatrix@R=6ex{
SFH(-M_{T},-{\Ga}_{n})\ar[rr]^{\psi_{\pm,n+1}^n}&&SFH(-M_{T},-{\Ga}_{n+1})\\
&SFH(-M_{T'},-\ga_{T'})\ar[ul]^{G_n}\ar[ur]_{G_{n+1}}&
}    
\end{equation*}
and
\begin{equation*}
\xymatrix@R=6ex{
SFH(-M_{T},-{\Ga}_{n})\ar[rr]^{\psi_{\pm,n+1}^n}\ar[dr]^{F_n}&&SFH(-M_{T},-{\Ga}_{n+1})\ar[dl]_{F_{n+1}}\\
&SFH(-M_{T'},-\ga_{T'})&
}    
\end{equation*}
\eprop
\bpf
It follows from the proof of Lemma \ref{lem_3: surgery triangle}.
\epf

\section{Proof of main theorems}\label{sec: decomposition}

In this section, we prove Theorem \ref{thm: main} and Theorem \ref{thm:shm=sfh} in the introduction.

\bthm[{Theorem \ref{thm: main}}]\label{thm: chi equal, main}
Suppose $(M,\ga)$ is a balanced sutured manifold. Suppose $H=H_1(M)$ and consider the (Turaev-type) torsion element $\tau(-M,-\ga)$ in Theorem \ref{thm: chi SFH}. Then we have\[\en(\shg(-M,-\ga))=\chi(SFH(-M,-\ga))=\tau(-M,-\ga)\in \mathbb{Z}[H]/\pm H.\]
\ethm
\bpf
By Theorem \ref{thm: chi SFH}, it suffices to prove the first equation. 

First, we consider the case that $(M,\ga)$ is strongly balanced. By discussion in Subsection \ref{subsection: sfh gradings}, we can construct a $\mathbb{Z}$-grading on $SFH$ associated to an admissible surface $S\subset (M,\ga)$. By discussion in Section \ref{section: Sutured Heegaard Floer homology}, this $\mathbb{Z}$-grading also satisfies properties in Subsection \ref{Subsection: Basic setups} for the $\mathbb{Z}$-grading on $\shg$ associated to $S$. Hence for a vertical tangle $T$ satisfies the conditions in Definition \ref{defn: decomposition}, we can define a vector space $\mathcal{SFH}_T(-M,-\ga)$ similar to $\mathcal{SHG}_T(-M,-\ga)$ in Definition \ref{defn: decomposition}. Note that for any $h_1\in H_1(M)$, the summand of $\mathcal{SFH}_T(-M,-\ga)$ is  $$SFH(-M_T,-\ga_T,\bs{S},-\bs{i}^\p)$$ with $j(\bs{\rho}^{-\bs{i}^\p})=h$, where $\bs{S},\bs{\rho},\bs{i}^\p$ come from Definition \ref{defn: enhanced euler}, the gradings come from Definition \ref{defn: SFH grading}, and $j: H_1(M_T)\to H_1(M)$ is the map induced by inclusion.

Similar to Proposition \ref{prop: isomorphism}, there is an isomorphism\begin{equation}\label{eq: sfh decomp}
    \mathcal{SFH}_T(-M,-\ga)\cong SFH(-M,-\ga).
\end{equation}Moreover, by the proofs of Lemma \ref{lem_3: surgery triangle} and Proposition \ref{prop_4: F_n is an isomorphism}, the isomorphism in (\ref{eq: sfh decomp}) is induced by contact 2-handle attachments along meridians of tangle components of $T$. Hence by Lemma \ref{lem: preserve spinc}, the isomorphism in (\ref{eq: sfh decomp}) respects spin$^c$ structures. This implies that there exists $\mathfrak{s}_0\in\spin(-M,-\ga)$, such that for any $h\in H_1(M)$, the summand of $\mathcal{SFH}_T(-M,-\ga)$ corresponding to $h$  is isomorphic to $SFH(-M,-\ga,s_0+h)$. In particular, we have
\[\en(SFH(-M,-\ga))\deq j_*(\chi(\mathcal{SFH}_T(-M,-\ga))=\chi(SFH(-M,-\ga))\in \mathbb{Z}[H]/\pm H,\]where $j_*:\mathbb{Z}[H_1(M_T)]\to \mathbb{Z}[H_1(M)]$.

By definition, the vector spaces $\mathcal{SFH}_T(-M,-\ga)$ and $\mathcal{SHG}_T(-M,-\ga)$ are direct summands of $SFH(-M_T,-\Ga)$ and $\shg(-M_T,-\Ga)$ for some $\Ga\subset \partial M_T$, respectively. By Lemma \ref{lem: torsion free}, the group $H_1(M_T)$ has no torsion. Hence by (\ref{eq: twisted euler}), we have \begin{equation*}
\begin{aligned}
     \gr(\mathcal{SHG}_T(-M,-\ga))=&\gr(\mathcal{SFH}_T(-M,-\ga))\\=&\chi(\mathcal{SFH}_T(-M,-\ga))\in\mathbb{Z}[H_1(M_T)]/\pm H_1(M_T).
\end{aligned}
\end{equation*}
Thus, we have
\begin{equation*}
\begin{aligned}\en(\shg(-M,-\ga))=&j_*(\gr(\mathcal{SHG}_T(-M,-\ga)))\\=&\en(SFH(-M,-\ga))\\=&\chi(SFH(-M,-\ga))\in\mathbb{Z}[H]/\pm H
\end{aligned}
\end{equation*}
Then we consider the case that $(M,\ga)$ is not strongly balanced. As mentioned in Remark \ref{rem: strongly balanced}. If $\partial M$ is not connected, we can construct a sutured manifold $(M^\p,\ga^\p)$ with connected boundary by attaching contact 1-handles (\textit{c.f.} \cite[Remark 3.6]{juhasz2008floer}). The product disks in $(M^\p,\ga^\p)$ corresponding to these 1-handles are admissible surfaces, and only one summand in the associated $\mathbb{Z}$-grading is nontrivial. Hence there is a canonical way to consider $\en(\shg(-M^\p,-\ga^\p))$ as an element in $\mathbb{Z}[H_1(M)]/\pm H_1(M)$. We can consider $(-M^\p,-\ga^\p)$ instead, and the above arguments about strongly balanced sutured manifolds apply to this case.
\epf

\bpf[{Proof of Theorem \ref{thm:shm=sfh}}]
We prove the theorem for $(-M,-\ga)$, and it suffices to prove the case that $(M,\ga)$ is strongly balanced. Consider the construction of $\mathcal{SFH}_T(-M,-\ga)$ in the proof of Theorem \ref{thm: chi equal, main}. We denote the monopole version of $\mathcal{SHG}_T(-M,-\ga)$ by $\mathcal{SHM}_T(-M,-\ga)$ and use it to provide a decomposition of $\shm(-M,-\ga)$.

By definition, the vector spaces $\mathcal{SFH}_T(-M,-\ga)$ and $\mathcal{SHM}_T(-M,-\ga)$ are direct summands of $SFH(-M_T,-\Ga)$ and $\shm(-M_T,-\Ga)$ for some $\Ga\subset \partial M_T$. By Lemma \ref{lem: torsion free}, the group $H_1(M_T)$ has no torsion. Hence it suffices to prove the theorem for $(-M_T,-\Ga)$.

We have\[SHM(-M_T,-\Ga)\cong \shm(-M_T,-\Ga)\cong SFH(-M_T,-\Ga)\otimes\Lambda\]with respect to the $H_1(M_T)$-grading induced by spin$^c$ structures or admissible surfaces, where $\Lambda$ is the Novikov ring over $\mathbb{Z}_2$. The first isomorphism comes from the construction of $SHM$ and $\shm$. The second isomorphism follows from \cite[Theorem 3.1]{Baldwin2020}, which essentially depends on the isomorphism between $\widecheck{HM}_\bullet$ and $HF^+$ for closed 3-manifolds (\textit{c.f.} \cite{kutluhan2010hf,taubes2010ech1,colin+}). The grading in the above isomorphism was discussed in {\cite[Corollary 3.45]{LY2021}}.

\epf

\section{Knots with small instanton knot homology}\label{sec: Knots with small instanton knot homology}

In this section, we prove detection results about $\widehat{HFK}$ and $KHI$ for null-homologous knots inside L-spaces.

\subsection{Restrictions on Euler characteristics}\quad\label{subsec: SFH}

In this subsection, we provide restrictions on Euler characteristics of $\widehat{HFK}$ and $KHI$. Since there is a canonical isomorphism between $\widehat{HFK}(Y,K)$ and $SFH(Y(K),\ga_K)$, we do not distinguish them later.

Suppose $(M,\ga)$ is a balanced sutured manifold and $H=H_1(M;\mathbb{Z})$. In Definition \ref{defn: chi(sfh)}, the Euler characteristic $\chi(SFH(M,\ga))$ has an ambiguity of $\pm H$. When $M=Y(K)$ for a knot $K$ in a rational homology sphere $Y$, we have $\partial M\cong T^2$. We can resolve the ambiguity of $\pm H$ as follows.

First, we resolve the sign ambiguity. Under the map $\mathbb{Z}[H]\to \mathbb{Z}[H_1(Y)]$ induced by inclusion, the Euler characteristic $\chi(SFH(M,\ga))$ becomes $$\pm \sum_{h\in H_1(Y)}h$$ by Theorem \ref{thm: chi SFH}. Hence we can fix the sign by choosing the one whose image is the positive one.

Then we resolve the ambiguity of $H$. We write $\spin(Y,K)$ for $\spin(Y(K),\ga_K)$. First, there is a natural identification $$SFH(Y(K),\ga_K)\xra{\cong} SFH(Y(K),-\ga_K),$$where $-\ga_K$ corresponds to $-K$, the knot obtained from $K$ by reversing the orientation. Second, since $\partial Y(K)\cong T^2$, the suture $-\ga_K$ is isotopic to the suture $\ga$, we have a map $$\iota:SFH(Y(K),\ga_K)\xra{\cong} SFH(Y(K),-\ga_K)\xra{=}SFH(Y(K),\ga).$$The square of this map is related to the basepoint moving map in \cite{sarkar15moving,zemke17moving}, though we do not use this fact. 

If there is a spin$^c$ structure $\mathfrak{s}\in \spin(Y,K)$ so that $\widehat{HFK}(Y,K,\mathfrak{s})$ that is invariant under $\iota$, then choose $\mathfrak{s}_0$ in Definition \ref{defn: chi(sfh)} so that $$PD(\mathfrak{s}-\mathfrak{t}_0)=e,$$where $e\in H$ is the identify element. If there is no summand that is invariant under $\iota$, and suppose $$\iota(\widehat{HFK}(Y,K,\mathfrak{s}))=\widehat{HFK}(Y,K,\mathfrak{s}^\p)$$ for some $\mathfrak{s},\mathfrak{s}^\p\in \spin(Y,K)$, then we define $\chi(\widehat{HFK}(Y,K))$ as an element in $(\frac{1}{2}\mathbb{Z})[H]$ so that the group elements corresponding to $\widehat{HFK}(Y,K,\mathfrak{s})$ and $\widehat{HFK}(Y,K,\mathfrak{s}^\p)$ are inverse elements. In particular, if $H\cong \mathbb{Z}$, then $(\frac{1}{2}\mathbb{Z})[H]\cong \mathbb{Z}[t^{\frac{1}{2}}, t^{-\frac{1}{2}}]$. Note that this definition is independent of the choices of $\mathfrak{s},\mathfrak{s}^\p$.

\bdefn\label{defn: canonical HFK}
Suppose $K$ is a knot in a rational homology sphere $Y$ and Suppose $H=H_1(Y(K))$. When fixing the $\mathbb{Z}_2$-grading and the spin$^c$ grading as above, the space $\widehat{HFK}(Y,K)$ is called the \textbf{canonical representative}. The corresponding spin$^c$ grading is called the \textbf{absolute Alexander grading}. For the canonical representative of $\widehat{HFK}(Y,K)$, the Euler characteristic $\chi(\widehat{HFK}(Y,K))$ is a well-defined element in $\mathbb{Z}[H]$ or $(\frac{1}{2}\mathbb{Z})[H]$. For $\mathfrak{s}\in \spin(Y)$, define $\widehat{HFK}(Y,K,[\mathfrak{s}])$ as the direct summand of all $\widehat{HFK}(Y,K,\mathfrak{s}^\p)$ with $\mathfrak{s}^\p\in \spin(Y,K)$ and $\mathfrak{s}^\p$ extends to $\mathfrak{s}$ on $Y$. Then $\chi(\widehat{HFK}(Y,K,\mathfrak{s}))$ is also a well-defined element in $\mathbb{Z}[H]$ or $(\frac{1}{2}\mathbb{Z})[H]$.
\edefn

Then we can state the main theorem of this subsection.

\bthm\label{thm: turaev}
Suppose $K_1$ and $K_2$ be two knots in a rational homology sphere $Y$ with $[K_1]=[K_2]\in H_1(Y)$. For $i=1,2$, suppose $m_i$ is the meridian of $K_i$. Then there exists an isomorphism $\phi:H_1(Y(K_1);\mathbb{Z})\cong H_1(Y(K_2);\mathbb{Z})$ so that $\phi([m_1])=[m_2]$. Using the isomorphism $\phi$, we write $H_1(Y(K_i);\mathbb{Z})$ as $H$ and write $[m_i]$ as $[m]\in H$. 

Consider the canonical representative of $\widehat{HFK}(Y,K_i)$. Then for any $\mathfrak{s}\in\spin(Y)$, there exists a Laurent polynomial $f_\mathfrak{s}(x)\in\mathbb{Z}[x,x^{-1}]$ and an element $h_{\mathfrak{s}}\in H$ such that\[\chi(\widehat{HFK}(Y,K_1,[\mathfrak{s}]))-\chi(\widehat{HFK}(Y,K_2,[\mathfrak{s}]))=([m]-1)^2f_\mathfrak{s}([m])h_{\mathfrak{s}},\]where both sides are elements in $\mathbb{Z}[H]$ or $(\frac{1}{2}\mathbb{Z})[H]$.
\ethm
Note that Theorem \ref{thm: turaev} is a generalization of \cite[Theorem 5.8]{Ye2020}. Indeed, the proof of \cite[Theorem 5.8]{Ye2020} applies without essential change. In the following, we prove generalizations of lemmas in \cite[Section 5]{Ye2020} and then sketch the proof of Theorem \ref{thm: turaev}.

\blem\label{lem: same homology}
Suppose $K_1$ and $K_2$ be two knots in a rational homology sphere $Y$ with $[K_1]=[K_2]\in H_1(Y)$. Then there exists an isomorphism $\phi:H_1(Y(K_1);\mathbb{Z})\cong H_1(Y(K_2);\mathbb{Z})$ so that $\phi([m_1])=[m_2]$.
\elem
\bpf
For a manifold $M$, let $T_1(M)$ be the torsion subgroup of $H_1(M)$ and let $B_1(M)=H_1(M)/T_1(M)$. By \cite[Theorem 3.1]{Brody1960}, since $[K_1]=[K_2]\in H_1(Y)$, there is a subgroup $B$ of $B_1(Y(K_1\cup K_2)$ so that for $i=1,2$, the map$$j_i:B_1(Y(K_1\cup K_2))\to B_1(Y(K_i))$$induced by the injection is an isomorphism on $B$. Moreover, the map $$k_i:T_1(Y(K_1\cup K_2))\to T_1(Y(K_i))$$induced by injection is an isomorphism. Since $H_1(M)\cong B_1(M)\oplus T_1(M)$ for any manifold $M$, we have an isomorphism $$\phi_0=(j_2\circ j_1^{-1},k_2\circ k_1^{-1}): H_1(Y(K_1))\cong B_1(Y(K_1))\oplus T_1(Y(K_1))\cong B_1(Y(K_2))\oplus T_1(Y(K_2)) \cong  H_1(Y(K_2)).$$Moreover, if $$l_i:H_1(Y(K_i))\to H_1(Y)$$is the map induced by injection, then $l_1=\phi_0\circ l_2$. 

For $i=1,2$, consider the long exact sequence about the pair $(Y,N(K_i))$:
$$H^1(Y)\to H^1(N(K_i))\to H^2(Y,N(K_i))\to H^2(Y)\to H^2(K_i)=0.$$Since $Y$ is a rational homology sphere, $$H^1(Y)\cong {\rm Hom}_\mathbb{Z}(H_1(Y),\mathbb{Z})=0.$$By excision theorem and the Poincar\'{e} duality, we have $$H^2(Y,N(K_i))\cong H^2(Y(K_i),\partial Y(K_i))\cong H_1(Y(K_i))\aand H^2(Y)\cong H_1(Y).$$Under the Poincar\'{e} duality, the image of $$H^1(N(K_i))\cong H_1(N(K_i),\partial N(K_i))\cong\mathbb{Z}$$ in $H_1(Y(K_i))$ is $[m_i]$. Since $l_1=\phi\circ l_2$, we have the following commutative diagram 

\begin{equation*}
    \xymatrix@R=6ex{
    0\ar[r]&\mathbb{Z}\ar[d]^{\cong}\ar[r]&H_1(Y(K_1))\ar[d]^{\phi_0}\ar[r]^{l_1} &H_1(Y)\ar[d]^{=}\ar[r]&0\\0\ar[r]&\mathbb{Z}\ar[r]& H_1(Y(K_2))\ar[r]^{l_2} &H_1(Y)\ar[r]&0
    }
\end{equation*}
Hence $\phi_0([m_1])=[m_2]^{\pm 1}$. If $\phi_0([m_1])=[m_2]$, let $\phi=\phi_0$. If $\phi_0([m_1])=[m_2]^{-1}$, let $\phi=\phi_0\circ \epsilon$, where $\epsilon$ maps an element to its inverse. Then $\phi:H_1(Y(K_1))\to H_1(Y(K_2))$ is an isormophism and $\phi([m_1])=([m_2])$.
\epf

\blem\label{lem: kernel 1-g}
Suppose $G$ is an abelian group and $g_0$ is an element in $G$. The quotient map $G\to G/(g_0)$ induces a map on the group ring$$\operatorname{pr}:\mathbb{Z}[G]\to \mathbb{Z}[G/(g_0)].$$Then the kernel of $\operatorname{pr}$ is generated by $1-g_0$.
\elem
\bpf
We can regard element in $\mathbb{Z}[G]$ as a function $f:G\to \mathbb{Z}$ that maps $g\in G$ to the coefficient of $g$. Note that $f(g)\neq 0$ for finitely many $g$. If $f\in ker(\operatorname{pr})$, then $$h(g)=\sum_{k\in\mathbb{N}}f(gg_0^{-k})$$is also a function $G\to \mathbb{Z}$ that is nonvanishing for finitely many $g$. It is straightforward to check that $h=(1-g_0)f$ as elements in $\mathbb{Z}[G]$.
\epf

\bpf[Proof of Theorem \ref{thm: turaev}]
The first part of this theorem is just Lemma \ref{lem: same homology}. Write $H=H_1(Y(K_i))$ and $[m]=[m_i]$. By Theorem \ref{thm: chi SFH}, we have $$\chi(\widehat{HFK}(Y,K_i))=(1-[m])\tau(K_i) \text{ in } \mathbb{Z}[H] \text{ or } (\frac{1}{2}\mathbb{Z})[H].$$\textit{A priori}, the Turaev torsion $\tau(Y(K_i))$ is not in $\mathbb{Z}[H]$, but the difference $\tau(Y(K_1))-\tau(Y(K_2))$ is. By Lemma \ref{lem: kernel 1-g}, we can apply the proof of \cite[Lemma 5.5]{Ye2020} to the case where $Y$ is a rational homology sphere. Note that we use the fact $b_1(Y(K_i))=1$ in that proof. From \cite[Lemma 5.5]{Ye2020}, we have $$\tau(Y(K_1))-\tau(Y(K_2))=(1-[m])g \text{ in } \mathbb{Z}[H] \text{ or } (\frac{1}{2}\mathbb{Z})[H] \text{ for some }g\in\mathbb{Z}[H].$$The ambiguity of $\pm H$ in the statement of \cite[Lemma 5.5]{Ye2020} is resolved because we consider absolute Alexander gradings on $\widehat{HFK}(Y,K_i)$. Then we have\begin{equation}\label{eq: 1-m2}
    \chi(\widehat{HFK}(Y,K_1))-\chi(\widehat{HFK}(Y,K_2))=(1-[m])(\tau(Y(K_1))-\tau(Y(K_2)))=(1-[m])^2g.
\end{equation}Suppose $H_1(Y;\mathbb{Z})=\{s_1,\dots,s_p\}$. Then the element $g$ can be written as the sum $$g=\sum_{j=1}^p g_j,$$where $g_j$ contains terms that are in the preimage of $s_j\in H_1(Y;\mathbb{Z})$ under the map $$H\to H/m\cong H_1(Y;\mathbb{Z}).$$For any $j$, there exists a Laurent polynomial $f_j(x)$ and an element $\tilde{s}_j\in H$ such that $g_j=f_j([m])\tilde{s}_j$. Since $\spin(Y)$ is an affine space on $H_1(Y)$, Thus, Equation (\ref{eq: 1-m2}) can be decomposed with respect to $\spin(Y)$, where $h_\mathfrak{s}$ corresponds to some $\tilde{s}_j$.

\epf

Finally, we deal with instanton knot homology. 

\bdefn\label{defn: canoncial, KHI}
Suppose $K$ is a knot in a rational homology sphere $Y$ and suppose $H=H_1(Y(K))$. Similar to the way for $\widehat{HFK}(Y,K)$, we can fix the $\mathbb{Z}_2$-grading and the $H$-grading from the enhanced Euler characteristic on $KHI(Y,K)$. When gradings are fixed, the space $KHI(Y,K)$ is called also the \textbf{canonical representative} and the corresponding $H$-grading is also called the \textbf{absolute Alexander grading}. For any element $s\in H_1(Y)$, let $[s]\subset H_1(Y(K))$ be the set of preimages of $s$ under the map $H\to H_1(Y)$ and define$$KHI(Y,K,[s])\deq\bigoplus_{h\in[s]}KHI(Y,K,h).$$Then $\chi(KHI(Y,K,[s])$ is also a well-defined element in $\mathbb{Z}[H]$ or $(\frac{1}{2}\mathbb{Z})[H]$.

\edefn

By Theorem \ref{thm: turaev} and Equation (\ref{eq: enhanced same}), we have the following corollary.

\bcor\label{cor: KHI turaev}

Suppose $K_1$ and $K_2$ be two knots in a rational homology sphere $Y$ with $[K_1]=[K_2]\in H_1(Y)$. For $i=1,2$, suppose $m_i$ is the meridian of $K_i$. Using the isomorphism $\phi$ in Lemma \ref{lem: same homology}, we write $H_1(Y(K_i);\mathbb{Z})$ as $H$ and write $[m_i]$ as $[m]\in H$. 

Consider the canonical representative of $KHI(Y,K_i)$. Then for any $s\in H_1(Y)$, there exists a Laurent polynomial $f_s(x)\in\mathbb{Z}[x,x^{-1}]$ and an element $h_{s}\in H$ such that\[\chi(KHI(Y,K_1,[s]))-\chi(KHI(Y,K_2,[s]))=([m]-1)^2f_s([m])h_{s},\]where both sides are elements in $\mathbb{Z}[H]$ or $(\frac{1}{2}\mathbb{Z})[H]$.

\ecor

\subsection{Detection results}\quad

In this subsection, we use Theorem \ref{thm: turaev} and Corollary \ref{cor: KHI turaev} to prove detection results in the introduction. 

\begin{conv}
Throughout this subsection, we suppose $K$ is a knot in a rational homology sphere $Y$ and suppose $H=H_1(Y(K))$. Moreover, we consider canonical representatives of $\widehat{HFK}(Y,K)$ and $KHI(Y,K)$ as in Definition \ref{defn: canonical HFK} and Definition \ref{defn: canoncial, KHI}. For simplicity, we write also write $(\frac{1}{2}\mathbb{Z})[H]$ as $\mathbb{Z}[H]$.
\end{conv}

First, we prove some lemmas.
\blem[]\label{lem: homology}
Suppose $K\subset Y$ is a knot so that $[K]=0\in H_1(Y)$. Then there exists a canonical isomorphism$$H_1(Y(K))\cong \mathbb{Z}\oplus H_1(Y),$$where the meridian of $K$ represents the generator of $\mathbb{Z}$.
\elem
\bpf
The isomorphism is induced by pairing with a Seifert surface of $K$.
\epf
Hence we can write elements in $H_1(Y(K))$ as $[m]^n\cdot s$ for $s\in H_1(Y)$ and $n\in\mathbb{Z}$.
\begin{lem}\label{lem: Euler characteristics of the unknot}
Suppose $Y$ is an instanton L-space and $U\subset Y$ is the unknot. Then for any $s\in H_1(Y)$, we have
$$KHI(Y,U,[s])\cong \mathbb{C}\aand \chi(KHI(Y,U,[s]))=s\in \mathbb{Z}[H].$$
\end{lem}
\bpf
The result $KHI(Y,U,[s])\cong \mathbb{C}$ follows directly from Equation (\ref{eq: enhanced same}) and the following isomorphisms:
$$KHI(Y,U)\cong I^{\sharp}(Y)\cong \mathbb{C}^{|H_1(Y)|}.$$The result $\chi(KHI(Y,U,[s]))=1$ follows from the fact that $g(U)=0$ and $KHI$ detects the genus of the knot \cite[Proposition 7.16]{kronheimer2010knots}.
\epf

\blem\label{lem: turaev, instanton}
Suppose $Y$ is an instanton L-space, and $K\subset Y$ is a knot so that $[K]=0\in H_1(Y).$ Suppose $m$ is the meridian of $K$. Then for any $s\in H_1(Y)$, the element
$$\chi(KHI(Y,K,[s]))-s\in\mathbb{Z}[H]$$
has a factor $([m]-1)^2$. See Definition \ref{defn: canoncial, KHI} for the definition of $\chi(KHI(Y,K,[s]))$.
\elem
\bpf
Since the unknot is also null-homologous, this lemma follows directly from Corollary \ref{cor: KHI turaev} and Lemma \ref{lem: Euler characteristics of the unknot}.
\epf
Then we prove the detect results in the introduction.
\bpf[Proof of Theorem \ref{thm: unknot detection}]
It is clear that
$$\dim_\mathbb{C}KHI(Y,K)=\dim_\mathbb{C}I^{\sharp}(Y)$$
when $K$ is the unknot in $Y$. Now suppose
$$\dim_\mathbb{C}KHI(Y,K)=\dim_\mathbb{C}I^{\sharp}(Y)$$
and we will show that $K$ must be the unknot. For any $s\in H_1(Y)$, Lemma \ref{lem: turaev, instanton} implies that
$$\chi(KHI(Y,K,[s]))\neq 0$$
and hence $KHI(Y,K,[s])\neq 0$. From the assumption, we have
$$\dim_\mathbb{C}KHI(Y,K)=\dim_\mathbb{C} I^{\sharp}(Y)=|H_1(Y)|.$$
Thus, we must have
$$KHI(Y,K,[s])\cong \mathbb{C}\aand \chi(KHI(Y,K,[s]))=[m]^n\cdot s\in\mathbb{Z}[H],$$where $m$ is the meridian of $K$ and $n$ is some integer. Applying Lemma \ref{lem: turaev, instanton} again, we know that $n$ must be $0$. Since $KHI$ detects the genus of the knot \cite[Proposition 7.16]{kronheimer2010knots}, we know that $g(K)=0$, which implies $K$ is the unknot.
\epf

\bpf[Proof of Theorem \ref{thm: trefoil detection}]
Applying Lemma \ref{lem: turaev, instanton}, for any $s\in H_1(Y)$, we have
$$\chi(KHI(Y,K,[s]))\neq 0.$$
Since $\chi(KHI(Y,K,[s]))$ and $\dim_{\mathbb{C}}KHI(Y,K,[s])$ have the same parity, we conclude that there exists $s_0\in H_1(Y)$ so that
$$KHI(Y,K,[s_0])\cong \mathbb{C}^3$$
and for any $s\neq s_0$,
$$KHI(Y,K,[s])\cong \mathbb{C}$$
Applying Lemma \ref{lem: turaev, instanton} again, for any $s\neq s_0$, we know that
$$\chi(KHI(Y,K,[s]))=s\in\mathbb{Z}[H].$$
For $s_0$, we know that
$$\chi(KHI(Y,K,[s_0]))-s_0\in\mathbb{Z}[H]$$
has a factor $([m]-1)^2$ and
$$\norm{\chi(KHI(Y,K,[s_0]))}\leq 3,$$where $m$ is the meridian of $K$ and the norm is defined before Example \ref{ex: not sharp}.

Hence there are only two possibilities.

{\bf Case 1}. $\norm{\chi(KHI(Y,K,[s_0]))}=1$ and we then conclude that
$$\chi(KHI(Y,K,[s_0]))=[m]^n\cdot s_0\in\mathbb{Z}[H].$$
Note that $\norm{\chi(KHI(Y,K,[s_0]))}=1$ implies that there is a $2$-dimensional summand of $KHI(Y,K)$ whose Euler characteristic is zero. Hence there are further two cases:


{\bf Case 1.1} $KHI(Y,K,[s_0])$ is supported in two different Alexander gradings. By assumption, we know that $KHI(Y,K,[s_0])$ has a $1$-dimensional summand at the Alexander grading $n$ and has a $2$-dimensional summand at the Alexander grading $n^\p$ for some $n^\p\neq n$. This contradicts the fact that $KHI(Y,K)$ is symmetric with respect to the Alexander grading.

{\bf Case 1.2}. $KHI(Y,K,[s_0])$ is support solely in the Alexander grading $n$. By symmetry of $KHI(Y,K)$, we must have $n=0$. Since $KHI$ detects the genus of the knot, we know $K$ is an knot, which contradicts Theorem \ref{thm: unknot detection} since $\dim_\mathbb{C}KHI(Y,K)=\dim_\mathbb{C}I^\sharp(Y)+2$.

{\bf Case 2}. $\norm{\chi(KHI(Y,K,[s_0]))}=3$. By symmetry on $KHI(Y,K)$, there exists $n\in\mathbb{N}_+$so that
$$\chi(KHI(Y,K,[s_0]))=([m]^{n}-1+[m]^{-n})\cdot s_0\in\mathbb{Z}[H].$$by \cite[Proposition 7.16 and Corollary 7.19]{kronheimer2010knots}, we know that $K$ is fibred of genus $n$. By \cite[Theorem 1.7]{baldwin2018khovanov}, we know that $n=1$. Hence $K$ is a genus-one-fibred knot.

\epf

The proofs of the following theorems are similar to those of Theorem \ref{thm: unknot detection} and Theorem \ref{thm: trefoil detection}, but using Theorem \ref{thm: turaev} and \cite[Theorem 1.1]{baldwin18fibred} instead of Corollary \ref{cor: KHI turaev} and \cite[Theorem 1.7]{baldwin2018khovanov}.

\bthm\label{thm: unknot detection, HFK}
Suppose $K$ is a null-homologous knot in a rational homology sphere $Y$. If \begin{equation}\label{eq: Lspace, HFK}\dim_\ft \widehat{HF}(Y)=|H_1(Y;\mathbb{Z})|,\end{equation}then $K$ is the unknot if and only if
\begin{equation}\label{eq: floer simple knot}
    \dim_\ft \widehat{HFK}(Y,K)=\dim_\ft \widehat{HF}(Y).
\end{equation}
\ethm
\bthm\label{thm: trefoil detection, HFK}
Suppose $K$ is a null-homologous knot in a rational homology sphere $Y$. If
\begin{equation}\label{eq: trefoil-like knot, HFK}\dim_\ft \widehat{HFK}(Y,K)=\dim_\ft \widehat{HF}(Y)+2=|H_1(Y;\mathbb{Z})|+2,\end{equation}
then $K$ must be a genus-one-fibred knot.
\ethm
\brem
Baldwin \cite{Baldwin08Lspace} classified L-spaces that contain null-homologous genus-one-fibred knots. He also computed knot Floer homologies of such knots, which only depend on their Alexander polynomials. The techniques in his classification involve the minus chain complex \cite{ozsvath2004holomorphic} and the mapping cone formula \cite{Ozsvath2008,Ozsvath2011}, which is not available in instanton theory yet. For knots in lens spaces, there are more results about genus-one-fibred knots \cite{Morimoto89,Baker09tunnel,Baker14gof}.
\erem

\bibliographystyle{alpha}

\input{main-torsion-v4-lastupdated.bbl}
\end{document}

%% file: command.tex
\usepackage{amsmath}
\usepackage{amssymb}
\usepackage{extarrows}
\usepackage{mathabx}
\usepackage[all]{xy}
\usepackage{amsbsy}
\usepackage{graphicx}
\usepackage{appendix}
\usepackage{float}
\usepackage{subfigure}
\usepackage[numbers,sort&compress]{natbib}
\usepackage{graphicx}
\usepackage{amsthm}
\usepackage{geometry}
\usepackage{fancyhdr}
\usepackage{color} 
\usepackage{overpic}
\usepackage{mathabx}
\usepackage{tikz-cd}
\usepackage{enumerate}
\usepackage{mathrsfs}
\usepackage{colonequals}
\usepackage{microtype}
\usepackage{commath}

\newtheorem{thm}{Theorem}[section]
\newtheorem{cor}[thm]{Corollary}
\newtheorem{prop}[thm]{Proposition}
\newtheorem{lem}[thm]{Lemma}

\theoremstyle{definition}
\newtheorem{defn}[thm]{Definition}
\newtheorem{exmp}[thm]{Example}

\newtheorem*{conv}{Convention}

\theoremstyle{remark}
\newtheorem{rem}[thm]{Remark}

\numberwithin{equation}{section}

\newcommand{\beq}{\begin{equation*}\begin{aligned}}
\newcommand{\eeq}{\end{aligned}\end{equation*}}
\newcommand{\bpf}{\begin{proof}}
\newcommand{\epf}{\end{proof}}
\newcommand{\bthm}{\begin{thm}}
\newcommand{\ethm}{\end{thm}}
\newcommand{\bprop}{\begin{prop}}
\newcommand{\eprop}{\end{prop}}
\newcommand{\bcor}{\begin{cor}}
\newcommand{\ecor}{\end{cor}}
\newcommand{\blem}{\begin{lem}}
\newcommand{\elem}{\end{lem}}
\newcommand{\bdefn}{\begin{defn}}
\newcommand{\edefn}{\end{defn}}
\newcommand{\bexmp}{\begin{exmp}}
\newcommand{\eexmp}{\end{exmp}}
\newcommand{\brem}{\begin{rem}}
\newcommand{\erem}{\end{rem}}
\newcommand{\benu}{\begin{enumerate}[(1)]}
\newcommand{\eenu}{\end{enumerate}}

\newcommand{\bdia}{\begin{displaymath}\xymatrix}
\newcommand{\edia}{\end{displaymath}}

\newcommand{\shm}{\underline{\rm SHM}}
\newcommand{\shi}{\underline{\rm SHI}}

\newcommand{\shib}{{\mathbf{SHI}}}
\newcommand{\shmb}{{\mathbf{SHM}}}

\newcommand{\shg}{\underline{\rm SHG}}

\newcommand{\khg}{\underline{\rm KHG}}

\newcommand{\hhat}{\underline{\rm HG}}

\newcommand{\deq}{\colonequals}
\newcommand{\sym}{{\rm Sym}}

\newcommand{\spin}{{\rm Spin}^c}

\newcommand{\al}{\alpha}
\newcommand{\be}{\beta}
\newcommand{\ga}{\gamma}
\newcommand{\Ga}{\Gamma}

\newcommand{\de}{\delta}

\newcommand{\p}{\prime}

\newcommand{\bs}{\boldsymbol}
\newcommand{\gr}{\chi_{\rm gr}}
\newcommand{\tor}{{\rm Tors}}
\newcommand{\en}{\chi_{\rm en}}
\newcommand{\aand}{~{\rm and}~}

\newcommand{\mch}{\mathcal{H}}

\newcommand{\intg}{\mathbb{Z}}
\newcommand{\ft}{{\mathbb{F}_2}}

\newcommand{\ra}{\rightarrow}

\newcommand{\xra}{\xrightarrow}

